\documentclass[11pt]{report}
\usepackage{latexsym}
\usepackage{amsmath,amssymb,amsfonts}
\usepackage{makeidx,pdfsync}
 
\newtheorem{Theorem}{Theorem}[part]
\newtheorem{Definition}{Definition}[part]
\newtheorem{Proposition}{Proposition}[part]

\newtheorem{Lemma}{Lemma}[part]
\newtheorem{Corollary}{Corollary}[part]
\newtheorem{Remark}{Remark}[part]
\newtheorem{Example}{Example}[part]

\def\txh{{\hat t,\hat x}}
\def\th{{\hat t}}
\def\xh{{\hat x}}
\def\txn{{t_n,x_n}}
\def\tn{{t_n}}

\def\we{}
\def\demi{\frac12}
\def\Dom{[0,T]\x \R^d\x \R\x\R^d\x\mathbb{S}^d}

\def\eqref#1{\reff{#1}} 
\def\be{\begin{eqnarray}}
\def\ee{\end{eqnarray}}
\def\bal{\begin{aligned}}
\def\eal{\end{aligned}}
\def\beq{\begin{equation}}
\def\eeq{\end{equation}}
\def\beqq{\begin{equation*}}
\def\eeqq{\end{equation*}}
\def\demi{\frac{1}{2}}

\def\we{\widetilde}

\def \Prod{\displaystyle\prod}

\def \be{\begin{eqnarray}}
\def \ee{\end{eqnarray}}
\def \b*{\begin{eqnarray*}}
\def \e*{\end{eqnarray*}}

\def \E{\mathbb{E}}
\def \F{\mathbb{F}}

\def \M{\mathbb{M}}
\def \N{\mathbb{N}}
\def \P{\mathbb{P}}
\def \Q{\mathbb{Q}}
\def \R{\mathbb{R}}

\def\p{+}

\def \[{[\,\!\![}
\def \]{]\,\!\!]}

\def \1{{\bf 1}}
\def \esssup{{\rm esssup}}

\def \proof{{\noindent \bf Proof. }}
\def \ep{\hbox{ }\hfill$\Box$}

\def\reff#1{{\rm(\ref{#1})}}

\def\Ac{{\cal A}}

\def\Cc{{\cal C}}

\def\Ec{{\cal E}}
\def\Fc{{\cal F}}
\def\Gc{{\cal G}}
\def\Hc{{\cal H}}

\def\Lc{{\cal L}}
\def\Mc{{\cal M}}

\def\Lc{{\cal L}}
\def\Oc{{\cal O}}
\def\Pc{{\cal P}}

\def\Tc{{\cal T}}
\def\Uc{{\cal U}}

\def\vp{\varphi}

\def\eps{\varepsilon}
 
 \def\vs#1{\vspace{#1mm}}
 \def\vv{\vs2}
\def \Prod{\displaystyle\prod}

\def\esssup{{\rm ess}\!\sup\limits}
\def \E{\mathbb{E}}
\def \F{\mathbb{F}}

\def \M{\mathbb{M}}
\def \N{\mathbb{N}}
\def \R{\mathbb{R}}

\def\P{\mathbb{P}}
\def\Q{\mathbb{Q}}

\def\Ac{{\cal A}}

\def\Cc{{\cal C}}

\def\Ec{{\cal E}}
\def\Fc{{\cal F}}
\def\Gc{{\cal G}}
\def\Hc{{\cal H}}

\def\Lc{{\cal L}}
\def\Mc{{\cal M}}

\def\Oc{{\cal O}}
\def\Pc{{\cal P}}

\def\Tc{{\cal T}}
\def\Uc{{\cal U}}

\def\ep{\hbox{ }\hfill$\Box$}
\def\reff#1{{\rm(\ref{#1})}}
\def\be{\begin{eqnarray}}
\def\ee{\end{eqnarray}}
\def\beq{\begin{equation}}
\def\eeq{\end{equation}}
\def\b*{\begin{eqnarray*}}
\def\e*{\end{eqnarray*}}
\def\x{\times}
\def \sge{\;\ge\;}

\def\={\;=\;}
\def\lra{\longrightarrow}

\def\proof{{\noindent \bf Proof. }}

\def\And{\;\mbox{ and }\;}

\def\Pas{\mathbb{P}-\mbox{a.s.}}

\def\.{\;.}

\def\eps{\varepsilon}
\def\vp{\varphi}
\def\1{{\bf 1}}
\def\Esp#1{\mathbb{E}\left[#1\right]}
\def\EspQ#1{\mathbb{E}^\mathbb{Q}\left[#1\right]}

\def\Pro#1{\mathbb{P}\left[{#1}\right]}

\def\eqref#1{\reff{#1}}

\makeatletter \@addtoreset{Theorem}{chapter}

\title{{\huge Portfolio management under risk contraints}\\~~ \\ Lectures given at MITACS-PIMS-UBC Summer School in Risk Management and
Risk Sharing}

\author{\large{Bruno Bouchard}
											  \\\small Universit\'e Paris-Dauphine-CEREMADE and ENSAE-CREST
											  \\~~~
\\											  Exercices prepared by\\											  \large{Ludovic Moreau} and 				  \large{Adrien Nguyen Huu}
											   \\\small Universit\'e Paris-Dauphine-CEREMADE}
\date{ This version : July 2010 }

\begin{document}
 
\maketitle

\tableofcontents

\chapter{Introduction et notations}

The aim of these lectures at MITACS-PIMS-UBC Summer School in Risk Management and Risk Sharing is to discuss risk controlled approaches for the pricing and hedging of financial risks. 

We will start with the classical dual approach for financial markets, which allows to rewrite super-hedging problems in terms of optimal control problems in standard form. Based on this, we shall then consider hedging and pricing problems under utility or risk minimization criteria. This approach will turn out to be powerful whenever linear (or essentially linear) problems are considered, but not adapted to more general settings with non-linear dynamics (e.g. large investor models, high frequency trading with market impact features, mixed finance/insurance issues).

In the second part of this lecture, we will develop on a new approach for risk control problems based on a stochastic target formulation. We will see how flexible this approach is and how it allows to characterize very easily super-hedging prices in term of suitable Hamilton-Jacobi-Bellman type partial differential equations (PDEs). We will then see how quantile hedging and expected loss pricing problems can be embeded into this framework, for a very large class of financial models. We shall finally consider a simple example of optimal book liquidation in which the control is a continuous non-decreasing process, as an illustration of possible practical developments in optimal trading under risk constraint.

These lectures are organized in small chapters, each of them being focused on a particular aspect. 

\section{Notations}

We first make precise some notations that will be used in all these notes. 

\vs2

In all these lectures notes, we shall consider a probability space $(\Omega,\Fc,\P)$ supporting a $d$-dimensional standard Brownian motion $W$. In the following, $\F=(\Fc_{t})_{0\le t\le T}$ will denote the completed right-continuous filtration generated  by $W$. Here, $T>0$ is a finite time horizon. If nothing else is specified, we shall assume that $\Fc_T=\Fc$.

Given a sub-algebra $\Gc\subset \Fc$ and a set $A\subset \R^d$, we write $L^0(A,\Gc)$ for the set of $A$-valued $\Gc$-measurable random variables. We similarly write $L^p(A,\Q,\Gc)$,  $\Q\sim \P$ and $p\in (0,\infty]$, to denote random variables in $L^0(A,\Gc) $ with finite $p$-moment under $\Q$, or essentially bounded if $p=\infty$. When $A$ or $\Gc$ are clearly given by the context, we shall omit them. 

For $p\ge 0$, we  write $L^p_b(\Q,\Gc)$ to denote the collection of element $G\in L^p(\Q,\Gc)$ such that $G\ge -c$ $\Q$-a.s. for some $c>0$. 

The set  predictable processes $\psi$ with values in $\R^d$ satisfying $\E^\Q[\int_0^T|\psi_s|^2ds]<\infty$ is denoted by  $L^2_{\Pc}(\Q)$, or simply $L^2_{\Pc}$ if $\Q=\P$. 

If nothing else is specified $\E$ denote the expectation operator under $\P$. Otherwise, we write $\E^\Q$ if we want to consider the expectation operator under $\Q\ne \P$.

In the following, inequality between random variables have to be understood in the $\Pas$ sense. 

We denote by $x^i$ the $i$-th component of a vector $x\in \R^d$, which will always be viewed as a column vector, with transposed vector $x'$.  We write $|\cdot|$ for the Euclydean norm, and $\M^d$ denotes the set of $d$-dimensional square matrices. We denote by $\mathbb{S}^d$ the subset of elements of $\M^d$  that are symmetric. For a subset $\Oc$ of $\R^d$,  we denote by cl$(\Oc)$ its closure, by int$(\Oc)$ its interior, by $\partial \Oc$ its boundary, and by dist$(x,\Oc)$ the Euclidean distance from $x$ to $\Oc$ with the convention dist$(x,\emptyset)=\infty$.  We denote by $B_r(x)$ the open ball of radius $r>0$ centered at $x\in \R^d$.   If $B=[s,t]\x \Oc$ for $s\le t$ and $\Oc\subset \R^d$, we write $\partial_p B:=([s,t)\x \partial \Oc )\cup (\{t\} \x {\rm cl}(\Oc))$ for its parabolic boundary.

Given a smooth function   $\vp: (t,x)\in \R_+\x \R^d\mapsto \vp(t,x)\in \R$, we denote by $\partial_t \vp$ its derivative with respect to its first variable, and by $D\vp$ and $D^2\vp$ its Jacobian and Hessian matrix with respect to the second one. For $\vp: (t,x_1,\ldots,x_k)\in \R_+\x \R^{kd}\mapsto \vp(t,x)\in \R$, we write $D_{(x_i,x_j)}\vp$ and $D^2_{(x_i,x_j)}\vp$ the Jacobian and Hessian matrix associated to the couple $(x_i,x_j)$.

\section{Financial market and wealth process}

In order to fix ideas and notations, we describe here the typical financial model we have in mind, also more general one will be considered later on. 
\vs2

As usual the financial market will consists in two types of assets. The first one is a risk free asset $B$, often called cash-account, whose dynamics is given by 
$$
B_t=1\p \int_0^t B_s r_s ds =e^{\int_0^t r_s ds}\;\;\;,\;t\ge 0\;,
$$
where $r$ is a predictable real valued process satisfying 
\be\label{hyp: taux sans risque}
\int_0^t  |r_s| ds < \infty \;\; \mbox{ for all } t\ge 0\;.
\ee 

For ease of notation, we also introduce the associated stochastic discount factor $\beta$:
$$
\beta_t:=1/B_t=e^{-\int_0^t r_s ds} \;\;\;,\;t\ge 0\;.
$$

Risky assets (bonds, stocks, derivatives, etc...) are modeled via a $d$-dimensional  process $X=(X^1,\ldots,X^d)$ satisfying 
\b*
X_t=X_0\p \int_0^t \mu_s ds \p\int_0^t \sigma_s dW_s
\e*
where $(\mu,\sigma)$ is a predictable process with valued in $\R^d\x \M^d$ that is bounded on $[0,T]$ $\Pas$  Each component $X^i$ of $X$ denotes a given risky asset.  
 
A financial strategy is described by an element of the set $\Ac$ of   $d$-dimensional predictable processes $\phi$ satisfying
\be\label{hyp: phi strat financiere}
\int_0^t   |\phi_s'|^2 ds <\infty  \;\;  \mbox{ for all } t\le  T\;.
\ee
Each component $\phi^i_t$ denotes the number of units of asset $X^i$ in the portfolio at time $t$.  

To an initial wealth $y\in \R$ and a strategy $\phi\in \Ac$, we associate the portfolio process $Y^{y,\phi}$ defined as
\be\label{eqdef: dyna richesse Y}
Y^{y,\phi}_t:= y \p \int_0^t \phi_s' dX_s \p \int_0^t (Y^{y,\phi}_s-\phi'_s X_s) r_s ds\;,\;t\le T\;.
\ee
In the following, we say that a strategy $\phi$ is admissible, and  we write $\phi \in \Ac_b$, if there exists a constant $c>0$  such that 
\be\label{eq: def admis Acb}
Y^{y,\phi}_t\ge -cB_t  \;\;\mbox{ for all } t\le T\;.
\ee
This condition means that the financial agent has a finite ``credit line'', i.e. his wealth can not go too negative. Note that the constant $c$ may depend on the chosen strategy and is not universal. Moreover, since $Y^{y,\phi}=yB\p Y^{\phi}$, see \reff{eqdef: dyna richesse Y}, the set $\Ac_b$ does not depend on the initial endowment $y$. 
\vs2

For later use observe that $\tilde Y^{y,\phi}:=\beta Y^{y,\phi}$ solves 
\b*
\tilde Y^{y,\phi}_t:= y \p \int_0^t \phi_s' d \tilde X_s \mbox{ and satisfies  }  \tilde Y^{y,\phi}_t\ge -c \;\;\mbox{ for all } t\ge 0\;
\e*
for some $c>0$, 
where $\tilde X:=\beta X$ is given by 
\b*
\tilde X_t=X_0\p \int_0^t (\tilde \mu_s-r_s\tilde X_s) ds \p \tilde \sigma_s dW_s
\e*
with $\tilde \mu:=\beta \mu$ and $\tilde \sigma:=\beta \sigma$. Here, $\tilde Y^{y,\phi}$ and $\tilde X$ can be interpreted as the discounted values of the wealth and financial assets processes.  

\begin{Remark}{\rm 
When dealing with PDE-oriented approaches, we shall specialize to models of the form  (for instance) $r=\rho(X)$, $\mu=\mu(X)$ and $\sigma=\sigma(X)$ where $\rho$, $\mu$ and $\sigma$ will be considered as deterministic functions. In this case, we shall write $(X_{t,x}(s))_s$ for $(X_s)_s$ to insist on the fact that $X$ takes values $x$ at time $t$. We will similarly write $(Y_{t,x,y}^\phi(s))_s$ for $(Y^\phi_s)_s$. More general cases where $\mu$ and $\sigma$ depend on $\phi$ will also be considered. In such a situation, we shall write $(X^\phi_{t,x}(s))_s$ to insist on the dependence of $X$ with respect to the strategy $\phi$.
}
\end{Remark}

\begin{Remark}{\rm Additional constraints will be  imposed later  on strategies. This will allow us to consider more general models where some of the components of $X$ will no more be considered as tradable assets but as non-tradable factors (e.g. stochastic volatility in Markovian models). 
}
\end{Remark}

\section{Hedging problem and hedging criteria}
 
The pricing and hedging problem is the following. We are given a random claim  $G\in L^0(\R,\Fc_T)$ that will impact the wealth of an investor at time $T$. This can be the payoff of a financial derivative that has been sold at time $0$, or any risk related to already engaged positions. 

The question is: what is the amount of money required today in order to be able to construct a financial strategy which will allow to reduce this risk in an appropriate way ? 

Many approaches can be considered depending on the market and the risk tolerance of the investor. 

The first approach consists in trying to make the risk completely disappear. This is the philosophy of the super-hedging point of view: evaluate the risk at its  super-hedging price
\b*
p(G):=\inf\left\{y \in \R~:~\exists\; \phi \in \Ac_b\mbox{ s.t. } Y^{y,\phi}_T\ge G\right\}\;.
\e*
Then, starting from $y=p(G)$, or $y>p(G)$ if the infimum above is not achieved, one can follows a strategy $\phi$ such that $Y^{y,\phi}_T\ge G$, i.e. the risk is completely covered. 

This approach is the most conservative. However, it has two important drawbacks:

1. The associated strategy may not be easy to implement in practice. For instance, it can lead to very large and too quickly varying financial positions. This is typically the case for digital or barrier options for which it can explode near the maturity or the barrier, see e.g. 
\cite{BCS98} and \cite{SSW02}.

2. The computed value may be too large and therefore non-reasonable, see e.g. \cite{CPT99}  for an example of stochastic volatility model in which the super-hedging price of a call is just the spot value of the underlying. 
\vs2

In order to answer the first criticism, one can add portfolio constraints in the model, and compute the corresponding super-hedging price under these  constraints.  
\vs2

As for the second criticism, we need to relax the $\Pas$ super-hedging criteria. 
One way to do this, consists in allowing to miss the hedge with a given probability, i.e. compute the so-called quantile hedging price, see \cite{FoLe99}:
\b*
\inf\left\{y \ge -c~:~\exists\; \phi \in \Ac_b\mbox{ s.t. } \Pro{Y^{y,\phi}_T\ge G}\ge p \right\}\;
\e*
for some $p\in [0,1)$ and $c\in \R_\p$. Here, the constant $c$ is added as a minimum requirement in order to avoid degenerate results.

Another way  consists in allowing to miss the hedge with a level of risk, see \cite{FoLe00}, which leads to problems of the form:
\b*
\inf\left\{y \ge -c~:~\exists\; \phi \in \Ac_b\mbox{ s.t. } \Esp{\ell(Y^{y,\phi}_T- G)}\ge l \right\}\;
\e*
for some $l\in$ Image$(\ell)$ and $c\in \R_\p$. Here, $\ell$ is typically a convex non-decreasing function viewed as a loss function. The map $(y,\phi)\mapsto -\Esp{\ell(Y^{y,\phi}_T- G)}$ has to be interpreted as a measure of the risk induced by starting with $y$ and following the policy $\phi$.

\section{Duality versus stochastic targets}

The above problems have been considered in the literature under the angle of the so-called dual approach. It is based on the relation between super-hedgeable claims and probability measures that  turn   discounted price processes into (local) martingales. This approach allows to appeal to the convex analysis machinery which turns out to be very powerful.

The main drawback of this approach is that it does not allow to consider models where the wealth dynamics in non-linear or in which the financial strategy may have an impact on the price process of financial assets.  

We shall see in these lectures how the recent theory of stochastic targets can handle in a direct way such situations.

%%%%%%%%%%%%%%%%%%%%%%%%%%%%%%%%%%%%%%%%%%%%%%%%%%%%%%%%%%%%%%%%%%%%%%%%%%%%%%%%%%%%%%%%%%%%%%%%%%%%%%%%%%
%%%%%%%%%%%%%%%%%%%%%%%%%%%%%%%%%%%%%%%%%%%%%%%%%%%%%%%%%%%%%%%%%%%%%%%%%%%%%%%%%%%%%%%%%%%%%%%%%%%%%%%%%%
%%%%%%%%%%%%%%%%%%%%%%%%%%%%%%%%%%%%%%%%%%%%%%%%%%%%%%%%%%%%%%%%%%%%%%%%%%%%%%%%%%%%%%%%%%%%%%%%%%%%%%%%%%
\part{Dual approach to  risk based pricing and hedging}

%%%%%%%%%%%%%%%%%%%%%%%%%%%%%%%%%%%%%%%%%%%%%%%%%%%%%%%%%%%%%%%%%%%%%%%%%%%%%%%%%%%%%%%%%%%%%%%%%%%%%%%%%%
%%%%%%%%%%%%%%%%%%%%%%%%%%%%%%%%%%%%%%%%%%%%%%%%%%%%%%%%%%%%%%%%%%%%%%%%%%%%%%%%%%%%%%%%%%%%%%%%%%%%%%%%%%
\chapter{Dual formulation for super-hedging and martingale representation}\label{CHAP: prix dual}

This first part is dedicated to the so-called dual approach. 

%%%%%%%%%%%%%%%%%%%%%%%%%%%%%%%%%%%%%%%%%%%%%%%%%%%%%%%%%%%%%%%%%%%%%%%%%%%%%%%%%%%%%%%%%%%%%%%%%%%%%%%%%%
\section{The complete market case}\label{Section: chap prix dual complet} 

We first consider the so-called complete market case where {\sl any} risk can be covered. 
\vs2

This corresponds to the situation where $\sigma$ is invertible   with bounded inverse on $[0,T]$ $\Pas$   and   the risk premium $\lambda$ defined by 
\b*
\lambda:=\tilde \sigma^{-1} (\tilde \mu - r\tilde X)=\sigma^{-1}(\mu-r X)
\e*
satisfies\footnote{This notation means that $H$ solves $H_t=1- \int_0^t H_s \lambda_s' dW_s\;,\;t\le T$.} 
\be\label{hyp: marche complet}
H:=\Ec\left(-\int_0^\cdot \lambda_s'dW_s\right) \mbox{ is a martingale.}
\ee
so that $\Q\sim \P$ defined by 
\b*
d\Q/d\P=H_T
\e*
is the unique element of the set $\Mc$ of $\P$-equivalent probability measures such that $\tilde X$ is a martingale. 

We then define the $\Q$-Brownian motion $W^\Q$ by 
$$
W^\Q_t=W_t\p \int_0^t \lambda_s ds\;,
$$
recall Girsanov's Theorem, so that 
$$
\tilde X_t=X_0\p \int_0^t \tilde \sigma_s dW^\Q_s 
$$
and therefore 
$$
\tilde Y^{y,\phi}_t=y\p \int_0^t \phi_s' \tilde \sigma_s dW^\Q_s 
$$

\begin{Remark}\label{rem: Y Q super mart}{\rm For $\phi \in \Ac$,  $\tilde Y^{y,\phi}$ is a   $\Q$-local  martingale, i.e. there exists a sequence of stopping times $(\tau_n)_{n\ge 1}$ such that $\tau_n\uparrow \infty$ $\Pas$ and $(\tilde Y^{y,\phi}_{\cdot \wedge \tau_n})$ is a $\Q$-martingale for each $n\ge 1$. Since, for $\phi \in \Ac_b$,   $\tilde Y^{y,\phi}$ is also bounded from below, a straightforward application of Fatou's Lemma shows that it is indeed a $\Q$-supermatingale. 
}
\end{Remark}

\vs2

Under the condition \reff{hyp: marche complet},    any random variable $G$ such that $\beta_TG \in L^1_b(\Q,\Fc_T)$  can be written as the time $T$ value of a wealth process. This is a consequence of the martingale representation theorem. 

\begin{Theorem}\label{thm: representation martingale} Given $G\in L^0$ such that $G\in L^1(\Q,\Fc_T)$ there exists a predictable process $\psi$ satisfying $\int_0^T|\psi_s|^2 ds <\infty$ such that 
$$
\E^\Q[G~|~\Fc_t]=\E^{\Q}[G]\p \int_0^t \psi_s' dW^\Q_s \.
$$
If $G\in L^2(\Q)$, then  $\psi\in L^2_{\Pc}(\Q)$.
\end{Theorem}

Otherwise stated, the $\Q$-martingale $(\E^\Q[G~|~\Fc_t])_{t\le T}$ can be represented in terms of a stochastic integral with respect to $W^\Q$. 

\begin{Corollary}\label{cor: form dual prix couverture complet} Fix   $G\in L^0$ such that  $\beta_T G\in L^1_b(\Q,\Fc_T)$. Then, 
$$
p(G)=\E^\Q[\beta_TG]
$$
and there exists $\phi \in \Ac_b$ such that 
$$
V^{p(G),\phi}_T=G\;. 
$$
If $G\in L^2(\Q)$, then  $\psi\in L^2_{\Pc}(\Q)$.
\end{Corollary}

\proof   For $y>p(G)$, there exists $\phi \in \Ac_b$ such that  $Y^{y,\phi}_T\ge G$. Since $\tilde Y^{y,\phi}$  is a $\Q$-supermatingale, by Remark \ref{rem: Y Q super mart}, this implies that $y\ge \E^{\Q}[\beta_T G]$. 
On the other hand, it follows from Theorem \ref{thm: representation martingale} that there exists a predictable process $\psi$ satisfying $\int_0^T|\psi_s|^2 ds <\infty$ such that 
$$
p(G)\p  \int_0^t \psi_s' dW^\Q_s =\E^\Q[\beta_T G~|~\Fc_t]\;.
$$
By taking $\phi$ defined as $\psi':=\phi'\tilde \sigma$, we obtain 
$$
\tilde Y^{p(G),\phi}_T =\beta_T G\;,
$$
where $\phi$ satisfies $\int_0^T|\phi_s  |^2 ds <\infty$,
 note that $\tilde \sigma^{-1}$ is bounded   on $[0,T]$ $ \P$-a.s.,  and $\tilde  Y^{p(G),\phi}=\E^\Q[\beta_T G~|~\Fc_\cdot]\ge -c$ for some $c>0$. 
\ep 

%%%%%%%%%%%%%%%%%%%%%%%%%%%%%%%%%%%%%%%%%%%%%%%%%%%%%%%%%%%%%%%%%%%%%%%%%%%%%%%%%%%%%%%%%%%%%%%%%%%%%%%%%%
\section{Incomplete markets and portfolio constraints}\label{Section: chap dualite incomplete market}

In order to take into account the incompleteness of the market and possible portfolio constraints, we shall restrict from now on to admissible strategies $\phi \in \Ac_b$ such that $\phi \in K$ $dt\times d\P$-a.e.,  where $K$ is a given convex set of $\R^d$. We denote by $\Ac_K$ the set of such elements. 

\begin{Example} Here are some relevant examples: 

{\rm 1.} Short selling constraints: $K=[0,\infty)^d$. 

{\rm 2.} ``Asset'' 1 can not be traded, no constraint on the others: $K=\{0\}\times \R^{d-1}$. 

{\rm 3.} Bounded positions in any asset: $K=\Prod_{i=1}^d [-m_i,M_i]$ for some $m_i,M_i\ge 0$. 
\end{Example}

\subsection{The general dual formulation}

The aim of this section is to extend the formulation of Corollary \ref{cor: form dual prix couverture complet} to the super-hedging price under constraint: 
$$
p_K(G):=\inf\{y\in \R~:~ \exists\;\phi \in \Ac_K\mbox{ s.t. } Y^{y, \phi}_T\ge G\}\;.
$$

In order to do this, we first need to characterize the set $K$ in term of the support function 
$$
\zeta\in \R^d\mapsto \delta_K(\zeta):=\sup_{\eta \in K} \eta'\zeta.
$$

\begin{Proposition}\label{prop: caract K} $$
\eta \in K \Longleftrightarrow \inf_{|\zeta|=1} \delta_K(\zeta)-\zeta'\eta \ge 0\;.  
$$

\end{Proposition}

\proof The implication $\Rightarrow$ follows from the definition. Conversely, if $\bar \eta\notin K$, which is convex and closed, then the Hahn-Banach separation theorem, see \cite{R70},  implies that there exists $\zeta \in \R^d$  such that $  \sup_{\eta\in K} \eta'\zeta<\bar \eta'\zeta$.   This implies that $\delta_K(\zeta)-\bar \eta'\zeta<0$, where $\zeta$ can always be chosen such that $|\zeta|=1$ by an obvious normalization. \ep
\\

In the following, we let $\Uc_b$ denote the set of $\R^d$-valued predictable processes such that, for some constant $c>0$, 
 $\sup_{s\le T}(|\nu_s|\p |\delta_K(\nu_s)|)\le c$ $\Pas$    For $\nu \in \Uc_b$, we define $\Q^\nu\sim \P$ by  
$$
d\Q^\nu/d\P:= H^\nu_T
$$
where 
$$
H^\nu:=\Ec\left(-\int_0^\cdot (\lambda^\nu_s)' dW_s\right) \;\mbox{ with }\;\lambda^\nu:= \sigma^{-1} (\mu - rX ) - \tilde \sigma^{-1} \nu  \;.
$$
We also define
$$
Z^\nu:=\int_0^\cdot \delta_K(\nu_s) ds \;\mbox{ and the $\Q^\nu$ Brownian motion} \; W^\nu:=W\p \int_0^\cdot \lambda^\nu_s ds \;.
$$

Observe that, for $\nu \in \Ac_K$, 
$$
d ( \tilde Y^{y,\phi}_t -Z^\nu_t)= \left(\phi'_t \nu_t -\delta_K(\nu_t)\right) dt \p \phi_t' \tilde \sigma_t dW^{\nu}_t  \;.
$$

In particular, it follows from Proposition \ref{prop: caract K} that  $\tilde Y^{y,\phi} -Z^\nu$ is a $\Q^\nu$-local supermartingale for any $\phi \in \Ac_K$. Note that, for some $c>0$, $\tilde Y^{y,\phi} -Z^\nu\ge -c-cT$. Hence, this $\Q^\nu$-local supermartingale is bounded from below and is therefore   a $\Q$-super-martingale. This leads to 
the following first result:

\begin{Proposition} Fix $G\in L^0$ such that $\beta_T G\in L^0_b(\Fc_T)$. Then, 
$$
p_K(G)=\inf\{y\in \R~:~\phi \in \Ac_K\mbox{ s.t. } Y^{y, \phi}_T\ge G\}\ge \sup_{\nu \in \Uc} \E^{\Q^\nu}[\beta_T G- Z^\nu_T]\;.
$$
\end{Proposition}

We shall now show that equality actually holds. 

\begin{Theorem}\label{thm: formulation duale sous contrainte}Fix $G\in L^0$ such that $\beta_T G\in L^0_b(\Fc_T)$. Then, 
$$
p_K(G)= \sup_{\nu \in \Uc_b} \E^{\Q^\nu}[\beta_T G- Z^\nu_T]\;.
$$
Moreover, if $p_K(G)<\infty$, then there exists $\phi \in \Ac_K$ such that $Y^{p_K(G),\phi}_T\ge G$.
\end{Theorem} 

We split the proof of the above result in various Lemma.  

Let us now define $P$ as the cadlag adapted process satisfying\footnote{We recall that $\esssup\Ec$, for a family $\Ec$ of random variables, is the smallest random variables which dominates all elements of $\Ec$, in the a.s. sense. } 
 $$
 P_t:=\esssup_{\nu \in \Uc }J_t^\nu \mbox{ where } J^\nu_t:=\E^{\Q^{\nu}}[\beta_TG- (Z^{\nu}_T-Z^\nu_t)~|~\Fc_t]\;,\;t\le T\;
  $$
Note that the existence of a cadlad process satisfying the above property is not obvious. Here, this follows from arguments developed in \cite{stflour} and we omit the details. 

The key argument for proving Theorem \ref{thm: formulation duale sous contrainte} consists in showing that $P$ is a supermatingale under any $\Q^\nu$, $\nu\in \Uc_b$, see Proposition \ref{prop: surmart family} below. 

We first show that the family $\{J_t^\nu\;,\nu \in \Uc_b\}$ is directed upward in the following sense. 

\begin{Definition} We say that a family of random variables $\Ec$ is directed upward is for any $\zeta_1,\zeta_2 \in \Ec$, there exists $\zeta_3 \in \Ec$ such that $\zeta_3\ge \max\{\zeta_1,\zeta_2\}$.
\end{Definition}

\begin{Proposition}\label{prop: J direct upward} For each $t$,  the family $\{J_t^\nu\;,\nu \in \Uc_b\}$ is directed upward. 
\end{Proposition}

\proof Fix $\nu^1,\nu^2\in \Uc_b$, and set $\nu^3=\nu^1\1_{[0,t)}\p \1_{[t,T]} \left(\nu^1\1_{A}\p \nu^2\1_{A^c}\right)$, where $A:=\{J_t^{\nu^1}\ge J_t^{\nu^2}\}$. Clearly, $J_t^{\nu^3}=\max \{J_t^{\nu^1}, J_t^{\nu^2}\}$. Moreover, if $c>0$ is such that $\sup_{s\le T}(|\nu^i_s|\p |\delta_K(\nu^i_s)|)\le c$ $\Pas$ for $i=1,2$, then the same inequality holds for $i=3$. Hence, $\nu^3\in \Uc_b$. 
\ep
\\

In order to prove Proposition \ref{prop: surmart family}, we now use the following well-know property of directed upward families, see e.g. \cite{Neveu}. 

\begin{Lemma}\label{lem: famille directed upward} If $\Ec$ is a family directed upward. Then there exists a sequence $(\zeta_n)_{n\ge 1}\subset \Ec$ such that $\esssup \Ec=\lim\limits_{n\to \infty}\uparrow \zeta_n$.
\end{Lemma}

We can now prove the supermartingale property.

\begin{Proposition}\label{prop: surmart family} For all $\nu \in \Uc_b$,  $P-Z^\nu$ is a  $\Q^\nu$-supermartingale.
\end{Proposition} 

\proof  Fix $t\ge s$ and  $\nu \in \Uc_b$. Let $(\nu_n)_{n\ge 1}$ be such that
 $J^{\nu^n}_t\uparrow P_t$ as $n\to \infty$, see Lemma \ref{lem: famille directed upward}  and Proposition \ref{prop: J direct upward}. 
For $\nu \in \Uc_b$, set $\bar \nu_n:=\nu\1_{[0,t)}\p \nu_n\1_{[t,T]}$. 
Then, 
\b*
\E^{\Q^{\nu}}[ P_t-Z^{\nu}_t~|~\Fc_s]&=&\E^{\Q^\nu}\left[\lim_{n\to \infty} \uparrow \E^{\Q^{\nu_n}}[\beta_TG- (Z^{\nu_n}_T-Z^{\nu_n}_t)~|~\Fc_t]-Z^\nu_t  ~|~\Fc_s\right]
\\
&=&\lim_{n\to \infty} \uparrow \E^{\Q^\nu}\left[\E^{\Q^{\nu_n}}[\beta_TG- (Z^{\nu_n}_T-Z^{\nu_n}_t)~|~\Fc_t]-Z^\nu_t  ~|~\Fc_s\right] 
\\
&=&\lim_{n\to \infty} \uparrow \E^{\Q^{\bar \nu_n}}\left[ \beta_TG- (Z^{\bar \nu_n}_T-Z^{\bar \nu_n}_s) ~|~\Fc_s\right] -Z^\nu_s 
\\
&\le& P_s-Z^\nu_s\;.
\e*
\ep 

\begin{Proposition} For each $\nu \in \Uc_b$, there exists a $\Q^\nu$-martingale $M^\nu$ and a non-decreasing process $A^\nu$ such that $A^\nu_0=0$ and  $P-Z^\nu=M^\nu-A^\nu$. 
\end{Proposition}

\proof This follows from the Doob-Meyer decomposition together with the previous proposition.
\ep\\

In order to conclude the proof, we now apply the martingale representation to $M^0$ to obtain some predictable process $\psi$ satisfying $\int_0^T|\psi_s|^2 ds <\infty$ such that 
$$
P_t=P_t-Z^0_t = P_0 \p \int_0^t \psi'_s dW^0_s-A^0_t\;. 
$$
By taking $\phi$ such that $\phi' \tilde \sigma=\psi'$, we obtain 
$$
P_t = \tilde Y^{P_0,\phi}_t - A^0_t =\esssup_{\nu \in \Uc_b }\E^{\Q^{\nu}}[\beta_TG- (Z^{\nu}_T-Z^\nu_t)~|~\Fc_t]\ge \E^\Q[\beta_TG~|~\Fc_t]\;,\;t\le T,
$$
which implies that $\phi \in \Ac_b$ and that $Y^{P_0,\phi}_T\ge G$, since $A^0\ge 0$. 
To conclude the proof, it remains to shows that $\phi \in K$ $dt\x d\P$-a.e. 
To see this, recall that $P$ can also be decomposed as $P-Z^\nu=M^\nu-A^\nu$. In particular, we must have 
\b*
P-Z^\nu&=&P_0 \p \int_0^\cdot \psi_s'dW^0_s- A^0- Z^\nu \\
&=&P_0 \p \int_0^\cdot \psi_s'dW^\nu_s\p  \int_0^\cdot  (\psi_s' \tilde \sigma_s^{-1} \nu_s -\delta_K(\nu_s) )ds- A^0 \\
&=& P_0 \p \int_0^\cdot \psi_s'dW^\nu_s\p  \int_0^\cdot  (\phi_s' \nu_s -\delta_K(\nu_s)) ds- A^0
\e* 
so that $A^\nu=A^0-\int_0^\cdot  (\phi_s'  \nu_s -\delta_K(\nu_s) ds$ which is therefore non-decreasing. It follows that 
$$
\int_0^\cdot  (\phi_s' \nu_s -\delta_K(\nu_s) )ds\le A^0_T
$$
for all $\nu \in \Uc_b$. By replacing $\nu$ by $n\nu$ and by sending $n\to \infty$, we deduce from the above inequality that 
$$
\int_0^\cdot  (\phi_s' \nu_s -\delta_K(\nu_s) )ds\le 0\;.
$$
Let us now define   $\bar \nu$ as $\bar \nu:=\arg\min_{|\zeta|=1} (\delta_K(\zeta)-\phi'\zeta)$. Taking $\nu:=\bar\nu\1_{\{\delta_K(\bar\nu)-\phi'_s\bar\nu<0\}}$ in the last inequality, shows  that 
 $\phi\in K$ $dt\x d\P$-a.e, recall  Proposition \ref{prop: caract K}.\ep 

\subsection{Examples}
We conclude this section with three examples of applications. The first one corresponds to a Brownian model with portfolio constraints, the second one to a Black-Scholes model with constraints on the amount of money   invested in the asset, the last one to a stochastic volatility model. 

\begin{Example}{\bf (Brownian model with portfolio constraint)}
Let us consider the case $d=1$ where  $X=X^1$ has the dynamics
$$
X_t=X_0\p \mu t \p \sigma W_t\;\;\;t\le T\;,
$$
and $r=0$. We want to hedge an option of payoff $g(X_T)$ paid at time $T$ under the constraints $K=[-m,M]$ with $M,m\ge 0$. We shall assume here that $g$ is non-decreasing.

In this case, $\delta_K(\zeta)=\zeta^\p M \p \zeta^- m$ so that dom$(\delta_K)=\R$. Let us define the function $\hat g$ by $\hat g(x):=\sup_{u\in \R} \left(g(x\p u)-(u^\p M \p u^- m)\right)$. Then, it follows from Theorem \ref{thm: formulation duale sous contrainte} that:
\b*
p_K(G)&=&\sup_{\nu \in \Uc_b}\E^{\Q^\nu}\left[g(  X_T)-\int_0^T   (\nu_s^\p M \p \nu_s^- m)\right]\;
\\
&=& \sup_{\nu \in \Uc_b}\E^{\Q^{  \nu}}\left[g\left(  X_0\p \int_0^T\nu_s ds \p \sigma W^{ \nu}_T \right)-\int_0^T   \delta_K(\nu_s) ds \right]
\\
&\le& \sup_{\nu \in \Uc_b}  \E^{\Q^{  \nu}} \left[  \hat g\left( X_0  \p \sigma W^{  \nu}_T  \right) \right]
\e*
where we used the fact that $g(x)=g(x\p u - u  )\le \hat g(x-u)\p  \delta_K(u)$.
It follows that 
\b*
p_K(G)&\le &\E^{\Q}\left[  \hat g\left(X_0  \p \sigma W^{\Q}_T\right)\right]\;.
\e*
We now observe that, by a formal identification of the law of $W^\Q$ under $\Q$ and $W^\nu$ under $\Q^\nu$, 
\b*
p_K(G)
&=& \sup_{\nu \in \Uc_b}\E^{\Q}\left[g\left(  X_0\p \int_0^T\nu_s ds \p \sigma W^{\Q}_T \right)-\int_0^T   \delta_K(\nu_s) ds \right],
\e*
 see  \cite{SSW02}  for a rigorous argument.
 Moreover, any bounded $\Fc_t$-measurable random variable, with $t<T$, can be written in the form $\int_0^T\nu_s ds$ with $\nu \in \Uc_b$. Indeed, give $\xi \in L^0(\Fc_t)$, one has $\int_0^T (\xi/(T-t))\1_{s\ge t} ds =\xi$. Given $\xi \in L^\infty(\Fc_T)$, one can then approximate it by the sequence $\Esp{\xi~|~\Fc_{T(1-1/n)}}_{n\ge 1}$. 
It follows that, for $g$ continuous and bounded from below, 
\b*
p_K(G)&\ge &\sup_{\xi \in L^\infty(\R_\p,\Fc_T)}\E^{\Q}\left[  g\left(X_0  \p \sigma W^{\Q}_T \p \xi  \right) - \xi^\p M \right]\;.
\e*
Here, we restrict to non-negative random variable because $g$ is non-decreasing and it should therefore be optimal to restrict to $\nu \ge 0$ or equivalently $\xi \ge 0$. 
Now, we clearly have 
\b*
&&\sup_{\xi \in L^\infty(\R_\p,\Fc_T)}\E^{\Q}\left[ g\left(X_0  \p \sigma W^{\Q}_T\p \xi  \right) - \xi^\p M \right]
\\
&&=
\E^{\Q}\left[\sup_{\zeta\in \R_\p} \left(  g\left(X_0  \p \sigma W^{\Q}_T\p  \zeta  \right) - \zeta M  \right)\right]\;.
\e*
This shows that 
\b*
p_K(G)&=& \E^{\Q}\left[ \hat g\left(X_0  \p \sigma W^{\Q}_T   \right)  \right]\;,
\e*
i.e., the price under constraint for the option $g$ is the usual unconstrained price  in the Brownian  model of the face-lifted payoff $\hat g$. 
\end{Example} 

\begin{Example}{\bf (Black-Scholes model with portfolio constraint)}

Let us now consider the Black-Scholes model where $X$ is given by 
$$
dX_t/X_t= \mu dt \p \sigma dW_t
$$
and $r=0$ for simplicity. In this example, we impose the constraint 
$$
\psi:=\phi X  \in K\;\;dt\times d\P{\rm -a.e.}
$$
i.e. the amount  invested in the risky asset belongs to $K$. Let $\hat A_K$ denote the set of processes $\phi \in \Ac_b$ such that the above constraint is satisfied.  
\vs2

We shall see how we can reduce the problem of super-hedging a claim $g(X_T)$ to the problem discussed in the previous example. 

\vs2

To do this, first observe that  
$$
Y^{y,\phi}_t=y\p \int_0^t \phi_s dX_s=y\p \int_0^t \psi_s  dX_s/X_s=y\p \int_0^t \psi_s  \mu ds \p \int_0^t \psi_s   \sigma dW_s
$$
where $\psi:=\phi X$, so that 
$$
Y^{y,\phi}_t=y \p \int_0^t \psi_s d\bar X_s\;
$$
with 
$$
\bar X_t:=\mu t \p \sigma W_t\;. 
$$
It follows that, at least for $g$ bounded from below, 
\b*
\hat p_K(g(X_T))&:=&\inf\left\{y~:~\exists\;\phi \in \hat \Ac_K \mbox{ s.t. }   Y^{y,\phi}_T\ge g(X_T)\right\}\\
&=&  \inf\left\{y~:~\exists\;\psi \in   \Ac_K \mbox{ s.t. }   y\p \int_0^T\psi_s d\bar X_s \ge \bar g(\bar X_T)\right\}
\e*
where $\bar g(x):=g(X_0e^{x-(\sigma^2/2) T})$. 

Letting $\bar p_K$ be defined as $p_K$ but for the model where the stock price is given by $\bar X$, the above arguments show that 
$$
\hat p_K(g(X_T))=\bar p_K(\bar g(\bar X_T))\;.
$$
In view of the previous example, one can then obtain an explicit formulation for $\hat p_K(g(X_T))$.
\end{Example} 

\begin{Example}\label{ex: dual stochastic vol}{\bf (Stochastic volatility)} In this example, we take $d=2$ and let $(X^1, X^2)$ be the solution of 
\b*
X^1_t&=&X^1_0\p \int_0^t X^1_s r ds \p  \int_0^t X^1_s \sigma(X^2_s) dW^1_s\\
X^2_t&=&X^2_0\p \gamma_1 W^1_t \p  \gamma_2 W^2_t
\e*
where   $\gamma_1,\gamma_2>0$, $\sigma\ge \eps$ for some $\eps>0$ and $\sigma$ is bounded. We impose the constraint $K:=\R\x\{0\}$, i.e. $X^2$ can not be traded. This corresponds to the simplest stochastic volatility model, in which $X^2$ should be considered as a factor driving the volatility of $X^1$, and not as an asset. 

In this case, we have $\delta_K(\zeta)=0$ is $\zeta^1=0$ and $\delta_K(\zeta)=\infty$ otherwise. It follows that 
$$
p_K(g(X^1_T))=\sup_{\lambda \in \Lambda} \E^{\Q^\lambda}\left[\beta_T g(X^1_T)\right]
$$
where $\Lambda$ denotes the set of real valued  predictable processes $\lambda$ satisfying $\sup_{s\le T}|\lambda_s|\le c$ for some $c>0$, and $\Q^\lambda$ is defined by 
$$
\frac{d\Q^\lambda}{d \P} = e^{-\frac12 \int_0^T (\gamma_2^{-1} \lambda_s)^2 ds \p \int_0^T  \gamma_2^{-1} \lambda_s dW^2_s }\;,
$$ 
which, up to the boundedness imposed on $\lambda$, corresponds to the family of all martingale measures for $X^1$. 

We shall come back to this example   in Chapter \ref{CHAP: pricing equation II} below. 
\end{Example}
%%%%%%%%%%%%%%%%%%%%%%%%%%%%%%%%%%%%%%%%%%%%%%%%%%%%%%%%%%%%%%%%%%%%%%%%%%%%%%%%%%%%%%%%%%%%%%%%%%%%%%%%%%
%%%%%%%%%%%%%%%%%%%%%%%%%%%%%%%%%%%%%%%%%%%%%%%%%%%%%%%%%%%%%%%%%%%%%%%%%%%%%%%%%%%%%%%%%%%%%%%%%%%%%%%%%%
\chapter{The pricing equation I: the complete market case}\label{CHAP: pricing equation complet}

In this chapter, we restrict to the Markovian setting where $X$ is given as the solution of an SDE of the form 
\be\label{eq: dyna X markov prix dual}
X_{t,x}(s)=x\p \int_t^s r_{t,x}(u) X_{t,x}(u) du \p \int_t^s \sigma(X_{t,x}(u)) dW^\Q_u \;,
\ee
for  a risk free  interest rate of the form   
$$
r_{t,x}=\rho(X_{t,x})
$$  
where  $\rho$, $\mu$ and $\sigma$ are assumed to be Lipschitz continuous, and $\rho$ is such that $\rho^-$ is bounded and  $x\mapsto \rho(x)x$ is Lipschitz continuous.
\vs2

For ease of notations, we shall only consider the case where $X$ can take any values in $\R^d$, also in most  financial models we should typically restrict to $(0,\infty)^d$. The arguments being the same in this last case. 

\vs2 

The aim of this section is to provide a PDE formulation for the price function of an option of payoff $g(X_{t,x}(T))$ paid at time $T$, depending on the initial time $t$ and the initial value of $X$ at this time. \vs2

In the following, $g$ will be assumed to be continuous with linear growth and uniformly bounded from below. 
%%%%%%%%%%%%%%%%%%%%%%%%%%%%%%%%%%%%%%%%%%%%%%%%%%%%%%%%%%%%%%%%%%%%%%%%%%%%%%%%%%%%%%%%%%%%%%%%%%%%%%%%%%
\section{Problem extension and dynamic programming}

Motivated by Section \ref{Section: chap prix dual complet} of Chapter \ref{CHAP: prix dual}, we now introduce the pricing function associated to the complete market case:
$$
(t,x)\in [0,T]\x \R^d \mapsto v(t,x):=\E^\Q[\beta_{t,x}(T) g(X_{t,x}(T))]
$$
where 
$$
\beta_{t,x}:=e^{-\int_t^\cdot  \rho(X_{t,x}(s))ds }\;.
$$

The key assertion for deriving a PDE associated to $v$ is the following {\sl dynamic programming} equation which relates the time $t$ value of the price to its time $\theta$ value, for any stopping time $\theta$ bigger than $t$. 
In the following, we shall denote by $\Tc_{[t,\tau]}$ the collection of stopping times taking values in $[t,T]$.

\begin{Proposition}\label{prop: PD FK} For all $\theta\in \Tc_{[t,T]}$, we have
	\be\label{eq: PD FK}
	v(t,x)=\E^\Q\left[\beta_{t,x}(\theta)v(\theta,X_{t,x}(\theta))\right]\;. 
	\ee
\end{Proposition} 

\proof   By the flow property of $X$ and the usual tower property,  we have
	\b*
	v(t,x)
	&=&
	\E^\Q\left[ \beta_{t,x}(\theta) \E^\Q\left[\beta_{\theta,X_{t,x}(\theta)}(T)g(X_{\theta,X_{t,x}(\theta)}(T))~|~\Fc_\theta\right]\right]\;.
	\e* 
It then follows from  the strong Markov property of $X$ defined by \reff{eq: dyna X markov prix dual}  that 
	\b*
  v(\theta,X_{t,x}(\theta))&=&\EspQ{\beta_{\theta,X_{t,x}(\theta)}(T)g(X_{\theta,X_{t,x}(\theta)}(T))~|~(\theta,X_{t,x}(\theta))}
  \\
  &=&\E^\Q\left[\beta_{\theta,X_{t,x}(\theta)}(T)g(X_{\theta,X_{t,x}(\theta)}(T))~|~\Fc_\theta\right],
  	\e*
hence the required result.
\ep

%%%%%%%%%%%%%%%%%%%%%%%%%%%%%%%%%%%%%%%%%%%%%%%%%%%%%%%%%%%%%%%%%%%%%%%%%%%%%%%%%%%%%%%%%%%%%%%%%%%%%%%%%%
\section{Feynman Kac representation in the smooth case}

Using the above proposition, we   can now show that, whenever it is smooth enough, $v$ solves the PDE
\be\label{eq pde FK} 
	\Lc^\Q v=\rho v 
	\ee
on $[0,T)\x \R^d$ with the boundary condition $v(T,\cdot)=g$. Here,  $\Lc^\Q$  is  the {Dynkin operator}  associated to $X$ under $\Q$:
	\b*
	\Lc^\Q \vp(t,x)&:=& 
	 {\partial_t} \vp(t,x) + \rho(x)x'D\vp(t,x)  + \frac12 {\rm Tr}\left[\sigma\sigma'(x)D^2\vp(t,x)\right] \;.
\e*

\subsection{Derivation} 
	
\begin{Theorem}(Feynman-Kac)\label{thm FK regu} Assume that    $v$ is continuous on $[0,T]\x \R^d$ and $v \in  C^{1,2}([0,T)\x \R^d)$. Then,  $v$ is a solution on $[0,T)\x \R^d$ of 
\reff {eq pde FK} 
and satisfies the boundary condition $\lim_{t\nearrow T,z\to x} v(t,z)=g(x)$ on $\R^d$. 
\end{Theorem}

\proof The boundary condition is a consequence of the continuity assumption on $v$. 
It remains to show that $v$ solves \reff{eq pde FK}.  
We now fix $(t,x) \in [0,T)\x \R^d$.    Let $\theta$ be the first time when $(s,X_{t,x}(s))_{s\ge t}$ exits a given bounded  open neighborhood of $(t,x)$. Set $\theta^h=\theta \wedge (t+h)$ for $h>0$ small.  
Using Proposition \ref{prop: PD FK} and  It\^{o}'s Lemma, we deduce that 
	\be\label{eq: E div par h}
	0=\Esp{\frac1h \int_t^{\theta^h}  \beta_{t,x}(s)\left(\Lc^\Q v (s,X_{t,x}(s))-(\rho v) (s,X_{t,x}(s))\right)ds}\;. 
	\ee
Now, we observe that $s\mapsto X_{t,x}(s)$ is $\Pas$ continuous, so that  $|X_{t,x}(s\wedge (t\p h)) - x|\to 0$ $\Pas$ as $h\to 0$ for each $s\ge t$.  Moreover,  $\theta>0$ $\Pas$ so that   $(\theta^h-t)/h\to 1$ $\Pas$ Using the mean value theorem and the continuity of $\Lc^\Q v-\rho v$, we then deduce that 
	\b*
	&&\frac1h \int_t^{\theta^h}  \beta_{t,x}(s)\left(\Lc^\Q v (s,X_{t,x}(s))-(\rho v) (s,X_{t,x}(s))\right)ds
	\\
	&&\to  (\Lc^\Q v - \rho v) (t,x)  \;\;\Pas
	\e*
as $h\to 0$. The required result is then obtained by applying the dominated convergence theorem to pass to the limit in \reff{eq: E div par h}, observe that $(s,X_{t,x}(s))_{s\ge t}$ is bounded on $[t,\theta]$ by definition of $\theta$.
\ep

\subsection{Comparison and uniqueness}

In order to show that Theorem \ref{thm FK regu} provides a full characterization of $v$, it remains to show that $v$ is the unique solution of \reff {eq pde FK} within a suitable class of functions. This is a consequence of the following {\sl comparison result}. 

\begin{Theorem}\label{thm: comparaison regu FK}(Comparison principle) Assume that   $U$ and $V$ are continuous on  $[0,T]\x \R^d$ and  $C^{1,2}$ on $[0,T)\x \R^d$. Assume further that, on  $[0,T)\x \R^d$, 
	\be\label{eq comp regu FK sur sous sol}
	\Lc^\Q U\le \rho U &\mbox{ and } &  \Lc^\Q V \ge \rho V
	\ee
and that $U(T,\cdot)\ge V(T, \cdot)$ on $\R^d$. Finally assume that $U$ and  $V$  have polynomial growth. Then, $U\ge V$ on $[0,T]\x \R^d$. 
\end{Theorem}
 
 \proof By possibly replacing $U$ and $V$ by $\tilde U(t,x):=e^{\kappa t} U(t,x)$ and $\tilde V(t,x):=e^{\kappa t} V(t,x)$ for a large $\kappa$, we can assume that $\rho\ge \eta$ on $\R^d$ for some $\eta>0$. Indeed, $\tilde U$ and $\tilde V$ would satisfy \reff{eq comp regu FK sur sous sol} with    $\rho U$ and $\rho V$ replaced by $(\rho\p \kappa) \tilde U$ and  $(\rho\p \kappa) \tilde V$, where   $\rho^-$ is bounded.  Assume now that, for some $(t_0,x_0) \in  [0,T]\x \R^d$, we have $U(t_0,x_0)<V(t_0,x_0)$. We shall show that this leads to a contradiction. Fix $\eps>0$, $\kappa>0$ and $p$ an integer greater that $1$ such that $\limsup_{|x|\to \infty} \sup_{t\le T} (|U(t,x)|+|V(t,x)|)/(1+|x|^p)=0$. Then, there is $(\hat t ,\hat x ) \in  [0,T]\x \R^d$ such that, for $\eps$ small enough,
	\b*
	0&<&V(\hat t,\hat x)-U(\hat t,\hat x)-\phi(\hat t,\hat x)
	=
	\max_{(t,x)\in  [0,T]\x \R^d} \left(V(t,x)-U(t,x)-\phi(t,x)\right)\;, 
	\e*
where 
	\b*
	\phi(t,x):=\eps e^{-\kappa   t}  (1+| x|^{2p})\;. 
	\e*
Since $U\ge V$ on $\{T\}\x\R^d$, we must have $\hat t<T$. Moreover, the one and second order conditions of optimality imply
	\b*
	\partial_t V(\hat t,\hat x)\le  (\partial_t U+\partial_t \phi)(\hat t,\hat x)\;,\; DV(\hat t,\hat x)= (DU+D\phi)(\hat t,\hat x)  \;
	\e*
 and 
	\b*
	D^2V(\hat t,\hat x) &\le& (D^2U+D^2\phi)(\hat t,\hat x) 
	\e*
in the sense of matrices. Combined with \reff{eq comp regu FK sur sous sol}, this leads to
	\b*
	\rho(V-U)(\hat t,\hat x) &\le& \Lc^\Q(V-U)(\hat t,\hat x) 
	\\
	&\le&\partial_t \phi(\hat t,\hat x) + \rho(\hat x)\hat x'D \phi(\hat t,\hat x)
	+     {\rm Tr}\left[\sigma\sigma'(\hat x)D^2  \phi(\hat t,\hat x)\right] 
	\\
	&\le&
 \Lc^\Q\phi(\hat t,\hat x) \;.
	\e*
 Since $x\mapsto\rho(x)x$ and $x\mapsto\sigma(x)$  have linear growth, we can choose $\kappa>0$ sufficiently large so that 
 \b*
 \Lc^\Q\phi  = -\kappa \phi \p \rho   x'D \phi + {\rm Tr}\left[\sigma\sigma'D^2  \phi\right] <0\mbox{ on } [0,T]\x \R^d\;.
 \e*
 This contradicts
	$(V-U)(\hat t,\hat x)>0$ since $\rho\ge \eta>0$.
\ep 

\begin{Corollary}\label{cor: unicite FK} Assume that  $v$ is  $C^{1,2}([0,T)\x \R^d)\cap C^{0}([0,T]\x \R^d)$, then it is the unique $C^{1,2}([0,T)\x \R^d)\cap C^{0}([0,T]\x \R^d)$ solution of \reff{eq pde FK} satisfying $v(T,\cdot)=g$ in the class of solutions with polynomial growth. If $g$ is bounded from below and Lipschitz continuous, then there exists $\phi \in \Ac_b$ such that $Y^{v(0,X_0),\phi}_T=g(X_{0,X_0}(T))$ and $\phi=Dv(\cdot,X_{0,X_0})$ on $[0,T)$.
\end{Corollary}

\proof Since  $g$ has  linear growth and $\rho$ is bounded, we deduce from standard    estimates    that $v$ has  linear growth too. The first result then follows from Theorems \ref{thm FK regu} and \ref{thm: comparaison regu FK}. Moreover, an application of It\^o's Lemma implies that:
\b*
v(0,X_0)\p \int_0^T \beta^0(t)Dv(t,X^0(t))\sigma(X^0(t)) dW^\Q_t
&=&
\beta^0(T)v(T,X^0(T))\\
&=&\beta^0(T)g(X^0(T))
\e* 
where $X^0:=X_{0,X_0}$ and $\beta^0:=\beta_{0,X_0}$, which is equivalent to $Y^{y_0,\phi}_T=g(X^0(T))$ with $\phi=Dv(\cdot,X^0)$ on $[0,T)$ and $y_0:=v(0,X_0)$. Since $g$ is bounded from below and $\rho^-$ is bounded, we have $\beta^0(T)g(X^0(T))$ bounded from below. Moreover, the fact that $g$ and all the parameters  are Lipschitz continuous implies, by standard estimates, that $v$ is Lipschitz continuous in $x$, uniformly in time. This implies that $Dv$ is bounded so that $\tilde Y^{y_0,\phi}$ is a martingale such that $\tilde Y^{y_0,\phi}(T)$ is bounded from below. Hence, it is bounded from below on the time interval $[0,T]$. 
\ep

\subsection{Verification theorem}

In practice, the regularity assumptions of the above theorem are very difficult to check and we have to rely on a weaker definition of solutions, like viscosity solutions (see e.g. \cite{CIL92} and below), or to use a  {verification theorem} which essentially consists in showing that, if a smooth solution  of \reff{eq pde FK}  exists, then it coincides with $v$. 

\begin{Theorem}\label{thm verification FK}(Verification) Assume that there exists a $C^{1,2}([0,T)\x \R^d)$ solution $\vp$ to \reff{eq pde FK} with polynomial growth  such that 
	\be\label{eq cond bord verif FK}
	\lim_{t\nearrow T, z\to x} \vp(t,z)= g(x) \;\;\;\mbox{ on $\R^d$}\;.
	\ee
Then, $v=\vp$. 
\end{Theorem}

\proof Given $n\ge 1$, set 
	\b*
	\theta_n := \inf\{s\ge t:~|X_{t,x}(s)| \ge n\} \;. 
	\e*
Note   that  $X_{t,x}$ is bounded on $[t,\theta_n\wedge T]$. 
By It\^{o}'s Lemma  and the fact that $\vp$ solves \reff{eq pde FK}, we obtain
	\begin{align}\label{eq verif n a jeter FK}
	\vp(t,x)=\EspQ{\beta_{t,x}(\theta_n\wedge T)\vp(\theta_n\wedge T,X_{t,x}(\theta_n\wedge T))}
	\end{align}
for each $n$. Now, observe that $\theta_n\to \infty$ as $n\to \infty$. In view of \reff{eq cond bord verif FK}, this  implies that  
	\b*
	\beta_{t,x}(\theta_n\wedge T)\vp(\theta_n\wedge T,X_{t,x}(\theta_n\wedge T)) 
	\lra \beta_{t,x}(  T)g(X_{t,x}(  T))  \;\;\Pas
	\e*
Moreover, standard estimates, based on the fact that $v$ has polynomial growth, that $\rho$ is bounded from below, and on the Lipschitz continuity of the coefficients, imply that 
the sequence $(\beta_{t,x}(\theta_n\wedge T)\vp(\theta_n\wedge T,X_{t,x}(\theta_n\wedge T)))_{n\ge 1}$ is uniformly integrable.
We then deduce   that $\vp=v$ by sending $n\to \infty$ in  \reff{eq verif n a jeter FK} and using the dominated convergence theorem.  
\ep

%%%%%%%%%%%%%%%%%%%%%%%%%%%%%%%%%%%%%%%%%%%%%%%%%%%%%%%%%%%%%%%%%%%%%%%%%%%%%%%%%%%%%%%%%%%%%%%%%%%%%%%%%%
\section{Feynman Kac representation in the viscosity sense}

Except when   $\sigma$ satisfies the following type of uniform ellipticity condition
	\be\label{eq sigma unif ellip}
	\exists \;c>0\mbox{ s.t. } \xi' \sigma\sigma' \xi\ge c|\xi|^2\; \mbox{ for all } \xi \in \R^d,
	\ee
it is difficult to show (and in general not true) that $v$ is $C^{1,2}$. Still, it can be shown to solve \reff{eq pde FK}  in a weak sense: the viscosity sense. In the subsections below, we explain this notion and show that 
$v$ is the unique  viscosity solution of \reff{eq pde FK} satisfying $v(T-,\cdot)=g$, in the class of continuous functions with polynomial growth. We refer to \cite{CIL92} for a general overview of the theory of viscosity solutions. 

%%%%%%%%%%%%%%%%%%%%%%%%%%%%%%%%%%%%%%%%%%%%%%%%%%%%%%%%%%%%%%%%%%%%%%%%%%%%%%%%%%%%%%%%%%%%%%%%%%%%%%%%%%
\subsection{Viscosity solutions: definition and main properties}

Let $F$ be an operator from  $[0,T]\x \R^d\x \R\x\R\x \R^d \x\mathbb{S}^d$ into $\R$, where $\mathbb{S}^d$ denotes the set of $d$-dimensional symmetric matrices. In this section, we will be mostly interested by the case
	\begin{equation}\label{eq: def F FK}
	F(t,x,u,q,p,A)=\rho(x)u- q -\rho(x)x'p  - \frac12 {\rm Tr}\left[\sigma\sigma'(x)A\right] \;,
	\end{equation}
 so that $v$ solves  \reff{eq pde FK} means
	\begin{equation}\label{eq: def F FK avec I}
	F(t,x,v(t,x),\partial_t v(t,x),Dv(t,x),D^2v(t,x))=0\;.
	\end{equation}

\vs2

We say that $F$ is {elliptique} if it is non increasing with respect to $A \in \mathbb{S}^d$. This is clearly the case for $F$ defined as in \reff{eq: def F FK}.  In the following, $F$ will always be assumed to be elliptic. 

\vs2

Let us assume for a moment that $v$ is smooth. Let $\vp$ be $C^{1,2}$ and $(\hat t,\hat x)$ be a (global) minimum point of $v-\vp$. After possibly adding a constant to $\vp$, one can always assume that $(v-\vp)(\hat t,\hat x)=0$. In this case, the first and second order optimality conditions imply
	$$
	(\partial_t v,Dv)(\hat t,\hat x)=(\partial_t \vp,D\vp)(\hat t,\hat x) \mbox{ and }  D^2v(\hat t,\hat x)\ge D^2\vp(\hat t,\hat x)\;. 
	$$
Since $F$ is elliptic and $v\ge \vp$ on the domain with equality at $(\hat t,\hat x)$, we deduce that 
	$$
	F(\hat t,\hat x, \vp(\hat t,\hat x),\partial_t \vp(\hat t,\hat x),D\vp(\hat t,\hat x),D^2\vp(\hat t,\hat x))\ge 	0\;
	$$ 
whenever 
	$$
	F(\hat t,\hat x,v(\hat t,\hat x),\partial_t v(\hat t,\hat x),Dv(\hat t,\hat x),D^2v(\hat t,\hat x))=0\;. 
	$$ 
Conversely, if  $(\hat t,\hat x)$ is a (global)  maximum point of $v-\vp$ then 
	$$
	F(\hat t,\hat x, \vp(\hat t,\hat x),\partial_t \vp(\hat t,\hat x),D\vp(\hat t,\hat x),D^2\vp(\hat t,\hat x))
	\le
	0\;.
	$$ 

\vs2

This leads to the following notion of viscosity solution. 

\begin{Definition}\label{def solution viscosite} Let $F$ be an elliptic operator as defined above.  We say that a l.s.c. (resp. u.s.c) function $U$ is a {\sl supersolution} (resp. {\sl subsolution}) of  \reff{eq: def F FK avec I} on $[0,T)\x \R^d$ if for all $\vp \in C^{1,2}$ and  $(\hat t,\hat x)\in [0,T)\x \R^d$ such that $0=\min_{[0,T]\x \R^d} (U-\vp)=(U-\vp)(\hat t,\hat x)$ (resp. $0=\max_{[0,T]\x \R^d} (U-\vp)=(U-\vp)(\hat t,\hat x)$), we have:
	\begin{align}\label{eq F ge 0 def visco fonction teste}
	&F(\hat t,\hat x, \vp(\hat t,\hat x),\partial_t \vp(\hat t,\hat x),D\vp(\hat t,\hat x),D^2\vp(\hat t,\hat x))\ge 0&\\
	& \;\;\;\;(\mbox{ resp. $\le 0$})\;.&\nonumber
	\end{align}
\end{Definition}

We shall say that a locally bounded function is a {discontinuous viscosity solution} of   $F=0$  if $U_*$ and $U^*$ are respectively super- and subsolution, where, for $(t,x) \in [0,T]\x \R^d$, 
	\b*
	U_*(t,x)=\liminf_{(s,y)\in [0,T)\x \R^d\to (t,x)}U(s,y)\;\mbox{ and } \;U^*(t,x)=\limsup_{(s,y)\in [0,T)\x \R^d\to (t,x)}U(s,y)\;.~~~~ 
	\e*
If $U$ is continuous, we simply say that it is a  {viscosity solution}. 
\vs2 

Note that a smooth solution $U$ is also a viscosity solution, as any point achieves a min (or max) of $U-U$. 

\begin{Remark}\label{rem min strict}{\rm If   $(\hat t,\hat x)\in [0,T)\x \R^d$ achieves a  minimum of  $U-\vp$ then it achieves a strict minimum of $U-\bar \vp$ where  $\bar \vp(t,x)=\vp(t,x)-  |x-\hat x|^4-|t-\hat t|^2$. Moreover, if $\bar  \vp$ satisfies \reff{eq F ge 0 def visco fonction teste} at  $(\hat t,\hat x)$ then $\vp$ satisfies the same equation. It is therefore clear that the notion of  minimum  can be replaced by that of strict minimum. Similarly, we can replace the notion of maximum by the one of strict maximum in the definition of subsolutions.   
}
\end{Remark}

%%%%%%%%%%%%%%%%%%%%%%%%%%%%%%%%%%%%%%%%%%%%%%%%%%%%%%%%%%%%%%%%%%%%%%%%%%%%%%%%%%%%%%%%%%%%%%%%%%%%%%%%%%
\subsection{Viscosity property}

We can now characterize $v$ as a continuous viscosity solution of   \reff{eq pde FK}. The continuity of $v$ follows from standard estimates and we omit the proof. 

\begin{Theorem}\label{thm visco FK} The value function  $v$   is  continuous on $[0,T]\x \R^d$ and is a viscosity solution on  $[0,T)\x \R^d$ of \reff{eq pde FK}.
\end{Theorem}

\proof We only prove the supersolution property of $v$. The proof of the subsolution property is symmetric. 
 Let $\vp \in C^{1,2}$ be such that   $0=\min_{[0,T]\x \R^d} (v-\vp)=(v-\vp)(\hat t,\hat x)$ for some $(\txh)\in [0,T)\x \R^d$. We proceed by contradiction, i.e. we assume that 
 	\b*
	\rho \vp(\txh)-\Lc^\Q \vp(\txh)<0\;
	\e*
and show that this contradicts \reff{eq: PD FK}. Indeed, if the above inequality holds at $(\txh)$, then 
	\b*
	\rho \vp(t,x)-\Lc^\Q \vp(t,x)\le 0\; 
	\e*
on a neighborhood  of $(\txh)$ of the form $B:= B_r(\hat t) \x B_r(\hat x)$, $r\in (\hat t,T-\hat t)$.  By Remark \ref{rem min strict}, we can then assume that there exists $\eta>0$ such   that
	\b*
	v\ge  \vp+\eta \;\;\mbox{ on } \partial_p B
	\e*
where $\partial_p B$ is the {parabolic boundary} of $B$, i.e. $(B_r(\hat t)\x \partial B_r(\hat x))\cup (\{\hat t+r\}\x {\rm cl} B_r(\xh))$. 
 
Let $\theta$ be the first exit time of $(t,X_{\txh}(t))_{t\ge \th}$ from  $B$.  
By   It\^{o}'s Lemma applied to $\vp$ and the above inequalities, we then obtain 
	\b*
	v(\txh)=\vp(\txh)
	&=&  
	\Esp{\beta_{\txh}(\theta)\vp(\theta,X_{\txh}(\theta))}
	\\
	&-& \Esp{\int_\th^{\theta} \beta_{\txh}(s) \left(\Lc^\Q\vp(s,X_{\txh}(s))-\rho\vp(s,X_{\txh}(s))\right) ds }
	\\
	&\le& 	
	\Esp{ \beta_{\txh}(\theta )\left(v(\theta ,X_{\txh}(\theta))-\eta \right)  }\;
	\\
	&<&\Esp{ \beta_{\txh}(\theta ) v(\theta ,X_{\txh}(\theta))  }\;,	
	\e*
a contradiction to \reff{eq: PD FK}.
\ep

%%%%%%%%%%%%%%%%%%%%%%%%%%%%%%%%%%%%%%%%%%%%%%%%%%%%%%%%%%%%%%%%%%%%%%%%%%%%%%%%%%%%%%%%%%%%%%%%%%%%%%%%%%
\subsection{Uniqueness} 

\subsubsection{An equivalent definition of viscosity solutions}

In order to complete the characterization of $v$, it remains to show that it is the unique solution of  \reff{eq pde FK} satisfying the boundary condition $v(T,\cdot)=g$. For this purpose, we need an alternative definition of viscosity solutions in terms  of super- et subjets. 

\vv

Note first that, if   $U$ is l.s.c., $\vp \in C^{1,2}$ and $(\hat t,\hat x)\in [0,T)\x \R^d$ is such that  $0=\min_{[0,T]\x \R^d} (U-\vp)=(U-\vp)(\hat t,\hat x)$, then a second order  Taylor expansion implies
	\b*
	U(t,x)&\ge& U(\txh)+\vp(t,x)-\vp(\txh)\\
	&=&
	U(\txh) + \partial_t \vp(\txh)(t-\th)\\
	&+&  ({x-\xh})'D \vp(\txh)+\frac12({x-\xh})'D^2 \vp(\txh)({x-\xh}) + o(|t-\hat t|+|x-\hat x|^2)\;. 
	\e*
This naturally leads to the notion of  {\sl subjet} defined as the set $\Pc^-U(\txh)$ of points $(q,p,A)\in \R\x\R^d\x\mathbb{S}^d$ satisfying 
	\b*
	U(t,x) \ge   
	U(\txh) + q(t-\th)+  ({x-\xh})'p +\frac12({x-\xh})'A(x-\xh)  + o(|t-\hat t|+|x-\hat x|^2)\;. 
	\e* 
We define similarly the {\sl superjet} $\Pc^+U(\txh)$ as the collection of points $(q,p,A)\in \R\x\R^d\x\mathbb{S}^d$ such that
	\b*
	U(t,x) \le   
	U(\txh) + q(t-\th)+  ({x-\xh})'p +\frac12({x-\xh})'A(x-\xh) + o(|t-\hat t|+|x-\hat x|^2)\;. 
	\e* 
For technical reasons related to Ishii's Lemma, see below,  we will also need to consider the ``limit'' super- and subjets. More precisely, we define $\bar \Pc^+ U(\txh)$ as the set of points $(q,p,A) \in \R\x\R^d\x\mathbb{S}^d$ for which there exists a sequence $ (t_n,x_n,q_n,p_n,A_n)_n$ of $\Dom$ such that $(t_n,x_n,q_n,p_n,A_n)\in \Pc^+ U(\txn)$ satisfying $(t_n,x_n,$ $U(t_n,x_n),$ $q_n,$ $p_n,A_n)$ $\to$ $(\txh,U(\txh),q,p,A)$. The set  $\bar \Pc^- U(\txh)$ is  defined similarly. 
\\ 
 
We can now state the alternative definition of viscosity solutions. 

\begin{Lemma}\label{lem def equi visco} Assume that $F$ is continuous.  A l.s.c. (resp. u.s.c.) function $U$ is  a {\sl supersolution} (resp. {\sl subsolution}) of   \reff{eq: def F FK avec I}  on $[0,T)\x \R^d$ if and only if for  all $(\hat t,\hat x)\in [0,T)\x \R^d$ and all  $(\hat q,\hat p, \hat A)\in \bar \Pc^-U(\txh)$   (resp. $\bar \Pc^+ U(\txh)$ ) 
	\be\label{eq F ge 0 def visco fonction teste bis}
	F(\hat t,\hat x, U(\hat t,\hat x),\hat q,\hat p, \hat A)\ge 0 \;\;\;\;(\mbox{ resp. $\le 0$})\;.
	\ee 
\end{Lemma}

\proof We only consider the supersolution property.   It is clear that the definition of the lemma implies the Definition   \ref{def solution viscosite}. Indeed, if $(\txh)\in \Dom$ is a minimum of $U-\vp$ then $(\partial_t \vp,D\vp,D^2\vp)(\txh) \in \bar \Pc^-U(\txh)$. It follows that
	\b*
	F(\hat t,\hat x, U(\hat t,\hat x),\hat q,\hat p, \hat A)\ge 0  
	\e*
with  $(\hat q,\hat p, \hat A)=(\partial_t \vp,D\vp,D^2\vp)(\txh)$. Since  $U\ge \vp$  and $F$ is elliptic, this implies the required result. 

We now prove the converse implication. Fix $(\hat t,\hat x)\in [0,T)\x \R^d$ and  $(\hat q,\hat p, \hat A)\in \bar \Pc^-U(\txh)$. It is clear that, if  $(\hat q,\hat p, \hat A)\in  \Pc^-U(\txh)$, then we can find $\vp$ locally $C^{1,2}$ such that   $(\hat q,\hat p, \hat A)=(\partial_t \vp,D\vp,D^2\vp)(\txh)$, $\vp=U$ at $(\txh)$ and $U\ge   \vp$ (see e.g.  \cite{FleSon93} page 225 for an example of construction). We then have 
	\b*
	F(\hat t,\hat x, U(\hat t,\hat x),\hat q,\hat p, \hat A )\ge 0  \;. 
	\e*
\ep

\subsubsection{Ishii's Lemma and Comparison Theorem}

The last ingredient to prove a comparison theorem is the so-called {Ishii's Lemma}. 

\begin{Lemma}{\bf (Ishii's Lemma)} Let $U$ (resp. $V$) be a  l.s.c. supersolution  (resp.  u.s.c. subsolution)  of \reff{eq: def F FK avec  I}  on $[0,T)\x \R^d$. Assume that $F$ is continuous and satisfies
	\b*
	 F(t,x,u,q,p,A)=F(t,x,u,0,p,A)-q \;\;   
	\e*
for all $(t,x,u,q,p,A)$.
Let $\phi \in C^{1,2,2}([0,T]\x\R^d\x\R^d)$ and $(\th,\xh,\hat y)\in [0,T)\x\R^d\x\R^d$ be such that 
	\b*
	W(t,x,y)&:=&V(t,x)-U(t,y)-\phi(t,x,y) \le W(\th,\xh,\hat y) \\
		&&\;\;\;\forall\;(t,x,y)\in [0,T)\x\R^d\x\R^d\;. 
	\e*
Then, for all $\eta>0$, there is $(q_1,p_1,A_1) \in \bar \Pc^+V(\txh)$ and $(q_2,p_2,A_2) \in \bar \Pc^-U(\th,\hat y)$ such that 
	\b*
	q_1-q_2=\partial_t \phi(\th,\xh,\hat y) & , & (p_1,p_2)=(D_x\phi,-D_y\phi)(\th,\xh,\hat y)
	\e*
and 
	\b*
	\left(\begin{array}{cc} A_1 & 0 \\ 0& -A_2 \end{array}\right) \le D_{(x,y)} \phi(\th,\xh,\hat y)  + \eta \left( D_{(x,y)} \phi(\th,\xh,\hat y)\right)^2\;. 
	\e*
\end{Lemma} 

\proof The proof is technical and long, we refer to \cite{CIL92} for details. 
\ep 

\vs2

We now prove the expected  {comparison theorem} also called {maximum principle}. 

\begin{Theorem}\label{thm comp FK}(Comparison) Let $U$ (resp. $V$) be a  l.s.c. supersolution  (resp.  u.s.c. subsolution)  with polynomial growth   of \reff{eq pde FK}  on $[0,T)\x \R^d$. If $U\ge V$ on $\{T\}\x \R^d$, then $U\ge V$ on $[0,T]\x \R^d$.
\end{Theorem}

\proof We can assume without loss of generality that $\rho>0$ (otherwise we replace $U$ and $V$ by $\tilde U(t,x):=e^{\kappa t}U(t,x)$ and  $\tilde V(t,x):=e^{\kappa t}V(t,x)$ for $\kappa$ large enough).  Assume now that there is some point $(t_0,x_0) \in  [0,T]\x \R^d$ such that  $U(t_0,x_0)<V(t_0,x_0)$. We shall prove that it leads to a contradiction. Let $\eps>0$, $\kappa>0$ and $p$ be an integer greater than $1$   such that  $\limsup_{|x|\to \infty} \sup_{t\le T} (|U(t,x)|+|V(t,x)|)/(1+|x|^p)=0$. Then there exists  $(\hat t ,\hat x ) \in  [0,T]\x \R^d$ such that  
	\b*
	0&<&V(\hat t,\hat x)-U(\hat t,\hat x)-\phi(\hat t,\hat x,\hat x)
	=
	\max_{(t,x)\in  [0,T]\x \R^d} \left(V(t,x)-U(t,x)-\phi(t,x,x)\right)\;, 
	\e*
where 
	\b*
	\phi(t,x,y):=\eps e^{-\kappa   t} (1+ | x|^{2p}+ | y|^{2p})\; 
	\e*
and $\eps$ is chosen small enough. 
Since  $U\ge V$ on $\{T\}\x\R^d$, it is clear that  $\hat t<T$.

For all  $n\ge 1$, we can  also  find $(t_n,x_n,y_n) \in [0,T]\x \R^d\x \R^d$ such that  
	\be\label{eq Gamma n >0 FK}
	0&<&\Gamma_n(t_n,x_n,y_n) 
	=
	\max_{(t,x,y)\in  [0,T]\x \R^d\x \R^d} \Gamma_n(t,x,y)
	\ee
where 
	\b*
	\Gamma_n(t,x,y)&:=  &
	V(t,x)-U(t,y)-\phi(t,x,y)-n|x-y|^2\\
	&-&(|t-\hat t|^2+ |x-\hat x|^4) \;. 
	\e*
It is easily checked that, after possibly passing to a subsequence, 
	\begin{equation}\label{eq comp FK visco conv tn xn yn}
	(t_n,x_n,y_n,\Gamma_n(t_n,x_n,y_n))\to (\hat t ,\hat x,\hat x,\Gamma_0(\hat t ,\hat x,\hat x)) \;\mbox{ and }\; n|x_n-y_n|^2 \to 0\;. 
	\end{equation}
Moreover, Ishii's Lemma implies that for all $\eta>0$, we can  find $(q^n_1,p^n_1,A^n_1) \in \bar \Pc^+V(\txn)$ and $(q^n_2,p^n_2,A^n_2) \in \bar \Pc^-U(\tn,y_n)$ such that 
	\b*
	q^n_1-q^n_2=\partial_t \vp_n(\tn,x_n,y_n) & , & (p_1,p_2)=(D_x\vp_n,-D_y\vp_n)(\tn,x_n,y_n)
	\e*
and 
	\b*
	\left(\begin{array}{cc} A^n_1 & 0 \\ 0& -A^n_2 \end{array}\right) \le D^2_{(x,y)} \vp_n(\tn,x_n,y_n)  + \eta \left( D^2_{(x,y)} \vp_n(\tn,x_n,y_n)\right)^2\;. 
	\e*
where 
	\b*
	\vp_n(t,x,y) := \phi(t,x,y)+n|x-y|^2+ |t-\hat t|^2+ |x-\hat x|^4\;.
	\e*
In order to obtain the required contradiction, it now suffices to appeal to  Lemma \ref{lem def equi visco} and to argue as in the proof of Theorem \ref{thm: comparaison regu FK}. Using \reff{eq comp FK visco conv tn xn yn}, we obtain that for all  $\eta>0$ 
	\b*
	\rho(V-U)(\hat t,\hat x)  
	&\le&\eps_n +\eta C_n +\Lc^\Q\phi(\hat t,\hat x,\hat x)
	\e*
where $\eps_n\to 0$ is independent of $\eta$ and $C_n$ does neither depend of $\eta$. By sending $\eta \to 0$, we deduce that
	\b*
	\rho(V-U)(\hat t,\hat x)  
	&\le&\eps_n    +\Lc^\Q\phi(\hat t,\hat x,\hat x)\;. 
	\e*
For $\kappa>0$ big enough so that the second  term in the right-hand side is strictly negative and $n$   large enough, we get $\rho(V-U)(\hat t,\hat x)  \le 0$. This contradicts the fact that 	$(V-U)(\hat t,\hat x)>0$ since $\rho$ is assumed to be (strictly) positive.
\ep

\begin{Corollary}\label{cor visco FK} The value function  $v$   is continuous and is the unique  viscosity solution on  $[0,T)\x \R^d$ of \reff{eq pde FK}
	satisfying $\lim_{s \uparrow T,\; y\to x} v(s,y)=g(x)$ in the class of discontinuous viscosity solutions with polynomial growth. 
\end{Corollary}

%%%%%%%%%%%%%%%%%%%%%%%%%%%%%%%%%%%%%%%%%%%%%%%%%%%%%%%%%%%%%%%%%%%%%%%%%%%%%%%%%%%%%%%%%%%%%%%%%%%%%%%%%%
%%%%%%%%%%%%%%%%%%%%%%%%%%%%%%%%%%%%%%%%%%%%%%%%%%%%%%%%%%%%%%%%%%%%%%%%%%%%%%%%%%%%%%%%%%%%%%%%%%%%%%%%%%
\chapter{The pricing equation II: the incomplete market case}\label{CHAP: pricing equation II}

In this section, we provide the pricing equation under portfolio constraints as studied in Section \ref{Section: chap dualite incomplete market} of Chapter \ref{CHAP: prix dual}. \\

We keep the notations and assumptions on the coefficients  of Chapter \ref{CHAP: pricing equation complet} except that we now assume that $\rho$ is bounded.

We define the value function:
$$
(t,x)\mapsto v(t,x):=\sup_{\nu \in \Uc_b} J(t,x;\nu)
$$
where 
$$
 J(t,x;\nu):=\E^{\Q^{\nu}_{t,x}}\left[\beta_{t,x}(T)g(X_{t,x}(T))- \int_t^T \beta_{t,x}(s) \delta_K(\nu_s) ds\right]
$$
and 
$$
 \frac{d\Q^{\nu}_{t,x}}{d\P}=\Ec^\nu_{t,x}(T)
$$
with 
\b*
\Ec^\nu_{t,x}(s)&=&e^{-\frac12 \int_t^s | \lambda^\nu_{t,x}(u)|^2 du - \int_t^s  \lambda^\nu_{t,x}(u) dW_u}\\
\lambda^\nu_{t,x}&:=&\sigma(X_{t,x})^{-1}(\mu(X_{t,x})-\rho(X_{t,x})X_{t,x}- \nu)\;.
\e*

Since $\rho$ is bounded, we have that $\beta_{t,x}$ and $\beta_{t,x}^{-1}$ are  bounded. By replacing $\nu \in \Uc_b$ by $\beta_{t,x} \nu \in \Uc_b$ and vice-versa, we deduce from  Section \ref{Section: chap dualite incomplete market} of Chapter \ref{CHAP: prix dual} that 
$$
p_K(G)=v(0,X_0)\;. 
$$

\begin{Remark}\label{rem: v prix surrep incomplet est sci}{\rm One easily checks that $J(\cdot;\nu)$ is l.s.c. for each $\nu \in \Uc_b$. It follows that $v$ is l.s.c. as well.}
\end{Remark}
 
%%%%%%%%%%%%%%%%%%%%%%%%%%%%%%%%%%%%%%%%%%%%%%%%%%%%%%%%%%%%%%%%%%%%%%%%%%%%%%%%%%%%%%%%%%%%%%%%%%%%%%%%%%
\section{Dynamic programming principle}

The key result for the derivation of a PDE associated to $v$ is the so-called dynamic programming principle. In the following, we denote by $\Tc^t_{[t,T]}$ the set of elements of $\Tc_{[t,T]}$ that are independent on $\Fc_t$. 

\begin{Theorem}\label{thm: dynamic programming pour pde} Fix $(t,x)\in [0,T)\x \R^d$ and let $\{\theta^\nu,\nu\in \Uc_b\}\subset  \Tc^t_{[t,T]}$ be such that $X_{t,x}$ is essentially bounded on $[t,\theta^\nu]$ for each $\nu\in \Uc_b$. Then, 
\b*
v(t,x)=\sup_{\nu \in \Uc_b} \E^{\Q_{t,x}^\nu}\left[\beta_{t,x}(\theta^\nu)v(\theta^\nu,X_{t,x}(\theta^\nu)- \int_t^{\theta^\nu} \beta_{t,x}(s) \delta_{K}(\nu_s) ds\right]\;.
\e*
\end{Theorem}

\proof   For ease of notations, we omit, if not necessary, the dependence of $\theta$ with respect to $\nu$. Let $\bar v(t,x)$ denote the right-hand side term in the above equation. 
 We first show that $v(t,x)\le \bar v(t,x)$. 
To see this, observe that
\b*
J(t,x;\nu)
&=&
 \E^{\Q_{t,x}^\nu}\left[
  \beta_{t,x}(\theta)
 \E^{\Q_{t,x}^\nu}\left[ 
 \beta_{\theta,\zeta}(T)
g(X_{\theta,\zeta}(T))- \int_{\theta}^T \beta_{\theta,\zeta}(s) \delta_{K}(\nu_s) ds
 ~|~\Fc_{\theta}\right]\right.\\
 &&~~~~~~~~~~~~~~~~~~~~~~~~\left.
 - \int_t^{\theta} \beta_{t,x}(s) \delta_{K}(\nu_s) ds
 \right]
\e*
where $\zeta:=X_{t,x}(\theta)$. We now observe that 
\b*
 &&\E^{\Q_{t,x}^\nu}\left[ 
 \beta_{\theta,\zeta}(T)
g(X_{\theta,\zeta}(T))- \int_{\theta}^T \beta_{\theta,\zeta}(s) \delta_{K}(\nu_s) ds
 ~|~\Fc_{\theta}\right]\\
 &&= \E^{\Q_{\theta,\zeta}^\nu}\left[ 
 \beta_{\theta,\zeta}(T)
g(X_{\theta,\zeta}(T))- \int_{\theta}^T \beta_{\theta,\zeta}(s) \delta_{K}(\nu_s) ds
 ~|~\Fc_{\theta}\right]
 \le v(\theta,\zeta)\,
\e*
Hence, the fact that   $v\le \bar v$. 

We now prove the converse inequality. To this purpose, given $k\ge 1$, we denote by $\Uc_{bk}$ the set of elements $\nu \in \Uc_b$ such that  $\sup_{s\le T}(|\nu_s|\p |\delta_K(\nu_s)|)\le k$ $\Pas$  and set $v_k(s,y):=\sup_{\nu \in \Uc_{bk}}J(s,y;\nu)$. Note that $v_k\uparrow v$ as $k\to \infty$. Then, for fixed $k\ge 1$, one easily checks that $J(\cdot;\nu)$ and $v_k$ are locally uniformly continuous in $(t,x)$, uniformly in $\nu \in \Uc_{bk}$.  

Let $\Uc_{bk}^t$ denotes the set of elements of $\Uc_{bk}$ that are independent of $\Fc_t$. Then, one easily checks, by using the fact that $X$ solves a Brownian SDE, that $v_k(t,x)=\sup_{\nu \in \Uc_{bk}^t} J(t,x;\nu)$, see \cite{BT10}.

Fix $\eps>0$. For $(s,y)\in [0,T]\x \R^d$, we can find $\nu^{s,y}\in \Uc^s_{bk}$ such that 
\be\label{eq: proof ppd 1}
 J(s,y;\nu^{s,y})\ge v(s,y)-\eps\;. 
 \ee
 
Let $A\subset \R^d$ be a compact set such that $X_{t,x}$ takes values in $A$ on $[t,\theta^\nu]$ for each $\nu\in \Uc_b$. 
It follows from the local uniform  continuity of $J$ and $v_k$ that there exists $\eta>0$ and a finite collection of points $(t_i,x_i)_{i\le I}\in [0,T]\x A$ such that $\cup_{i\le I} [t_i-\eta,t_i]\x B_\eta(x_i)\supset A$.  
\be\label{eq: proof ppd 2}
 |J(\cdot;\nu^{t_i,x_i})-J(t_i,x_i;\nu^{t_i,x_i}) \p |v-v(t_i,x_i)| \le \eps\;\; \mbox{ on } [t_i-\eta,t_i]\x B_\eta(x_i)\;.
 \ee
 Combining \reff{eq: proof ppd 1} and \reff{eq: proof ppd 2} leads to 
\be\label{eq: proof ppd 3}
 J(\cdot;\nu^{t_i,x_i})\ge  v-3\eps\;\; \mbox{ on } [t_i-\eta,t_i]\x B_\eta(x_i)\supset A_i\;,
 \ee
 where the $A_i$ can be constructed in such a way that they form a partition of $A$.

Given $\nu \in \Uc_{bk}$, we now define 
$$
\bar \nu:=\nu \1_{[0,\theta)}\p \1_{[\theta,T]} \sum_{i\le I} \nu^{t_i,x_i} \1_{(\theta,X_{t,x}(\theta) \in A_i}
$$
Then, using the fact that $\nu^{s,y}$ is independent of $\Fc_s$, for all $(s,y) \in [0,T]\x \R^d$, we obtain 
\b*
J(t,x;\bar \nu)&=&  \E^{\Q_{t,x}^\nu}\left[
  \beta_{t,x}(\theta)
 \E^{\Q_{t,x}^{\bar \nu}}\left[ 
 \beta_{\theta,\zeta}(T)
g(X_{\theta,\zeta}(T))- \int_{\theta}^T \beta_{\theta,\zeta}(s) \delta_{K}(\bar \nu_s) ds
 ~|~\Fc_{\theta}\right]\right.\\
 &&~~~~~~~~~~~~~~~~~~~~~~~~\left.
 - \int_t^{\theta} \beta_{t,x}(s) \delta_{K}(\nu_s) ds
 \right]
 \\
 &=&
  \E^{\Q_{t,x}^\nu}\left[
  \beta_{t,x}(\theta) \left(\sum_{i\le I} J(\theta,\zeta;\nu^{t_i,x_i})\1_{(\theta,\zeta)\in A_i}\right) - 
  \int_t^{\theta} \beta_{t,x}(s) \delta_{K}(\nu_s) ds
 \right]
\e*
so that, by  \reff{eq: proof ppd 3},
\b*
v_k(t,x)&\ge& J(t,x;\bar \nu) 
\\
&\ge &
  \E^{\Q_{t,x}^\nu}\left[
  \beta_{t,x}(\theta) v_k(\theta,\zeta) - 
  \int_t^{\theta} \beta_{t,x}(s) \delta_{K}(\nu_s) ds
 \right]-3\eps \E^{\Q_{t,x}^\nu}[
  \beta_{t,x}(\theta)]\;.
\e* 
Sending $\eps \to 0$ and using the arbitrariness of $\nu$, then shows that 
\b*
v_k(t,x)&\ge&\sup_{\nu \in \Uc_{bk}}
  \E^{\Q_{t,x}^\nu}\left[
  \beta_{t,x}(\theta) v_k(\theta,\zeta) - 
  \int_t^{\theta} \beta_{t,x}(s) \delta_{K}(\nu_s) ds
 \right]\;.
\e*
The result then follows by sending $k\to \infty$ and by using the monotone convergence theorem. 
\ep 

\begin{Remark}{\rm 
In the above proof, we used the approximation   $v_k$ in order to reduce to the case where the value function is u.s.c.  A more direct approach, based on test functions, could also be adopted, see \cite{BT10}. In particular, it would allow to provide a weak version of the dynamic programming principle of Theorem \ref{thm: dynamic programming pour pde} even if $v$ was not known to be measurable a-priori. It would then take the form: 
\b*
v(t,x)\le \sup_{\nu \in \Uc_b} \E^{\Q_{t,x}^\nu}\left[\beta_{t,x}(\theta^\nu)v^*(\theta^\nu,X_{t,x}(\theta^\nu)- \int_t^{\theta^\nu} \beta_{t,x}(s) \delta_{K}(\nu_s) ds\right]\;,
\e*
and 
\b*
v(t,x)\ge \sup_{\nu \in \Uc_b} \E^{\Q_{t,x}^\nu}\left[\beta_{t,x}(\theta^\nu)v_*(\theta^\nu,X_{t,x}(\theta^\nu)- \int_t^{\theta^\nu} \beta_{t,x}(s) \delta_{K}(\nu_s) ds\right]\;,
\e*
for all family of stopping times $\{\theta^\nu,\nu\in \Uc_b\}\subset  \Tc^t_{[t,T]}$   such that $X_{t,x}$ is essentially bounded on $[t,\theta^\nu]$, for each $\nu\in \Uc_b$. 

In the above assertion, $v^*$ could be replaced by $v$ if it is known to be measurable, in particular if it is l.s.c. Hence, for $v=v_*$, it coincides with the formulation of Theorem \ref{thm: dynamic programming pour pde}.

In view of the arguments below, the later formulation would already be enough to provide a PDE characterization for $v_*$ and $v^*$.  
}
\end{Remark}

\section{Hamilton-Jacobi-Bellman pricing equation}

In this section, we use the dynamic programming principle of Theorem \ref{thm: dynamic programming pour pde} to show that $v$ is a (discontinuous) viscosity solution of 
\be\label{eq: HJB marche incomplet}
\min\left\{\rho v- \Lc^\Q v \;,\;\min_{|\zeta|=1} \delta_K(\zeta)-\zeta'Dv\right\}=0\mbox{ on } [0,T)\x \R^d\;,
\ee
and provide a suitable boundary condition at $t=T$, which is related to the face-lifting phenomenon observed in Section \ref{Section: chap dualite incomplete market} of   Chapter \ref{CHAP: prix dual}.

\subsection{PDE characterization in the domain}

We first discuss the supersolution property. 

\begin{Proposition}\label{prop: sursol domaine v marche incomplet} The function $v$ is a viscosity supersolution of \reff{eq: HJB marche incomplet}.
\end{Proposition}

\proof Fix $\nu=u$ for some $u\in \R^d$ such that $\delta_K(u)<\infty$. Then, it follows from Theorem \ref{thm: dynamic programming pour pde} that 
\b*
v(t,x)&\ge &\E^{\Q_{t,x}^\nu}\left[\beta_{t,x}(\theta^h)v(\theta^h,X_{t,x}(\theta^h))- \int_t^{\theta^h} \beta_{t,x}(s) \delta_{K}(u) ds\right]
\\
&=& 
\E\left[\Ec_{t,x}^\nu(\theta^h)\left(\beta_{t,x}(\theta^h)v(\theta^h,X_{t,x}(\theta^h))- \int_t^{\theta^h} \beta_{t,x}(s) \delta_{K}(u) ds\right)\right]
\;,
\e*
where $\theta^h:=\inf\{s\ge t~:~|X_{t,x}(s)-x|\p |\Ec^\nu_{t,x}(s)-1|\ge 1\}\wedge (t\p h)$.
Let $\vp$ be a smooth function such that $(t,x)$ achieves a minimum of $v-\vp$, recall Remark \ref{rem: v prix surrep incomplet est sci}. We can always assume that $(v-\vp)(t,x)=0$. Thus, 
\b*
\vp(t,x)&\ge &
\E\left[\Ec_{t,x}^\nu(\theta^h)\left(\beta_{t,x}(\theta^h)\vp(\theta^h,X_{t,x}(\theta^h))- \int_t^{\theta^h} \beta_{t,x}(s) \delta_{K}(u) ds\right)\right]
\;.
\e*
By following  the same arguments as in the proof Theorem \ref{thm FK regu} and using the arbitrariness of $u$, we deduce that $\vp$ satisfies:
\b*
(\rho \vp- \Lc^\Q \vp)(t,x) \p \delta_K(u)-u'D\vp(t,x)\ge 0\;.
\e*
Since $u$ is abritrary and the set $\{u\in \R^d~:~\delta_K(u)<\infty\}$ is a cone which contains $0$, this proves the required result.
\ep

\begin{Proposition}\label{prop: soussol domaine v marche incomplet} The function $v^*$ is a viscosity subsolution of \reff{eq: HJB marche incomplet}.
\end{Proposition}

\proof Let $\vp$ be a smooth function such that $(t,x)$ achieves a strict local maximum of $v^*-\vp$. We can always assume that $(v-\vp)(t,x)=0$. We argue by contradiction and assume that 
 \b*
\min\left\{\rho \vp- \Lc^\Q \vp \;,\;\min_{|\zeta|=1} \delta_K(\zeta)-\zeta'D\vp\right\}(t,x)>0\;.
\e*
Then, 
 \be\label{eq: contra sous sol hjb marche incomplet}
\min\left\{\rho \vp- \Lc^\Q \vp \;,\;\min_{|\zeta|=1} \delta_K(\zeta)-\zeta'D\vp\right\}>0\; \mbox{ on } B_\eta(t,x)
\ee
for some $\eta>0$ small enough. 
Let $(t_n,x_n)_n$ be a sequence in $B_\eta(t,x)$ that converges to $(t,x)$ and such that $v(t_n,x_n)\to v^*(t,x)$. Let $\theta_n$ be the first exit time of $B_\eta(t,x)$ by $(s,X_{t_n,x_n}(s))_{s\ge t_n}$.  Fix $\nu \in \Uc_{b}$.  Using It\^o's Lemma and \reff{eq: contra sous sol hjb marche incomplet}, we then deduce that 
\b*
\vp(t_n,x_n)\ge \E^{\Q_{t_n,x_n}^{\nu}}\left[\beta_{t_n,x_n}(\theta_n)\vp(\theta_n,X_{t_n,x_n}(\theta_n))- \int_{t_n}^{\theta_n} \beta_{t_n,x_n}(s) \delta_{K}(\nu_s) ds\right]\;.
\e*
Moreover, since $(t,x)$ achieves a strict local maximum of $v^*-\vp$, we have $v-\vp\le v^*-\vp\le  -\xi$ on $\partial_p B_{\eta}(t,x)$ for some $\xi>0$. Hence, 
\b*
\vp(t_n,x_n)&\ge& \E^{\Q_{t_n,x_n}^\nu}\left[\beta_{t_n,x_n}(\theta_n)v(\theta_n,X_{t_n,x_n}(\theta_n))- \int_{t_n}^{\theta_n} \beta_{t_n,x_n}(s) \delta_{K}(\nu_s) ds\right]
\\
&&\p \xi \E^{\Q_{t_n,x_n}^\nu}[\beta_{t_n,x_n}(\theta_n)]\;.
\e*
Since $\rho$ is bounded, one easily checks that $\E^{\Q_{t_n,x_n}^\nu}[\beta_{t_n,x_n}(\theta_n)]\ge c$ for some $c>0$, for all $n$ and $\nu \in \Uc_{b}$. 
We then obtain 
\b*
v(t_n,x_n)&\ge& \E^{\Q_{t_n,x_n}^\nu}\left[\beta_{t_n,x_n}(\theta_n)v(\theta_n,X_{t_n,x_n}(\theta_n))- \int_{t_n}^{\theta_n} \beta_{t_n,x_n}(s) \delta_{K}(\nu_s) ds\right]
\\
&&\p \xi c \p v(t_n,x_n)-\vp(t_n,x_n)\;.
\e*
Since $ v(t_n,x_n)-\vp(t_n,x_n)\to 0$ as $n\to \infty$, we obtain a contradiction to Theorem \ref{thm: dynamic programming pour pde} for $n$ large enough.
\ep

\subsection{Boundary condition at $t=T$}

In order to complete the characterization of $v$, it remains to provide a terminal condition. We shall show below that $$v(T-,\cdot)=\hat g$$ where $\hat g$ is defined as in Chapter \ref{CHAP: prix dual}:
$$
\hat g(x):=\sup_{\zeta \in \R^d} g(x\p \zeta)-\delta_K(\zeta)\;.
$$

We split the proof in two separate propositions.   
\begin{Proposition}  For all $x\in \R^d$, $v(T,x)\ge \hat g(x)$. 

\end{Proposition}

\proof Let $(t_n,x_n)_{n\ge1}$ be a sequence such that $(t_n,x_n) \rightarrow (T,x)$ and $v (t_n,x_n) \rightarrow v (T,x)$. By the definition of $v$, we have

\begin{equation*}
    v (t_n,x_n) \ge \E^{\Q^{\nu^n}_{t_n,x_n}} \left[ \beta^n (T) g \left( X^n (T) \right) - \int_{t_n}^T \beta^n (s) \delta_K \left( \nu^n_s \right) ds \right]
\end{equation*}

where $\left( \beta^n,X^n \right) := \left( \beta_{t_n,x_n} , X_{t_n,x_n} \right)$ and $\nu^n_s := \frac{1}{T - t_n} u$, for some $u \in \text{dom} \left(  \delta_K \right)$. 
Now observe that $\delta_K (\lambda u) = \lambda \delta_K (u)$ for every $\lambda>0$, so that 
\begin{equation*}
    \int_{t_n}^T \beta^n (s)\delta_K \left( \nu^n_s \right) ds = \delta_K (u) \frac{1}{T - t_n} \int_{t_n}^T \beta^n (s) ds \underset{n \rightarrow \infty}{\longrightarrow} \delta_K (u) \;\Pas
\end{equation*}
since $\rho$ is bounded. Hence, using the fact that $\rho$ is bounded again and a dominated convergence argument, we obtain 
\b*
 v(T,x)&= & \lim_{n\to \infty} v \left( t_n,x_n \right) \\
 &\ge&  \liminf_{n\to \infty} \E^{\Q^{\nu^n}_{t_n,x_n}} \left[ \beta^n (T) g \left( X^n (T) \right)  \right]- \delta_K (u) .
\e*
To conclude the proof, it remains show that $E_n:=\E^{\Q^{\nu^n}_{t_n,x_n}} \left[ \beta^n (T) g \left( X^n (T) \right)  \right]\to g(x\p u)$, and use the arbitrariness of $u\in$ dom$(\delta_K)$. 

To see this, first observe that  
\begin{equation*}
    X^n   = u\frac{\cdot-t_n}{T-t_n} + x_n + \int_{t_n}^ \cdot \rho \left( X^n (s) \right) X^n (s) ds + \int_{t_n}^ \cdot \sigma \left( X^n (s) \right) dW^{\nu^n}_s,
\end{equation*}
so that 
\begin{equation*}
    E_n = \Esp{\beta^n (T) g \left( Z^n (T) \right)},
\end{equation*}
where $Z^n$ satisfies
\begin{equation*}
    Z^n   = u\frac{\cdot-t_n}{T-t_n}  \p x_n + \int_{t_n}^\cdot \rho \left( Z^n (s) \right) Z^n (s) ds + \int_{t_n}^ \cdot \sigma \left( Z^n (s) \right) dW_s.
\end{equation*}
 
Clearly, the sequence $(Z^n (T))_{n\ge 1} $ is bounded in $L^2$ and converges to $x\p u$ $\Pas$ Since $g$ is continuous with linear growth, the dominated convergence theorem implies 
\begin{equation*}
   \lim_{n\to \infty} E_n=g(x\p u).
\end{equation*}
\ep

\begin{Proposition} Assume that $\hat g$ is upper-semicontinuous with linear growth. Assume further that $\sigma$ is bounded. Then, for all $x\in \R^d$, $v^*(T,x)\le \hat g(x)$. 
\end{Proposition}

\proof  Let $(t_n,x_n)_n$ be a sequence which
converges to $(T,x_0)$ and such that $v(t_n,x_n)\to v^*(T,x_0)$.
Set $(\beta^n, X^n) $ $=$ $(\beta_{t_n,x_n}, X_{t_n,x_n})$. By definition of $v$, there is some
$\nu^n\in \Uc_b$ such that
    \b*
    v(t_n,x_n)
    &\le &
    \E^{\Q^{\nu^n}}\left[\beta^n(T) g(X^n(T))-\int_{t_n}^T\beta^n(s)\delta_K(\nu^n_s) ds \right]+n^{-1}\;.
    \e*
Since dom$(\delta_K)$ is a convex cone and  $\delta_K$ is $1$-homogeneous, we have
    \b*
      \beta^n(T)g(X^n(T))
     &\le&
    \beta^n(T) \hat g\left(X^n(T) -\int_{t_n}^{T}\beta^n(T)^{-1}\beta^n(s)  \nu^n_sds  \right)
    \\
    &&\p \int_{t_n}^{T}  \beta^n(s)\delta_K(\nu^n_s)ds\;.
    \e*
  This implies that 
      \b*
    v(t_n,x_n)
    &\le &
    \E^{\Q^{\nu^n}} \left[\beta^n(T) \hat g\left(X^n(T) -\int_{t_n}^{T}\beta^n(T)^{-1}\beta^n(s)  \nu^n_sds \right) \right]+n^{-1}\;.
    \e*
In view of the above inequalities and the definition of
$(t_n,x_n)$, it remains to show that
    \be\label{eq prop cond T pour v^* 1}
    \limsup_{n\to \infty}
     \E^{\Q^{\nu^n}} \left[\beta^n(T) \hat g\left(X^n(T) -\int_{t_n}^{T}\beta^n(T)^{-1}\beta^n(s)  \nu^n_sds\right) \right] &\le& \hat g(T,x_0)\;.
    \ee
From now on, we assume that $\hat g$ is uniformly Lipschitz continuous. We shall explain at the end of the proof how to handle cases where it is not true. 
If $\hat g$ is $L$-Lipschitz, then
\b*
&&\left|\beta^n(T) \hat g\left(X^n(T) -\int_{t_n}^{T}\beta^n(T)^{-1}\beta^n(s)  \nu^n_sds \right)-\beta^n(T)\hat g(x_0)\right|
\\
&&\le L\left|  \beta^n(T) X^n(T) -\int_{t_n}^{T}\beta^n(s)  \nu^n_sds-\beta^n(T)x_0\right|
\\
&&=L\left| \int_{t_n}^T \beta^n(s) \sigma(X^n(s)) dW^{\nu^n}_s \p x_n-\beta^n(T)x_0 \right|
\e*
where,  since $\sigma$ and $\rho$ are bounded,  
$$
\E^{\Q^{\nu^n}} \left[\left| \int_{t_n}^T \beta^n(s) \sigma(X^n(s)) dW^{\nu^n}_s  \right|\right]\le C(T-t_n)^\frac12\;
$$
for some $C>0$ independent of $n$.
This proves the required result for $\hat g$ Lipschitz. 

In the case, we $\hat g$ is not Lipschitz, then we   construct,  for each $\eps>0$,   a
Lipschitz function $\Psi_\eps$ such that $|\hat
g(x_0)-\Psi_\eps(x_0)| \le \eps$ and $\Psi_\eps\ge \hat g$.
It follows that, for each $\eps$, we can find some finite
$L_\eps>0$ such that
    \b*
   && \limsup_{n\to \infty}
    \E^{\Q^{\nu_n}}\left[ \beta^n(T)\hat g\left( X^n(T) -\int_{t_n}^{T}\beta^n(T)^{-1}\beta^n(s)  \nu^n_sds\right)\right]
    \\
    &&\le
    \limsup_{n\to \infty}
    \E^{\Q^{\nu_n}}\left[  \beta^n(T)\Psi_\eps\left( X^n(T) -\int_{t_n}^{T}\beta^n(T)^{-1}\beta^n(s)  \nu^n_sds \right)\right]
    \\
    &&\le
    \Psi_\eps(x_0) + \limsup_{n\to \infty} L_\eps \;C\; \left( |x_n-x_0|  +  (T-t_n)^{1/2}\right)
    \\
    &&=\Psi_\eps(x_0)\le \hat g(x_0)\p \eps\;
    \e*
and  the proof   is concluded by sending $\eps$ to $0$.
 
  We conclude this proof by constructing the sequence
of functions $(\Psi_\eps)_{\eps>0}$.  For $x
\in \R^d$, we   define
    \b*
    G_k( x)&=&
    \sup_{z \in \R^d}
    \left[ \hat g( z)-k  |z -x | \right]
    \;\;,\;k\ge 1\;.
    \e*
Recall that $g$ has linear growth. Clearly, $G_k \ge \hat g$ and $G_k$ is $k$-Lipschitz. Moreover,
taking $k$ large enough, it follows from the linear growth and
upper-semicontinuity assumptions on $\hat g$ that,  for
all $x\in \R^d$, the maximum is attained in the
above definition by some $x_k(x)$. In particular,
    \b*
    G_k(x)&=& \hat g(x_k(x))
    -k |x_k(x)-x| \;\ge\;\hat g( x )
    \;.
    \e*
Using the linear growth of $\hat g$ again, we deduce  that
$ x_k(x) \to x$ as $k\to \infty$ after possibly
passing to a subsequence. Since $\hat g$ is upper-semicontinuous,
this also implies that
  \b*
  \hat g(x_0) &\ge&
  \limsup_{k\to \infty} \hat g( x_k(x_0))
  \sge
  \limsup_{k\to \infty} G_k(x_0)
   \sge \hat g(  x_0 ) \;.
    \e*
We can then choose $k_\eps$ such that  $|G_{k_\eps}( x_0)-\hat
g( x_0)|\le \eps$ and set $\Psi_\eps:=G_{k_\eps}$. 
\ep

\section{Example: non-hedgeable stochastic volatilty}

As an example of application, let us come back to the model of Example \ref{ex: dual stochastic vol} of Chapter \ref{CHAP: prix dual}. Note here that the volatility of $X^1$ is given by $X^1\sigma(X^2)$ which is not bounded. However, the above argument holds. Moreover, $X^1$ takes values in $(0,\infty)$ but it does not change anything in the above proofs. 

In this case, the price function $v$ is therefore a (discontinuous) viscosity solution on $[0,T)\x (0,\infty)\x \R$ of 
 \be\label{eq: pde vol sto}
\min\left\{r \vp- \Lc^\Q\vp\;,\;\min_{|\zeta|=1} \delta_K(\zeta)-\zeta'D\vp\right\}(t,x)=0\;
\ee
where 
$$
\Lc^\Q\vp =\partial_t \vp\p rx^1 \partial_{x^1} \vp \p\frac12 \left[ (x^1\sigma(x^2))^2  \partial^2_{x^1x^1} \vp\p \gamma^2  \partial^2_{x^2x^2} \vp \p 2 x^1\sigma(x^2)\gamma_1 \partial^2_{x^1x^2} \vp\right] 
$$
with $\gamma^2:=\gamma_1^2\p \gamma_2^2$. 
Moreover, it satisfies $v_*(T,x^1,x^2)=v^*(T,x^1,x^2)=\hat g(x^1)$.

\vs2

Since $\delta_K(\zeta)=0$ is $\zeta^1=0$ and $\delta_K(\zeta)=\infty$ is $\zeta^1\ne 0$, we deduce that $$v_*(T,x^1,x^2)=v^*(T,x^1,x^2)=\hat g(x^1)=g(x)$$ and, using the right-hand side term in \reff{eq: pde vol sto}, that $v$ is a supersolution of $\partial_{x^2}\vp=0$ and $-\partial_{x^2}\vp=0$ on $[0,T)\x (0,\infty)\x \R$. As for smooth functions, this implies that $v$ does not depend on $x^2$. We therefore now simply write $v(t,x^1)$. As for smooth function again, this also implies that $v$ is a supersolution on $[0,T)\x (0,\infty)$ of 
\be\label{eq: before BS Barenblatt}
\inf_{x^2} \Hc_{x^2}=0 \;
\ee
where 
$$
 \Hc_{x^2}\vp:=r \vp- \partial_t \vp- rx^1 \partial_{x^1} \vp - \frac12   (x^1\sigma(x^2))^2  \partial^2_{x^1x^1}\vp\;,
$$
i.e. 
\be\label{eq: BS Barenblatt}
r \vp- \partial_t \vp-rx^1 \partial_{x^1} \vp - \frac12 (x^1)^2 \left[ \underline \sigma^2\1_{  \partial^2_{x^1x^1}\vp < 0} \p \bar \sigma^2\1_{  \partial^2_{x^1x^1}\vp\ge  0} \right]\partial^2_{x^1x^1}\vp =0\;,
\ee
where $\bar \sigma:=\sup_{x^2}\sigma(x^2)$ and $\underline \sigma:=\inf_{x^2} \sigma(x^2)$. This is the so-called Black-Scholes-Barenblatt equation. 

\vs2

When $\bar \sigma<\infty$, $\sigma$ is continuous  and $g$ is continuous with linear growth, it is possible to show that this equation admits a comparison principle in the class of functions with linear growth. 
In particular, if there exists a smooth solution, say $\vp $, satisfying $\vp(T-,\cdot)=g$, then $v\ge \vp $. But on the other hand, \reff{eq: BS Barenblatt} and the previous boundary condition imply  that 
\b*
 \vp(0,X^1_0)\p \int_0^T \beta_s D \vp(s,X^1_s)dX_s^1
 &=& \beta_T \vp(T,X_T^1)\p \int_0^T \beta_s \Hc_{X^2_s} \vp(s,X^1_s)ds 
 \\
 &\ge& \beta_T \vp(T,X_T^1)\\
 & =& \beta_T g(X_T^1)\;,
\e*
where $X:=X_{0,X_0}$ and $\beta=\beta_{0,X_0}$. This shows that $v=\vp$. 

\vs2

In the limiting case where $\bar \sigma=\infty$, then \reff{eq: before BS Barenblatt} implies that $v$ is concave in $x^1$. If moreover, $\underline \sigma=0$ and $r=0$, then it should be non-increasing in time. This implies that $v\ge \bar g$, where $\bar g$ denotes the concave envelope of $g$. On the other hand, it is clear that, starting with $\bar g(X_0)$ allows to find a super-hedging strategy, which is actually of buy-and-hold type. Hence, $v=\bar g$. 
Note that $\bar g(x^1)=x^1$ for  $g(x^1)=[x^1-\kappa]^\p$ ! The same holds for $r\ne 0$ up to passing to discounted quantities. 

%%%%%%%%%%%%%%%%%%%%%%%%%%%%%%%%%%%%%%%%%%%%%%%%%%%%%%%%%%%%%%%%%%%%%%%%%%%%%%%%%%%%%%%%%%%%%%%%%%%%%%%%%%%%%%%%%%%%%%%%%%%%%%%%%%%%%%%%%%%%%%%%%%%%%%%%%%%%%%%%%%%%%%%%%%%%%%%%%%%%%%%%%%%%%%%%%%%%%%%%%%%%%%%%%%%%
\chapter{Approximate hedging and risk control}\label{CHAP: quantile et loss par dualite}

In this section, we discuss two approximate hedging technics that were discussed in \cite{FoLe99} and \cite{FoLe00}. We shall restrict here to the case of complete markets without constraints, because it is {\sl essentially} the only case where explicit formulations can be obtained by standard convex duality technics, and it already provides the general form of the solution. Extensions to incomplete markets are considered in the above mentioned papers. More general models will be discussed in Chapter \ref{CHAP: cible en proba}, in a Markovian setting.

\section{Quantile hedging}\label{Section: quantile hedging dualite}

We first discuss the case of a trader who wants to  hedge a random payoff $G\in L^0(\R_+)\setminus\{0\}$ from an initial wealth $y>0$. However, because for instance the hedging price  $p(G)=\EspQ{\beta_TG}$ is too high (which can be due to the fact that it was face-lifted in order to avoid explosion of the hedging strategy near the maturity, see Section \ref{Section: chap dualite incomplete market} of Chapter \ref{Section: chap prix dual complet}), his initial wealth is strictly less than $p(G)$. 

\subsection{Minimizing the probability of missing the hedge}

The first criteria we discuss here is the so-called {\sl quantile hedging} criteria. Namely, we try to find the optimal solution to the problem
	\be\label{eq: prob quantile}
	\inf_{\phi \in \Ac_+(y)} \Pro{G>Y^{y,\phi}_T} \mbox{ for some } 0<y<p(G),
	\ee
where $\Ac_+(y)$ is the restriction  of  $\Ac_b$ to strategies leading to non-negative wealth processes. 
\vs2

As shown in \cite{FoLe99}, this problem can be reduced to a standard test problem in mathematical statistics, which can then be solved by using the Neyman and Pearson's Lemma which we recall below.

To see this, we first note that the problem \reff{eq: prob quantile} can be reduced as follows. 
\begin{Proposition} The following holds:
 	\be\label{eq: prob quantile bis}
	\sup_{\phi \in \Ac_+(y)} \Pro{Y^{y,\phi}_T\ge G}= \sup\left\{ \Esp{\vp}\;,\;\vp \in L^0(\{0,1\})\mbox{ s.t. } \E^\Q[\beta_TG \vp]\le y\right\}\;.
	\ee
\end{Proposition}
\proof 	Let us first fix $\phi \in \Ac_+(y)$. Then,  $\vp:=\1_{Y^{y,\phi}_T\ge G}$ satisfies $\Pro{ Y^{y,\phi}_T\ge G}=\E[\vp]$ and $   G \vp\le Y^{y,\phi}_T$ so that $ \E^\Q[  \beta_TG \vp]\le  \E^\Q[ \beta_TY^{y,\phi}_T]$ $\le$ $y$, see Chapter \ref{CHAP: prix dual}. This shows that the left-hand side term in  \reff{eq: prob quantile bis} is smaller than the right-hand side term. Conversely, if $\vp \in L^0(\{0,1\})$ is such that $\E^\Q[\beta_TG \vp]\le y$, then it follows from Chapter \ref{CHAP: prix dual} that there exists $\phi \in \Ac_b(y)$ such that $Y^{y,\phi}_T\ge G \vp $. Since $G \vp\ge 0$, the super-martingale $Y^{y,\phi}$ remains non-negative so that $\phi \in \Ac_\p(y)$. Moreover, $Y^{y,\phi}_T\ge G$ on $\{\vp= 1\}$. Since    $\vp \in L^0(\{0,1\})$, this implies that $\Pro{Y^{y,\phi}_T\ge G}\ge   \Esp{\vp}$. 
\ep
\vs2

We next observe that the right-hand side problem in \reff{eq: prob quantile bis} can be interpreted as a statistical test problem:
 	\be\label{eq: prob quantile ter}	
	\sup\left\{ \Esp{\vp}\;,\;\vp \in L^0([0,1])\mbox{ s.t. } \E^{\Q_G}[\vp]\le y/p(G)\right\}\;,
	\ee
where $\Q_G$ is defined by 
	$$
	\frac{d\Q_G}{d\P}:=\frac{d\Q}{d\P} \frac{\beta_T G}{\E^\Q[\beta_T G]},
	$$
except that we look for a solution of the above test problem in  $L^0(\{0,1\})$. \vs2

The solution to this problem is given by  Neyman and Pearson's Lemma which we now state. 
\begin{Lemma}{\bf (Neyman and Pearson)}  Let $\P_0$ and $\P_1$  be two probability measures that are absolutely continuous with respect to $\P$. Given $\alpha\in ]0,1[$, the solution to the problem 
	$$
	\sup\left\{\E^{\P_1}\left[\xi\right]~:~\xi\in L^0([0,1]),\;\E^{\P_0}\left[\xi\right]\le \alpha\right\}\;,
	$$
is given by any random variable of the form
	$$
	\hat \xi:=\1_{\frac{d\P_1}{d\P}>\hat a \frac{d\P_0}{d\P}}+\hat \gamma \1_{\frac{d\P_1}{d\P}=\hat a \frac{d\P_0}{d\P}}
	$$
where
	$$
	\hat a:=\inf\left\{a>0~:~ \P_0\left[\frac{d\P_1}{d\P}>  a \frac{d\P_0}{d\P}\right]\le y \right\}
	$$
and $\hat \gamma \in [0,1]$ is such that  $\E^{\P_0}\left[\hat \xi\right]= \alpha$. 
\end{Lemma}

\begin{Remark}{\rm 
In the above Lemma,  $\xi$ has to be interpreted has a random test of Hyp$_0$ : $\P_0$ against Hyp$_1$ : $\P_1$. If the state of nature $\omega$ is such that $\xi(\omega)=p$, then one accepts  Hyp$_0$ with probability $1-p$. The quantity $\E^{\P_0}\left[\xi\right]$ corresponds  to the probability to reject  Hyp$_0$ while Hyp$_0$ is true (this is the risk of first kind), and  $\E^{\P_1}\left[\xi\right]$ corresponds  to the probability to reject Hyp$_0$ while   Hyp$_0$  is indeed false (this is called the power of the test).The test  $\hat \xi$ is called  UMP (uniformly most powerful) of size $\alpha$.}
\end{Remark}

Applying the above lemma to the problem \reff{eq: prob quantile ter} leads to an optimal solution $\hat \vp$ of the following  form:
\begin{Theorem} Assume that 
	$$
	\hat c:=\inf\left\{c>0~:~ \E^\Q\left[\beta_TG \1_{\frac{d\P}{d\Q}>  c \frac{d\Q_G}{d\Q}} \right]\le y\right\}
	$$
is such that 
$$
\E^\Q\left[\beta_TG \1_{\frac{d\P}{d\Q}>  \hat c \frac{d\Q_G}{d\Q}} \right]= y\;.
$$
 Then, the optimal solution to the problem \reff{eq: prob quantile} is given by the strategy $\hat \phi \in \Ac_\p(y)$ satisfying
$$
Y^{y,\hat \phi}_T=G \hat \vp
$$
where 
\b*
\hat\vp=\1_{\frac{d\P}{d\Q}>\hat c \frac{d\Q_G}{d\Q}}\;.
\e*
\end{Theorem} 
 
In most applications, $\hat c>0$ is such that    $\E^\Q\left[\beta_TG \1_{\frac{d\P}{d\Q}>  c \frac{d\Q_G}{d\Q}} \right]= y$, recall that $y<p(G)=\E^\Q[\beta_TG]$, so that the optimal strategy $\hat \phi$ satisfies 
$$
Y^{y,\hat \phi}_T=G\1_{\hat A}
$$
for $\hat A:=\{ d\P/d\Q>  \hat c\; d\Q_G/d\Q\}$. It means that the optimal solution consists in hedging a digital type option which pays $G$ on $\hat A$ and $0$ otherwise. 

Such a behavior is certainly not nice in practice since it may lead, as in general for discontinuous payoffs, to an explosion of the number of assets to have in the portfolio near to the maturity.

 Note that, right from the beginning, one could criticize the criteria which is only concerned with the probability of not missing the hedge but does not take into account of the sizes of the potential losses.

\vs2

In the case where $\hat c>0$ only satisfies    $\E^\Q\left[\beta_TG \1_{\frac{d\P}{d\Q}>  c \frac{d\Q_G}{d\Q}} \right]< y$, the above Theorem does not apply. However, the same reasoning can still be applied for the optimal {\sl success ratio} problem
	\be\label{eq: prob success ratio}
	\sup_{\phi \in \Ac_+(y)} \Esp{\frac{Y^{y,\phi}_T}{G}\wedge 1 } \mbox{ for some } 0<y<p(G),
	\ee
with the convention $z/0=\infty$ for $z\in \R$.  

\begin{Theorem}  The optimal solution to the problem \reff{eq: prob success ratio} is given by the strategy $\hat \phi \in \Ac_\p(y)$ satisfying
$$
Y^{y,\hat \phi}_T=G \hat \vp
$$
where 
\b*
\hat\vp=\1_{\frac{d\P}{d\Q}>\hat c \frac{d\Q_G}{d\Q}}+\hat \gamma \1_{\frac{d\P}{d\Q}=\hat c \frac{d\Q_G}{d\Q}}
\e*
with 
	$$
	\hat c:=\inf\left\{c>0~:~ \E^\Q\left[\beta_TG \1_{\frac{d\P}{d\Q}>  c \frac{d\Q_G}{d\Q}} \right]\le y\right\}
	$$
and where  $\hat \gamma \in [0,1]$ is such that  $\E^{\Q}\left[\beta_T G \hat \vp\right]= y$. 
\end{Theorem} 
 
Note that, when $\hat \gamma=0$, then it coincides with the solution of the quantile hedging problem.  

\subsection{Quantile hedging price}\label{subSection: quantile hedging price}

The quantile hedging price of the payoff $G$ is the minimal initial wealth that allows to hedge the option with a given probability of success, namely 
$$
p(G;\alpha):=\inf\left\{y\ge 0~:~\exists\;\phi \in \Ac_\p(y) \mbox{ s.t. } \P[Y^{y,\phi}_T\ge G]\ge \alpha\right\},\mbox{ for } \alpha \in [0,1]\;.
$$
Clearly, $p(G;1)=p(G)$ and $p(G;0)=0$. For $\alpha \in (0,1)$, it can be computed thanks to the results of the previous section. Indeed, given $0<y<p(G)$, one can find $\alpha(y) \in (0,1)$ such that 
$$
\alpha(y)=\sup_{\phi \in \Ac_+(y)} \Pro{Y^{y,\phi}_T\ge G}\;.
$$
Then, by definition, 
$$
p(G;\alpha)=\inf\left\{y\ge 0~:~\alpha(y)\ge \alpha\right\}\;.
$$ 

We shall see in Chapter \ref{CHAP: cible en proba} how the quantile hedging price can be directly related to a PDE, without having to invert the value function of an optimization problem, as suggested here. 

%%%%%%%%%%%%
\section{Hedging under expected loss constraints}\label{Section: expected shortfall dualite}

\subsection{Minimizing the expected shortfall}

In order to better take into account the amount of possible losses, we now consider   a risk control criteria of the form $\ell((G-V^{x,\phi}_T)^+)$, i.e.  we try to minimize
	\be\label{eq prob min loss}
	\inf_{\phi \in \Ac_+(y)} \Esp{\ell((G-Y^{y,\phi}_T)^+)}  \mbox{ for some } 0<y<p(G)\;,
	\ee
where, as above, $\Ac_+(y)$ is the restriction  of  $\Ac_b$ to strategies leading to non-negative wealth processes. 
Here, the {\sl loss function}   $\ell$ is $C^1$ strictly convex, increasing and defined on  $\R_+$, $\ell(0)=0$, and such that $\nabla \ell(+\infty)=\infty$, $\nabla \ell(0+)=0$.  We note $I$ $:=$ $(\nabla \ell)^{-1}$, the inverse of the derivative of $\ell$. As above, we assume that $G\in L^0(\R_\p)\setminus\{0\}$

\begin{Theorem}\label{thm: dualite expected shorfall} There exists a solution  $\hat \phi\in \Ac_\p(y)$ to the problem \reff{eq prob min loss}. It satisfies
	\b*
	Y^{y,\hat \phi}_T=\hat \vp(\hat c)G
	\e*
where, for $c>0$, 
	\b*
	\hat \vp(c):= \1_{G>0} \left(1- \frac{I(  c \beta_Td\Q/d\P)}{G}\wedge 1\right)\;,
	\e*
and $\hat c>0$ is the unique positive solution of 
	\b*
	\EspQ{\beta_T \hat \vp(c)G}=y\;. 
	\e*
\end{Theorem}

\proof 1. First of all, one can observe that  
	\b*
	\Esp{\ell((G-Y^{y,\phi}_T)^+)}=\Esp{\ell(G(1-\vp^\phi))}
	\e*
where $\vp^\phi:=[(Y^{y,\phi}_T/G)\wedge 1]\1_{G>0}$ satisfies $\EspQ{\beta_T \vp^\phi G}\le y$. Conversely, if $\vp\in L^0([0,1])$ satisfies the above constraint, then $\vp G$ can be reached by a financial portfolio starting from $y$ whose discounted value is a $\Q$-martingale and therefore remains non negative, since $G\ge 0$, see Section \ref{Section: chap prix dual complet} in Chapter \ref{CHAP: prix dual}.  The above problem is thus equivalent to 
	\be\label{eq probleme loss vp 1}
	\inf_{\vp\in L^0([0,1])} \Esp{\ell( (1-\vp)G )} \mbox{ under the constraint $\EspQ{\beta_T \vp G}\le y$.}
	\ee

\vs2

2.  We now check that existence holds in the above problem by using the following technical lemma which we state without proof. 
\begin{Lemma}{\bf (Komlos Lemma)} Let $(\zeta_n)_n$ be a sequence of random variables that are uniformly bounded in $L^1(\P)$. Then, there exists a sequence $(\bar \zeta_n)_n$ and a random variable  $\bar \zeta$ in $L^1(\P)$ such that $\bar \zeta_n\to \bar \zeta$ $\Pas$ and 
	$$
	\bar \zeta_n \in \mbox{conv}\left(\zeta_k,\;k\ge n\right) \;\;\Pas
	$$  
for all $n\ge 1$, where conv denotes the convex envelope.  
\end{Lemma}
Since $\vp\mapsto \Esp{\ell(G(1-\vp))}$ is convex, one deduces from the preceding Lemma that there exists a minimizing sequence $(\vp_n)_n$ which converges  $\Pas$ to some $\hat \vp$ in $L^0([0,1])$. One concludes by using Fatou's Lemma and the fact that $\ell\ge 0$. 

\vs2

3. We now check that $\hat \vp$ has the form given in the Theorem. Given   $\vp\in L^0([0,1])$ and $\eps\in [0,1]$, let us set
	\b*
	\vp_\eps:=\eps \vp + (1-\eps) \hat \vp
	\e*
and
	$$
	F_\vp(\eps):= \Esp{\ell( (1-\vp_\eps)G )}\;. 
	$$
Recall that $\ell$ is convex so that its derivative is non-decreasing. Using a monotone convergence argument, one then easily checks that the right-derivative $\nabla F_\vp(0+)$ of $F_\vp$ at $0$ exists and satisfies
	$$
	\nabla F_\vp(0+)= \Esp{\nabla\ell( (1-\hat \vp)G ) (\hat \vp-\vp) G}\;. 
	$$
Since $F_\vp$ is convex, because $\ell$ is convex, $\hat \vp$ should satisfy the first order optimality condition $\nabla F_\vp(0+)\ge 0$  for any $\vp \in L^0([0,1])$. 
This amounts to say that $\hat \vp$ satisfies 
	\be\label{eq prob test 1}
	 \E^{\Q_{\hat \vp}}\left[\hat \vp\right]\ge \E^{\Q_{\hat \vp}}\left[ \vp\right]
	\ee
for any $\vp \in L^0([0,1])$ such that, recall \reff{eq probleme loss vp 1}, 
	\be\label{eq prob test 2}
	\E^{\Q_G}\left[\vp\right]\le \frac{y}{p(G)}=:\alpha
	\ee
where $\Q_{\hat \vp}$ and  $\Q_{G}$ are the probability measures associated to the densities
	\b*
	\frac{d\Q_{\hat \vp}}{d\P}&=&\nabla\ell( (1-\hat \vp)G )    G / \Esp{\nabla\ell( (1-\hat \vp)G )   G}
	\\
	\frac{d\Q_{G}}{d\P} &=& \frac{d\Q}{d\P}  \beta_T G/\EspQ{\beta_T G}\;.
	\e*
As in the previous section, this can be interpreted as a random test: test the hypothesis $\Q_{\hat \vp}$ against  $\Q_{G}$ 
with a level $\alpha$.  It then follows from  Neyman and Pearson's Lemma, see above, that the optimal test $\hat \vp$ takes the value  $0$ if $d\Q_{\hat \vp}/d\P<c\;d\Q_G/d\P$ and the value   $1$ if $d\Q_{\hat \vp}/d\P>c\;d\Q_G/d\P$, for a given positive constant  $c$ which depends on the size of the test. First note that one should have $\hat \vp<1$ on $\{G>0\}$ since $\nabla \ell(0)=0$ and therefore  $d\Q_{\hat \vp}/d\P=0<d\Q_G/d\P$ when $\hat \vp=1$ and $G>0$. This implies that   $d\Q_{\hat \vp}/d\Q_G\le c$ on $\{G>0\}$. On $\{G=0\}$, one has $d\Q_{\hat \vp}/d\P=d\Q_G/d\P=0$,  and we set $\hat \vp=1$, see step 3 below. This leads to the definition of  $\hat \vp$ given in the Theorem. 

3. In order to justify that we can take  $\hat \vp=1$ on $\{G=0\}$, it suffices   to check that   $\hat c>0$  is such that $\E^{\Q_{G}}[\hat \vp(\hat c)]=y/p(G)$.  To see this recall that $\nabla \ell$ is increasing, continuous  and satisfies $\nabla \ell(+\infty)=\infty$ as well as $\nabla \ell(0+)=0$, by assumption. It follows that $I$ is increasing, continuous and satisfies $\nabla I(+\infty)=\infty$, $\nabla I(0+)=0$. This implies that $\hat \vp((0,\infty))=[0,\1_{G>0})$ $\Pas$ and that $c\in (0,\infty)\mapsto \hat \vp(c)$ is $\Pas$ continuous. Using the monotone convergence theorem, we then deduce that $c\in (0,\infty) \mapsto k(c):=\E^\Q[\beta_T G\hat \vp(c)]$ is continuous and satisfies $k((0,\infty))\supset(0,\E^\Q[\beta_T G\1_{G>0}])=(0,p(G))$. The uniqueness of $\hat c$ follows from the fact that $I$ is strictly increasing  and that $y/p(G)<1$ so that $\Pro{I(  \hat c \beta_Td\Q/d\P)<G}>0$. 
\ep

\subsection{Expected shortfall price}

As for  the quantile hedging approach, one can define an expected shortfall price:  
$$
\inf\left\{y\ge 0~:~\exists\;\phi \in \Ac_\p(y) \mbox{ s.t. } \E[\ell((G-Y^{y,\phi}_T)^\p)]\le l\right\},\mbox{ for } l\in \ell(\R_\p)\;.
$$
It can be deduced from the result of Theorem \ref{thm: dualite expected shorfall} by following the arguments of Section \ref{subSection: quantile hedging price} above. As for the quantile hedging price, we shall see in Chapter \ref{CHAP: cible en proba} how it can be directly related to a PDE. 
%%%%%%%%%%%%%%%%%%%%%%%%%%%%%%%%%%%%%%%%%%%%%%%%%%%%%%%%%%%%%%%%%%%%%%%%%%%%%%%%%%%%%%%%%%%%%%%%%%%%%%%%%%
%%%%%%%%%%%%%%%%%%%%%%%%%%%%%%%%%%%%%%%%%%%%%%%%%%%%%%%%%%%%%%%%%%%%%%%%%%%%%%%%%%%%%%%%%%%%%%%%%%%%%%%%%%
\part{The stochastic target approach}

\vs2

\chapter{Super-hedging problems}\label{CHAPTER: cible pas}

\section{Model and problem formulation}\label{Section: model formulation cible}

In this part, we consider a more general model in which the trading strategy $\phi$ may have an impact on the wealth process, namely the dynamics of the risky assets is given by
\be\label{eq: dyna X cible}
X^\phi_{t,x}(s)=x\p \int_t^s \mu(X^\phi_{t,x}(u),\phi_u) du \p \int_t^s \sigma(X^\phi_{t,x}(u),\phi_u) dW_u \;, 
\ee
where $W$ is the Brownian motion under the original probability measure $\P$. 

As in the previous chapter, the risk free interest rate is a function $\rho$ which depends only on $x$.

It follows that the wealth dynamics is given by 
\be\label{eq: dyna Y cible}
Y^\phi_{t,x,y}(s)=y\p \int_t^s \mu_Y(Z^\phi_{t,x,y}(u),\phi_u) du \p \int_t^s \sigma_Y(X^\phi_{t,x}(u),\phi_u) dW_u \;, 
\ee
where  
$$
\mu_Y(x,y,a):=a'\mu(x,a)\p (y-a'x)\rho(x) \mbox{ and } \sigma_Y(x,a):=a'\sigma(x,a)\;.
$$

The aim of this Chapter is to provide a PDE characterization of the  hedging price under constraint without appealing to the dual formulation of Chapter \ref{CHAP: prix dual}, which will anyway not be correct for the above model whenever $\mu$ and $\sigma$ depends in a non-trivial way of the strategy $\phi$. 

We recall that the associated value function is given by 
$$
v(t,x):=\inf\left\{y\in \R~:~\exists\; \phi \in \Ac_K\mbox{ s.t. } Y^\phi_{t,x,y}(T)\ge g(X_{t,x}^\phi(T))\right\}\,
$$
where the payoff function is assumed to be continuous, with linear growth and uniformly bounded from below. 
\vs2

We   assume in all this part  that $\mu$, $\sigma$ and $\rho$ are locally Lipschitz continuous, that \reff{eq: dyna X cible} admits a unique strong solution for any $\phi \in \Ac_K$, that there exists a unique solution $\psi(x,p)$ to the root problem
$$
\sigma_Y(x,a)=p'\sigma(x,a) \mbox{ for some } a\in \R^d,
$$
for any $(x,p) \in \R^d \x \R^d$, and that 
\be\label{hyp: psi loc lipschitz}
(x,p) \in \R^d\x \R^d\mapsto \psi(x,p) \mbox{ is locally Lipschitz. }
\ee

In order to prove the supersolution property stated below, we shall also assume 
\be\label{eq: explosion delta vol a l'infini}
\limsup_{|a|\to \infty} \inf_{(x,p)\in A} |\sigma_Y(x,a)-p'\sigma(x,a)|=\infty
\ee
for all compact set $A\subset \R^d\x \R^d$. 

\vs2

The results given below could be obtained in much more general situations, however this would require a substantially more technical analysis, see \cite{BET08}. 

\section{Geometric dynamic programming principle}

The main tool for providing a PDE characterization of $v$ is the {\sl geometric dynamic programming principle} of Soner and Touzi \cite{SoTo02b}, see also \cite{SoTo02} and \cite{BoNa08} for an extension to American type options. 

\begin{Theorem}\label{thm: geometric dpp} Fix $(t,x,y)\in [0,T]\x \R^{d\p 1}$.  Let $(\theta^\phi,\; \phi\in \Ac_K)$ denote a family of stopping times in $\Tc^t_{[t,T]}$. Then the following holds: 

{\rm \bf (DP1):} If $y>v(t,x)$, then there exists $\phi \in \Ac_K$ such that 
$$
Y_{t,x,y}^\phi(\theta^\phi) \ge v(\theta^\phi,X_{t,x}^\phi(\theta^\phi))\;.
$$

{\rm \bf (DP2):} If $y<v(t,x)$, then 
$$
\Pro{Y_{t,x,y}^\phi(\theta^\phi) > v(\theta^\phi,X_{t,x}^\phi(\theta^\phi))}<1\; \;\forall\;\phi \in \Ac_K\;.
$$
\end{Theorem}

We shall not provide a rigorous proof of this result and refer to \cite{SoTo02b} and the remarks in \cite{BoNa08}. We only explain the main argument. 

If $y>v(t,x)$, then, by definition of $v$, there exists $\phi \in \Ac_K$ such that $Y^\phi_{t,x,y}(T)\ge g(X_{t,x}^\phi(T))$. On the other hand, if $Y_{t,x,y}^\phi(\theta^\phi) < v(\theta^\phi,X_{t,x}^\phi(\theta^\phi))$ on a set of non zero measure, then starting from time $\theta^\phi$ is in not always possible to find a strategy $\tilde \phi$ such that $Y^{\tilde \phi}_{\theta^\phi,X_{t,x}^\phi(\theta^\phi),Y_{t,x,y}^\phi(\theta^\phi)}(T)\ge g(X_{\theta^\phi,X_{t,x}^\phi(\theta^\phi)}^{\tilde \phi}(T))$. This, combined with the flow property, contradicts the fact that  $Y^\phi_{t,x,y}(T)\ge g(X_{t,x}^\phi(T))$ $\Pas$ 

On the other hand, if  $Y_{t,x,y}^\phi(\theta^\phi) > v(\theta^\phi,X_{t,x}^\phi(\theta^\phi))$ $\Pas$, then starting from the time $\theta^\phi$, one can construct a strategy which allows to super-hedge the claim. This should imply that $y\ge v(t,x)$.

%%%%%%%%%%%%%%%%%%%%%%
\section{Derivation  of the pricing equation}

\subsection{PDE characterization}\label{Section: pde chara cible}

Before to provide the rigorous characterization of $v$, let us explain the main idea. Assume that  $v$ is smooth and that (DP1) above holds with $y=v(t,x)$, which would be the case if the infimum in the definition of $v$ was achieved. Then, one can find $\phi \in \Ac_K$ such that
$$ 
Y_{t,x,y}^\phi(\theta) \ge v(\theta,X_{t,x}^\phi(\theta))\;,
$$ 
for any stopping time $\theta \in \Tc^t_{[t,T]}$. Applying this formally for $\theta=t\p$, this implies that 
$$ 
Y_{t,x,y}^\phi(t\p) \ge v(t\p,X_{t,x}^\phi(t\p))\;,
$$ 
so that, by It\^o's Lemma, 
$$
\mu_Y(x,y,\phi_t) dt \p \sigma_Y(x,\phi_t)dW_t\ge \Lc^{\phi_t}v(t,x) dt \p Dv(t,x)'\sigma(x,\phi_t) dW_t
$$
where, for $a\in \R^d$ and a smooth function $\vp$, 
$$
\Lc^a\vp:=\partial_t \vp \p \mu(\cdot,a)'D\vp \p \frac12 {\rm Tr}[\sigma\sigma'(x,a) D^2\vp]\;. 
$$
Since the $dW$ term behaves  like $\sqrt{dt}\eps$ for a standard Gaussian random variable $\eps$, this necessarily implies that 
$$
\sigma_Y(x,\phi_t)=Dv(t,x)'\sigma(x,\phi_t) 
$$
so that
$$
\phi_t=\psi(x,Dv(t,x))
$$
by the definition of $\psi$ above.
Coming back to the previous inequality and recalling that $y=v(t,x)$ then leads to 
$$
\Gc v(t,x)\ge 0
$$
where, for a smooth function $\vp$,
$$
\Gc\vp:=\mu_Y(\cdot,\vp,\psi(\cdot, D\vp))- \Lc^{\psi(\cdot,D\vp)}\vp\;. 
$$
Moreover, $\phi_t$, and therefore $\psi(x,Dv(t,x))$, should take values in $K$. Recalling Proposition \ref{prop: caract K} in Section \ref{Section: chap dualite incomplete market} of Chapter \ref{CHAP: prix dual}, this implies that 
$$
\Hc v(t,x)\ge 0
$$
where
$$
\Hc\vp:=\inf_{|\zeta|=1} \left(\delta_K(\zeta)-\zeta'\psi(\cdot,Dv)\right)\ge 0
$$
for a smooth function $\vp$.
\vs2

The optimality included in the definition of $v$ defined as an infimum should actually show that one of the above inequalities is sharp, i.e. $v$ solves 
\be\label{eq: edp cible pas}
\min\left\{\Gc\vp\;,\;\Hc\vp\right\}=0\mbox{ on } [0,T)\x \R^d\;. 
\ee
This can be checked by using the second part (DP2) of the geometric dynamic programming principle. 

\begin{Theorem}\label{thm: charact visco cible pas} Assume that $v$ is locally bounded. Then, $v_*$ and $v^*$ are respectively viscosity super- and subsolutions of  \reff{eq: edp cible pas}. 
\end{Theorem}

The proof is divided in two parts. 

%%%%%%%%%%
\subsubsection{Viscosity supersolution property}

Before to prove the supersolution property of Theorem \ref{thm: charact visco cible pas}, we formulate the following remark to which we shall appeal in the proof. 

\begin{Remark}\label{rem: caract distance noyau par fenchel support K} {\rm Fix $(x,p)\in  A\subset \R^d\x \R^d$, with $A$ compact. Assume that there exists $\eps>0$ such that 
\be \label{eq: rem caract distance noyau par fenchel support K -1}
\inf_{|\zeta|=1} \left(\delta_K(\zeta)-\zeta'\psi\right)\le -\eps\;\mbox{ on } A\;.
\ee
Then, there exists $c_\eps>0$ such that 
\be \label{eq: rem caract distance noyau par fenchel support K}
\inf_{a\in K} |a-\psi| \ge c_\eps\;\mbox{ on } A\;.
\ee
This follows from the fact that $\inf_{a\in K} |a-\psi|=0$ implies $\psi\in K$, since $K$ is closed, which, together with Proposition \ref{prop: caract K}    in Section \ref{Section: chap dualite incomplete market} of Chapter \ref{CHAP: prix dual}, would imply that $\inf_{|\zeta|=1} \left(\delta_K(\zeta)-\zeta'\psi\right)\ge 0$.

Moreover, \reff{eq: explosion delta vol a l'infini} and \reff{eq: rem caract distance noyau par fenchel support K} implies that there exists $k_\eps>0$ such that 
\be \label{eq: rem caract distance noyau par fenchel support K \p 1}
\inf_{a\in K} |\sigma_Y(x,a)-p'\sigma(x,a)| \ge k_\eps\;\mbox{ for } (x,p)\in A\;.
\ee
Otherwise, we would find $(a,x,p)$ in a compact subset of $ K\x A$ such that $\sigma_Y(x,a)-p'\sigma(x,a)=0$, which would imply $a=\psi(x,p)$, a contradiction. }
\end{Remark}

We can now provide the proof of the supersolution property. 

\vs2

{\bf 1.} Fix   $(t_0,x_0) \in [0,T)\x \R^d$  and let $\vp$ be a smooth function such that 
	\be\label{eq:sursol min strict}
	{\rm (strict)}\min_{[0,T]\x \R^d}(v_*-\vp)= (v_*-\vp)(t_0,x_0)=0\;.
	\ee
Assume to the contrary that $\min\{\Gc\vp,\Hc\vp\}(t_0,x_0)< 0$. Then, by continuity of the operators, there exists $r,\eps>0$ such that $B_0:=B_r(t_0,x_0)\subset  [0,T)\x \R^d$ and 
$$
\min\{\Gc(\vp\p \zeta)\;,\;\Hc\vp\}\le -2\eps \mbox{ for } |\zeta|\le r \mbox{ on } B_r(t_0,x_0)   \;.
$$
Recalling Remark \ref{rem: caract distance noyau par fenchel support K} and the very definition of $\Hc$, this implies that 
\be\label{eq contra preuve supersolution 1}
&\mu_Y(x,y,a)-\Lc^a\vp(t,x) \le -\eps \;\forall\;(t,x,y,a) \in B_0\x\R\x K&\\
& \mbox{ s.t. }  |\sigma_Y(x,a)-Dv(t,x)'\sigma(x,a)| \le k_\eps \mbox{ and } |y-\vp(t,x)|\le r&\nonumber
\ee
for some $k_\eps>0$. 

For later use, observe  that, by \reff{eq:sursol min strict} and the definition of $ \vp$,  		
	\be
	\zeta := \min_{\partial_p B_\eps(t_0,x_0)} (v_*- \vp)>0\;,\label{eq sursol interieur def zeta} 
	\ee
where $\partial_p B_\eps(t_0,x_0)$ denotes the parabolic boundary of $B_\eps(t_0,x_0)$.

{\bf 2.} Let $(t_n,x_n)_{n\ge 1}$ be a sequence in $B_0$ which converges to $(t_0,x_0)$ and such that $v(t_n,x_n)\to v_*(t_0,x_0)$. Set $y_n=v(t_n,x_n)+n^{-1}$ and observe that 
	\be\label{eq iotan to 0}
	\gamma_n:=y_n-\vp(t_n,x_n)\to 0\;.
	\ee
For each $n\ge 1$, we have $y_n>v(t_n,x_n)$. It thus follows from (DP1) of Theorem  \ref{thm: geometric dpp}, that  there exists some $\phi^n\in \Ac_K$ such that
	\be\label{eq nun for superhedge 1}
	Y^n(t\wedge \theta^{n})\ge v(t\wedge \theta^n,X^n(t\wedge \theta^n)) \;\;\mbox{ for } t\ge t_n,
	\ee
where
	\b*
	Z^n=(X^n,Y^n):= \left(X^{\phi^n}_{t_n,x_n},Y^{\phi^n}_{t_n,x_n,y_n}\right)\;
	 \mbox{ and } \theta^n:=\theta^{o}_n\wedge \theta^1_n\;,
	\e*
with 
	\b*
	\theta^{o}_{n}&:=&\left\{s\ge t_n:(s,X^{\phi^n}_{t_n,x_n}(s))\notin B_0\right\}\\
	 \theta^1_{n}&:=&\left\{s\ge t_n:|Y_{t_n,x_n,y_n}^{\phi^n}(s)-\vp(s,X_{t_n,x_n}^{\phi^n}(s))|\ge r\right\}\;.
	\e*
Let us define
	\be
	A_n&:=&\left\{s\in [t_n,\theta_n]~:~ \mu_Y \left(Z^n(s),\phi^n_s\right) -\Lc^{\phi^n_s}\vp\left(s,X^n(s)\right)> -\eps \right\}\label{eq def Cn proof sur sol 0T V}\;,
	\ee 
and observe that \reff{eq contra preuve supersolution 1} implies that the process
$$
\delta^n_s:=\sigma_Y(X^n(s),\phi^n_s)-D\vp(s,X^n(s))'\sigma(X^n(s),\phi^n_s)
$$
satisfies
	\be\label{eq N>delta}
	\left|\delta^n_s\right|>k_\eps \mbox{ for } s\in A_n.
	\ee
{\bf 3.} Using \reff{eq nun for superhedge 1},   the definition of $\zeta$ in \reff{eq sursol interieur def zeta}  and the definition of $\theta_n$, we then obtain
	\b*
	Y^n(t\wedge \theta_n)
	&\ge&  
	   \vp\left(t\wedge \theta^n,X^n(t\wedge \theta^n)\right) +\left(\zeta\1_{\{\theta^o_n=\theta^n \}}+r\1_{\{\theta^o_n>\theta^n \}}\right)\1_{\{t\ge \theta^n\}}
	\\
	&\ge& 
	   \vp\left(t\wedge \theta^n,X^n(t\wedge \theta^n)\right)
           +\left(\zeta\wedge r\right)\1_{\{t\ge \theta^n\}}\;,\;\;t\ge
           t_n\;.
	\e*
Since $\varphi$ is smooth, it follows from It\^{o}'s Lemma, \reff{eq contra preuve supersolution 1}, \reff{eq iotan to 0} and  the definition of $\delta^n$  that
	\be
	-\left(\zeta\wedge r\right)\1_{\{t<\theta^n\}}
	&\le& K^n_t
	  \;, \label{supersolgirsanov}
	\ee
where   	
	\b*
	K^n_t&:=&\gamma_n-(\zeta\wedge r)
	+ \int_{t_n}^{t\wedge \theta^n} b^n_s ds\p \int_{t_n}^{t\wedge \theta^n} \delta^n_s dW_s\;, \nonumber
	\e*
with
	\b*
	b^n_s&:=& \left[\mu_Y\left(Z^n(s),\phi^n_s\right)-\Lc^{\phi^n_s}\vp\left(s,X^n(s)\right)\right]\1_{A_n}(s)    \;.
	\e*
Let $M^n$ be the exponential local martingale defined by $M^n_{t_n}=1$ and, for $s\ge t_n$,
	\b*
	dM^n_s &=& -M^n_s   b^n_s |\delta^n_s|^{-2} \delta^n_s  dW_s
	\,,
	\e*
which is well defined by \reff{eq N>delta} and the Lipschitz continuity of the coefficients. 
By It\^o's formula  and \reff{supersolgirsanov}, we see that $M^nK^n$ is a  local  martingale which is bounded from below by the submartingale  $-\left(  \zeta\wedge r\right)  M^n$. Then, $M^n K^n$  is a supermartingale, and    it follows from \reff{supersolgirsanov} that 
	$$
	0
	\le 
	\Esp{ M^n_{\theta^n} K^n_{\theta^n}}
	\le  \gamma_n-(\zeta\wedge r) <0
\;,
	$$
for $n$ large enough,
recall \reff{eq iotan to 0}, which leads to a contradiction. \ep

%%%%%%
\subsubsection{Viscosity subsolution property}

{\bf 1.}  Fix   $(t_0,x_0) \in [0,T)\x \R^d$  and let $\vp$ be a smooth function such that 
	\be\label{eq:soussol max strict}
	{\rm (strict)}\max_{[0,T]\x \R^d}(v^*-\vp)= (v^*-\vp)(t_0,x_0)=0\;.
	\ee
    We assume to the contrary that 
\be\label{eq: soussol point int SG}
\min\{\Gc\vp\;,\; \Hc\vp\}(t_0,x_0)\ge  2\eta
\ee
for some $\eta>0$, and  work towards a contradiction. 

 Under the above assumption, we may find $r>0$  such that 
 	\be
	 & \mu_Y (\cdot,\vp\p \zeta ,\psi(\cdot,D\vp)) -\Lc^{\psi(\cdot,D\vp)} \vp >  \eta \mbox{ for } |\zeta|\le r\mbox{ on } B_0:=B_r(t_0,x_0)\;. &\label{eq contra preuve soussolution 1}
	 \ee
For later use note  that, by \reff{eq:soussol max strict} and the definition of $ \vp$,  		
	\be
	-\zeta := \max_{\partial_p B_r(t_0,x_0)} (v^*- \vp)<0\;.\label{eq soussol interieur def zeta} 
	\ee
Moreover, we can find a sequence $(t_n,x_n)_{n\ge 1}$ in $B_0$ which converges to $(t_0,x_0)$ and such that $v(t_n,x_n)\to v^*(t_0,x_0)$. Set $y_n=v(t_n,x_n)-n^{-1}$ and observe that 
	\be\label{eq iotan to 0 soussol}
	\gamma_n:=y_n-\vp(t_n,x_n)\to 0\;.
	\ee
{\bf 2.} We now let $Z^n:=(X^n,Y^n)$ denote the solution of \reff{eq: dyna X cible}-\reff{eq: dyna Y cible} associated to the Markovian control $\hat \phi^n:=\psi(\cdot,D\vp(\cdot,X^n))$ and  the initial condition $Z^n(t_n)=(x_n,y_n)$, recall that $\psi$ is assumed to be locally Lipschitz. We next define the stopping times 
	\b*
	\theta^o_n
	&:=&
	\inf\left\{s\ge t_n~:~(s,X^n(s))\notin B_0\right\},
	\\
	\theta_n
	&:=&
	\inf\left\{s\ge t_n~:~|Y^n(s)-\vp(s,X^n(s))|\ge r\right\}\wedge \theta^o_n\;.
	\e*
Note that, by definition of $\hat \phi^n$ and \reff{eq contra preuve soussolution 1}, $Y^n-\vp(\cdot,X^n)$ is non-decreasing on $[t_n,\theta_n]$, so that 
 \be\label{eq: preuve soussol int y ge vp}
 Y^n(\theta_n)-\vp(\theta_n,X^n(\theta_n))\ge y_n-\vp(t_n,x_n)=\gamma_n>-r
 \ee
 for $n$ large enough, recall \reff{eq iotan to 0 soussol}.
 Since $\varphi\ge v^*\ge v$, it follows that 
 \b*
 Y^n(\theta_n)-v\left(\theta_n,X^n(\theta_n)\right)
 &\ge&
 \1_{\{\theta_n<\theta^o_n\}}\left\{Y^n(\theta_n)-\varphi\left(\theta_n,X^n(\theta_n)\right)\right\}
 \\
 & &\hspace{10mm}
 +\1_{\{\theta_n=\theta^o_n\}}\left\{Y^n(\theta^o_n)-v^*\left(\theta^o_n,X^n(\theta^o_n)\right)\right\}
 \\
 &=&
 r\1_{\{\theta_n<\theta^o_n\}}+\1_{\{\theta_n=\theta^o_n\}}\left\{Y^n(\theta^o_n)-v^*\left(\theta^o_n,X^n(\theta^o_n)\right)\right\}
 \\
 &\ge&
 r\1_{\{\theta_n<\theta^o_n\}}+\1_{\{\theta_n=\theta^o_n\}}\left\{Y^n(\theta^o_n)+\zeta-\varphi\left(\theta^o_n,X^n(\theta^o_n)\right)\right\}
 \\
 &\ge&
 r\wedge\zeta
 +\1_{\{\theta_n=\theta^o_n\}}\left\{Y^n(\theta^o_n)-\varphi\left(\theta^o,X^n(\theta^o_n)\right)\right\}\;.
 \e* 
In view of \reff{eq: preuve soussol int y ge vp}, this leads to
 \b*
 Y^n(\theta_n)-v\left(\theta_n,X^n(\theta^n)\right)
 &\ge&
 (r\wedge\zeta)/2
 \e*
for $n$ large enough, since   $\gamma_n\to 0$.
Recalling   that   $y_n=v(t_n,x_n)-n^{-1}<v(t_n,x_n)$, this is clearly in contradiction with (DP2) of Theorem \ref{thm: geometric dpp}. 
\ep 

%%%%%%%%%%%%%%%%%%%%%%%%%%%%%%%
\subsection{Boundary condition as $t=T$}

Note that by construction $v(T,\cdot)=g$. However, it follows from the previous sections that $v$ satisfies $\Hc v\ge 0$, in the viscosity sense, which implies that $Dv$ is constrained on $[0,T)$. This constraint should propagate up to $T$. Hence, $v$ should solve
\be\label{eq: cible cond bord T}
\min\{\vp-g\;,\;\Hc \vp\}=0\mbox{ on } \{T\}\x \R^d\;.
\ee
We shall see below that this boundary condition is naturally related to the face-lifting phenomenon observed in    Chapters \ref{CHAP: prix dual} and  \ref{CHAP: pricing equation II}.

\begin{Theorem}\label{thm: condition en T cible pas} Assume that $v$ is locally bounded. Then, $v_*$ and $v^*$ are respectively super- and subsolution of \reff{eq: cible cond bord T}.
\end{Theorem}

\proof The proofs follow from similar arguments as in the previous section, up to the standard trick which consists in adding a term of the form $\pm \sqrt{T-t\p \alpha}$ to the test function $\vp$, so that, for $t$ close to $T$ and $\alpha>0$ small enough, it satisfies
$$
\pm \Gc(\vp\pm \sqrt{T-\cdot\p \alpha})\ge 0\;. 
$$
We only explain the argument for the subsolution property. The supersolution property is proved by using the same trick combined with the arguments used to prove the supersolution property in  Section \ref{Section: pde chara cible}. 
  
 Let $x_0 \in \R^d$ and $\vp$ be a smooth function such that 
	$$
	{\rm (strict)}\max_{  [0,T]\x \R^d}(v^*-\vp)= (v^*-\vp)(T,x_0)=0\,.
	$$
Assume that  
 \be\label{eq delta vp g V}
\min\{v^*-g\;,\;\Hc\vp\}(T,x_0)\ge 4\eta\;.
 \ee
Set $\tilde \vp(t,x):=\vp(t,x)+\sqrt{T-t\p \alpha}-\sqrt{\alpha}$. Since $\partial_t \tilde \vp(t,x) \to -\infty$ as $t\to T$ and $\alpha\to 0$, we deduce that,  for $r,\alpha>0$ small enough,  
 	\be
	 & \min\left\{\tilde \vp-g\;,\;\mu_Y\left(\cdot,\tilde \vp \p \zeta ,\psi(\cdot,D\tilde \vp)\right) -  \Lc^{\psi(\cdot,D\tilde \vp)}\tilde \vp\;,\;\Hc\tilde \vp\right\}\ge  \eta &
	 \nonumber \\
	 & \mbox{ for } |\zeta|\le r\mbox{ on } B_0:=[T-r,T]\x B_r(x_0)\;.&\label{eq contra preuve subsolution bord b}
	\ee   
Also observe that, since   $(v^*-\tilde \vp)(T,x_0)=0$ and $(T,x_0)$ achieves a strict maximum, we can choose $r>0$ so that
	\be\label{eq V-vp le rho sur 2}
	v^*(t,x)\le \tilde \vp(t,x)- \eps/2\;\;\mbox{ for all }\; (t,x) \in [T-r,T]\x \partial B_r(x_0)\;,
	\ee	
which, together with   $v(T,\cdot)=g$  and \reff{eq contra preuve subsolution bord b},  leads to
\be\label{eq contra preuve subsolution bord zeta}
	 v(t,x)-\tilde \vp(t,x)   \le -\zeta   \;\mbox{ for all }\; (t,x) \in \partial_p B_0\;
	\ee 
for some $r,\eps,\zeta>0$ small enough but so that the above inequalities still hold. 

By following the arguments of the previous  section,  we deduce that \reff{eq contra preuve subsolution bord b} and    \reff{eq contra preuve subsolution bord zeta} lead to a contradiction of (GDP2). 

\ep

%%%%%%%%%
\section{Extension to more general dynamics}

It should be noted that the above proofs and results do not depend on the specific form of $\mu_Y$ and $\sigma_Y$ defined in Section \ref{Section: model formulation cible}, but only on the general assumptions we made. 
\vs2

This implies that much more general dynamics could be considered. In particular, we could set 
$$
\mu_Y(x,y,a):=y\left(a'[x]^{-1}\mu(x,a)\p (1-a'{\bf 1})\rho(x)\right) \mbox{ and } \sigma_Y(x,a):=ya'[x]^{-1}\sigma(x,a)\;,
$$
with ${\bf 1}=(1,\ldots,1)$ and $[x]$ denoting the diagonal matrix with $i$-th diagonal component given by $x^i$. In this case,  the dynamics of $Y$ is given by 
\be\label{eq: def Y en proportion richesse cible}
dY^\phi_{t,x,y}(s)&=&  Y^\phi_{t,x,y}(s)\phi_s'[X^\phi_{t,x}(s)]^{-1}dX^\phi_{t,x}(s)\\
 \nonumber &&\p  \left(Y^\phi_{t,x,y}(s)-Y^\phi_{t,x,y}(s)\phi_s'{\bf 1}\right)\rho(X_{t,x}(s)) ds   \;.
\ee
This corresponds to a model where $\phi^i$ denotes the proportion of the wealth invested in the $i$-th risky asset. In this case, we have to put restrictions on the coefficient $\mu$ and $\sigma$ in order to ensure that $X$ has positive components whenever the initial condition belongs to $(0,\infty)^d$, and the viscosity solution properties have to be stated on $(0,\infty)^d$ instead of $\R^d$. 

%%%%
\section{Examples in the Black and Scholes model}

{\bf a.} If $\sigma(x)=x\sigma$ and $\mu(x)=x\mu$ where $\sigma>0$ and $\mu$ is a real constant, then $\psi(x,p)=p$. Moreover, if $K=\R$, then $\delta_K(\zeta)=\infty$ for $\zeta\ne 0$ and therefore if $|\zeta|=1$. It follows that $v$ is a discontinuous viscosity solution of 
\b*
0&=&D\vp x\mu\p (\vp-D\vp x) \rho-\partial_t \vp - x\mu D\vp -\frac12 x^2\sigma^2 D^2\vp
\\
&=&\rho \vp -\partial_t \vp - \rho x D\vp  -\frac12 x^2\sigma^2 D^2\vp
\e*
which is \reff{eq pde FK} of Chapter \ref{CHAP: pricing equation complet}. Moreover, the boundary condition at $t=T$ is simply given by $g$. 
\vs2

{\bf b.} If $K\ne \R$, then the same computations lead to \reff{eq: HJB marche incomplet} of Chapter \ref{CHAP: pricing equation II}. As for the boundary condition at $t=T$, we obtain 
$$
\min\left\{\vp(T,\cdot)-g\;,\;\inf_{|\zeta|=1} \delta_K(\zeta)- \zeta D \vp(T,\cdot)\right\}=0\;.
$$
in the discontinuous viscosity sense. One can show that $\hat g$ defined in  Chapter \ref{CHAP: prix dual} is the minimal supersolution of this equation, and that the above characterization of $v$ actually implies that  $v^*(T,\cdot)=v_*(T,\cdot)=\hat g$ as demonstrated in Chapter \ref{CHAP: pricing equation II}.
\vs2

{\bf c.} Let us now consider the case of the Black and Scholes model with $Y$ defined as in \reff{eq: def Y en proportion richesse cible}, i.e. where $\phi$ represents the proportion of the wealth invested in each asset. Then, for $g\ge 0$ so that $v>0$,  the PDE \reff{eq: edp cible pas} reads: 
\b*
\min\left\{\rho \vp -\partial_t \vp  - \rho x D\vp -\frac12 x^2\sigma^2 D^2\vp\;,\; \inf_{|\zeta|=1} \delta_K(\zeta)\vp- \zeta xD\vp  \right\}=0 
\e*
$\mbox{ on } [0,T)\x (0,\infty)$,
and the boundary condition is given by 
\b*
\min\left\{\vp(T,\cdot) -g\;,\; \inf_{|\zeta|=1} \delta_K(\zeta)\vp(T,\cdot)- \zeta xD\vp(T,\cdot)  \right\}=0 
\e*
$\mbox{ on }   (0,\infty)$. Under suitable assumptions, one can show that the smaller supersolution of the above equation is given by 
$$
\check g(x):=\sup_{\zeta \in K}e^{-\delta_K(\zeta)}g(xe^\zeta)
$$
and that $v(T-,\cdot)$ actually coincides with $\check g$, see e.g. \cite{SoTo02c}.

\chapter{Approximate hedging with  controlled risk}\label{CHAP: cible en proba}

We now turn to quantile and shortfall based pricing problems.  More precisely, we let $\Psi$ be a given real valued measurable function on  $\R^d\x \R_\p$, satisfying
$$
y\in \R_\p \mapsto \Psi(x,y) \mbox{ is non-decreasing for all } x \in \R^d\;,
$$
 and define
$$
v(t,x,p):=\min\left\{y\ge 0~:~\exists\; \phi \in \Ac_K\;\mbox{ s.t. } \Esp{\Psi(X^\phi_{t,x,y}(T),Y^\phi_{t,x,y}(T))}\ge p\right\}\;.
$$

For $\Psi(x,y)=\1_{y\ge g(x)}$, this corresponds to the quantile hedging problem discussed in Section \ref{Section: quantile hedging dualite} of Chapter \ref{CHAP: quantile et loss par dualite}.  For $\Psi(x,y)=-\ell((g(x)-y)^\p)$, this corresponds to the expected loss pricing rule of Section \ref{Section: expected shortfall dualite} of Chapter \ref{CHAP: quantile et loss par dualite}.

\vs2

The aim of this chapter is to show how such problems can be embedded into the class of general stochastic target problems as  discussed in Chapter \ref{CHAPTER: cible pas}.

\vs2

In the rest of this chapter, we shall often write $Z^\phi_{t,x,y}$ for $(X^\phi_{t,x},Y^\phi_{t,x,y})$. 

\section{Problem reduction}

The key point is the following observation made in \cite{BET08}. 
\vs2

\begin{Proposition}\label{prop: reformulation cible en moment}Fix $(t,x)\in [0,T]\x \R^d$ and assume that $\Psi(Z^\phi_{t,x,y}(T))\in L^2$ for any $\phi\in \Ac_K$ and $y\ge 0$. Then, 
\begin{equation}\label{eq: reformulation cible en moment}
v(t,x,p)=\min\left\{y\ge 0~:~\exists\; (\phi,\alpha) \in \Ac_K\x L^2_{\Pc}\;\mbox{ s.t. } \Psi(Z^\phi_{t,x,y}(T))\ge P^\alpha_{t,p}(T)\right\}
\end{equation}
where 
$$
P_{t,p}^\alpha:=p\p \int_0^\cdot \alpha_s' dW_s\;.
$$
Moreover, the above terms are also equal to 
$$
\min\left\{y\ge 0~:~\exists\; (\phi,\alpha) \in \Ac_K\x L^2_{\Pc}\;\mbox{ s.t. } Y^\phi_{t,x,y}(T)\ge \Psi^{-1}(X^\phi_{t,x}(T), P^\alpha_{t,p}(T))\right\}
$$
where $\Psi^{-1}$ denotes the right inverse of $\Psi$ in the $y$-variable.
\end{Proposition}

\proof Let $\bar v(t,x,p)$ denote the right-hand side in \reff{eq: reformulation cible en moment}. Then, for $y>\bar v(t,x,p)$, there exists $(\phi,\alpha) \in \Ac_K\x L^2_{\Pc}$ such that 
$\Psi(Z^\phi_{t,x,y}(T))\ge P^\alpha_{t,p}(T)$. Since $P^\alpha_{t,p}$ is a martingale, taking expectation leads to $\Esp{\Psi(Z^\phi_{t,x,y}(T))}\ge p$. This implies that $\bar v(t,x,p)\ge v(t,x,p)$. Conversely, if $y>v(t,x,p)$ then 
there exists $\phi  \in \Ac_K$ such that $M_0:=\Esp{\Psi(Z^\phi_{t,x,y}(T))}\ge p$. Let us define the martingale $M:=\Esp{\Psi(Z^\phi_{t,x,y}(T))~|~\Fc_\cdot}$. It follows from the martingale representation theorem, see Theorem \ref{thm: representation martingale} of Chapter \ref{CHAP: prix dual} or \cite{KS91},  that there exists $\alpha \in L^2_{\Pc}$ such that $M=P_{t,M_0}^\alpha$, with $P$ defined as in the proposition. In particular, 
$\Psi(Z^\phi_{t,x,y}(T))=P^\alpha_{t,M_0}(T)\ge  P^\alpha_{t,p}(T)$, since $M_0\ge p$. This proves that  $\bar v(t,x,p)\le v(t,x,p)$ and concludes the proof. 
\ep
\vs2

Otherwise stated, it suffices to consider an augmented system $(X,Y,P)$ with an augmented control $(\phi,\alpha)$ and apply the technics introduced in Chapter \ref{CHAPTER: cible pas} above. In particular, the geometric dynamic programming of Chapter \ref{CHAPTER: cible pas}  applies here.

\begin{Theorem}\label{thm: geometric dpp en moment}
Fix $(t,x,y,p)\in [0,T]\x \R^{d\p 1}\x \R$.  Let $(\theta^{\phi,\alpha},\; (\phi,\alpha) \in \Ac_K\x L^2_{\Pc})$ denote a family of stopping times in $\Tc^t_{[t,T]}$. Then the following holds: 

{\rm \bf (DP1):} If $y>v(t,x,p)$, then there exists $(\phi,\alpha) \in \Ac_K\x L^2_{\Pc}$ such that 
$$
Y_{t,x,y}^\phi(\theta^{\phi,\alpha}) \ge v(\theta^{\phi,\alpha},X_{t,x}^\phi(\theta^{\phi,\alpha}),P_{t,p}^\alpha(\theta^{\phi,\alpha}))\;.
$$

{\rm \bf (DP2):} If $y<v(t,x)$, then 
$$
\Pro{Y_{t,x,y}^\phi(\theta^{\phi,\alpha}) > v(\theta^{\phi,\alpha},X_{t,x}^\phi(\theta^{\phi,\alpha}),P_{t,p}^\alpha(\theta^{\phi,\alpha}))}<1\; \;\forall\;  (\phi,\alpha) \in \Ac_K\x L^2_{\Pc}\;.
$$
\end{Theorem}

%%%%%%%%%%%%%%%%%
\section{Pricing equation}

In view of Theorem \ref{thm: geometric dpp en moment}, one can now apply the same arguments as in Section \ref{Section: pde chara cible} of Chapter  \ref{CHAPTER: cible pas}. The only difference is that we now have to take into account a new control $\alpha$ and a new state process $P$. 

\subsection{In the domain}

Before to state the PDE characterization in the domain, let us first  introduce the notations  corresponding to our stochastic target problem.  
\vs2

First, we assume that the equation 
$$
\sigma_Y(x,a)=p'\sigma(x,a) \p q b'  \mbox{ for some } a \in \R^d
$$
admits a unique solution $\psi(x,p,q,b)$ which is locally Lipschitz continuous. The Dynkin operator associated to $(X,P)$ for the value of the control $(a,b)$ is denoted by 
$$
\Lc^{a,b}\vp:=\Lc^a \vp\p  D^2_{xp}\vp' \sigma(\cdot,a) b  \p  \frac12 b^2 D^2_p\vp
$$
where $\Lc^a$ is defined as in the previous chapter, $D^2_{xp}\vp$ stands for the second order cross derivatives $(\partial^2\vp/\partial_{x^ip})_{i\le d}$ and $  D^2_p\vp$ is the second derivative with respect to $p$. 

We then consider the counterparts of  operators $\Gc$ and $\Hc$ associated to $a=\psi(\cdot,b)$, for $b$ given: 
\b*
&\Gc^b\vp:=\mu_Y(\cdot,\psi_\vp^b)-\Lc^{\psi_\vp^b,b}\vp\;\mbox{ and }\;\Hc^b\vp:=\inf_{|\zeta|=1}(\delta_K(\zeta)-\zeta'\psi_\vp^b )&\\
&\mbox{ with } \psi_\vp^b:=\psi(\cdot,D_x\vp,D_p\vp,b)\;.&
\e*

In the following, we set 
$$
\Oc:=\{p\in \R~:~0<v(t,x,p)<\infty \mbox{ for all } (t,x)\in [0,T]\x \R^d\}
$$
and we assume that 
$$
\Oc \mbox{ is non-empty, convex and closed. }
$$
Note that the convexity is obvious, and is indeed not an assumption, whenever $\Oc$ is non-empty. 
\vs2

In what follows $v_*$ and $v^*$ are defined as the semicontinuous envelopes of $v$ in the three variables $(t,x,p)$ when approximated by a sequence $(t_n,x_n,p_n)\in [0,T)\x \R^d\x {\rm int}(\Oc)$

\begin{Theorem}\label{thm: cara visco cible en moment}  The following holds: 

{\rm (i)} If $K$ is compact, then 
$v_*$ is a viscosity supersolution of
 \be\label{eq: pricing equation cible en moment supersol}
 \max\left\{\sup_{b\in \R^d}\min\left\{\Gc^b\vp\;,\;\Hc^b \vp\right\}\;,\;-|D_{p}\vp|\right\}=0\mbox{ on } [0,T)\x \R^d\x \Oc\;.
\ee
{\rm (ii)} $v^*$  is a viscosity subsolution of
 \be\label{eq: pricing equation cible en moment subsol}
 \sup_{b\in \R^d}\min\left\{\vp\;,\;\Gc^b\vp\;,\;\Hc^b \vp\right\}=0\mbox{ on } [0,T)\x \R^d\x \Oc\;.
\ee
\end{Theorem}

\proof  We do not provide the entire proof because, thanks to Theorem \ref{thm: geometric dpp en moment},  it follows exactly the line of arguments of Section \ref{Section: pde chara cible} of Chapter \ref{CHAPTER: cible pas}. We only explain an additional technical point which should be taken into account in order to derive the supersolution property. Namely, in Section \ref{Section: pde chara cible} of Chapter \ref{CHAPTER: cible pas} we used the assumption \reff{eq: explosion delta vol a l'infini} in order to deduce \reff{eq: rem caract distance noyau par fenchel support K \p 1} from \reff{eq: rem caract distance noyau par fenchel support K -1}. Here the problem comes from the new control $b$ that is a-priori not bounded. However, if $D_p \vp\ne 0$ and $a\in K$ solves $\sigma_Y(x,a)=D_x\vp'\sigma(x,a) \p D_p\vp b'$, then the fact that $K$ is compact along with the regularity assumptions on $\sigma_Y$ and $\sigma$ imply that $b$ has to belong to a compact set. As a conclusion, the proof of the supersolution property can be reproduced without difficulty when $D_p \vp\ne 0$. When $D_p \vp= 0$, the viscosity supersolution property is satisfied by construction.  As for the subsolution property, nothing changes except that we need to have $y_n=v(t_n,x_n)-n^{-1}\ge 0$ in the proof of Section \ref{Section: pde chara cible} of Chapter \ref{CHAPTER: cible pas}, see just before \reff{eq iotan to 0 soussol}, since the initial wealth should be non-negative in the definition of our criteria. In order to ensure this, we need to have $v^*>0$ at the point where the maximum of the difference with the test function is achieved. 
\ep

%%%%
\subsection{Boundary condition at $t=T$}

By similar arguments, the boundary condition  of Theorem \ref{thm: condition en T cible pas} of Chapter \ref{CHAPTER: cible pas} extends to this context. 

\begin{Theorem}  Assume   that  $\Psi^{-1}$ is continuous on $\R^d\x  \Oc$. Then the following holds: 

{\rm (i)} If $K$ is compact, then 
$v_*$ is a viscosity supersolution of
 \be\label{eq: pricing equation cible en moment supersol en T}
 \max\left\{\min\left\{\vp-\Psi^{-1}\;,\;\sup_{b\in \R^d} \Hc^b \vp\right\}\;,\;-|D_{p}\vp|\right\}=0\mbox{ on } \{T\}\x \R^d\x\Oc\;.
\ee
{\rm (ii)}  $v^*$  is a viscosity supersolution of
 \be\label{eq: pricing equation cible en moment subsol en T}
\min\left\{\vp\;,\;\vp-\Psi^{-1}\;,\;\sup_{b\in \R^d} \Hc^b \vp\right\}=0\mbox{ on } \{T\}\x\R^d\x\Oc\;.
\ee
\end{Theorem}

In our context of financial mathematics, one can usually say a little bit more  on the boundary condition whenever their exists a well-behaved martingale measure. To see  this, 
 let  $H_{t,x}$ be defined by 
 \b*
 H_{t,x}(s)=1-\int_t^s H_{t,x}(u) \lambda(X_{t,x}(u)) dW_u \;\mbox{ with } \lambda(x):=\sigma^{-1}(\mu(x)-\rho(x) x)\;,
 \e*
where we implicitly assume that $\sigma$ is invertible and that $H$ is well-defined as a martingale for any initial conditions $(t,x)$. 

\begin{Proposition}\label{prop: cond en T par env convexe de G-1} Assume that $\Oc$ is compact.  Fix $(x,p) \in \R^d\x\Oc$   and assume that for all sequence $(t_n,x_n,p_n)_n\subset   [0,T)\x \R^d\x {\rm int}(\Oc)$ that converges to $(T,x,p)$, and for all sequence $({\phi_n})_n \subset \Ac_K$, we have 
\b*
& \Esp{|H_{t_n,x_n}(T) \beta_{t_n,x_n}(T)\widehat \Psi^{-1}(X^{\phi_n}_{t_n,x_n}(T),p)-\widehat \Psi^{-1}(x,p)|}\to 0\;&\\
& \mbox{ and }\;\;\;
\Esp{ |H_{t_n,x_n}(T)  \beta_{t_n,x_n}(T) \nabla^\p    \widehat \Psi^{-1}(X^{\phi_n}_{t_n,x_n}(T),p)-\nabla^\p    \widehat \Psi^{-1}(x,p)|}\to 0\;,&
\e*
where $\widehat \Psi^{-1}$ denotes the convex envelope of $\Psi^{-1}$ with respect to $p$, and $\nabla^\p \widehat \Psi^{-1}$ its right-derivative with respect to $p$. Then, 
$$
v_*(T,x,p)\ge \widehat \Psi^{-1}(x,p)\;.
$$
\end{Proposition}

When $\widehat \Psi^{-1}$ and $ \nabla^\p    \widehat \Psi^{-1}$ are continuous with   polynomial growth in $x$,   $\mu$ and $\sigma$ are uniformly Lipschitz in $x$, and $\lambda$ is bounded, then the above assumptions are trivially satisfied.
\vs2

Note that, if   \reff{eq: pricing equation cible en moment subsol en T} holds on $ \{T\}\x\R^d\x \Oc$,     $v^*$ can be shown to be convex and if $v_*$ is strictly increasing in $p$, then this implies that $v_*$ and $v^*$ are super- and subsolutions of 
$$
\min\left\{\vp\;,\;\vp-\widehat \Psi^{-1}\;,\;\sup_{b\in \R^d} \Hc^b \vp\right\}=0\mbox{ on } \{T\}\x\R^d\x\Oc\;,
$$
using the fact that $v_*\ge 0$ by construction.
In the limiting case where $K=\R^d$, and therefore $\Hc^b \vp\equiv \infty$, then the boundary condition simply reads 
$$
\vp(T,\cdot)=\widehat \Psi^{-1}\vee 0\;.
$$
Otherwise stated, a first face-lift of the natural terminal condition $\Psi^{-1}$ is due to the additional state process $P$. When $K\ne \R^d$, then an additional face-lift is required  as explained in Chapter \ref{CHAPTER: cible pas}. We shall provide two examples in Sections \ref{sec: example cible en moment 1} and \ref{sec: example cible en moment} below. 
\\

We conclude this section with the proof of the above proposition\footnote{There is a slight error in the proof of the corresponding result in \cite{BET08}, see their Proposition 3.2 in which $P$ is $\P$-martingale and not a $\Q$-martingale. We take this opportunity to correct it and we thank Nizar Touzi for the discussions we had on this point.  A rigorous version is given in Moreau \cite{moreau}}.\\
{\bf Proof of Proposition  \ref{prop: cond en T par env convexe de G-1}.}  Let us set $(X^n,Y^n,P^n,\beta^n,H^n):=$ $(X^{\phi_n}_{t_n,x_n},$ $Y^{\phi_n}_{t_n,x_n,y_n},$ $P^{\alpha_n}_{t_n,p_n},$ $ \beta_{t_n,x_n}(T),$ $H_{t_n,x_n})$ for $y_n:=v(t_n,x_n,p_n)\p 1/n$ and $(\phi_n,\alpha_n)\in \Ac_K\x L^2_{\Pc}$ such that 
$$Y^n(T)\ge \Psi^{-1}(X^n(T),P^n(T))\;.$$  Then, by the supermartingale property of $H^n\beta^nY^n$, one has 
$$
y_n\ge\Esp{H^n(T) \beta^n(T)\Psi^{-1}(X^n(T),P^n(T))}
$$
and, by choosing $(t_n,x_n,p_n)_n$ such that  $v(t_n,x_n,p_n)\to v_*(T,x,p)$, we obtain
\b*
v_*(T,x,p)&\ge& \liminf_{n\to \infty} \Esp{H^n(T) \beta^n(T)\Psi^{-1}(X^n(T),P^n(T))}\\
&=& \widehat \Psi^{-1}(x,p) \p \liminf_{n\to \infty} \delta_n\;,
\e*
where 
$$
\delta_n:= \Esp{H^n(T)\beta^n(T) \Psi^{-1}(X^n(T),P^n(T))}-\widehat \Psi^{-1}(x,p)\;.
$$
It remains to show that $\liminf_{n}\delta_n\ge  0$. To see this, first observe that $\Psi^{-1}\ge \widehat \Psi^{-1}$ so that 
$$
\delta_n\ge \Esp{H^n(T) \beta^n\widehat \Psi^{-1}(X^n(T),P^n(T))-\widehat \Psi^{-1}(x,p)}\;.
$$
Moreover, by convexity of $ \widehat \Psi^{-1}$, we have 
$$
  \widehat \Psi^{-1}(X^n(T),P^n(T))\ge   \widehat \Psi^{-1}(X^n(T),p)\p  \nabla^\p    \widehat \Psi^{-1}(X^n(T),p)(P^n(T)-p)
$$
so that 
\b*
  \delta_n&\ge& \Esp{ \nabla^\p    \widehat \Psi^{-1}(x,p)P^n(T)- H^n(T)\beta^n(T)  \nabla^\p    \widehat \Psi^{-1}(X^n(T),p)p }\\
  &&- \Esp{|H^n(T)\beta^n(T) \widehat \Psi^{-1}(X^n(T),p)-\widehat \Psi^{-1}(x,p)|}
  \\
  && -|\Oc|_{\infty}  \Esp{ |H^n(T)  \beta^n(T)\nabla^\p    \widehat \Psi^{-1}(X^n(T),p)-\nabla^\p    \widehat \Psi^{-1}(x,p)|}\;,
\e*
where $|\Oc|_{\infty}:=\max\{|q|,\;q\in \Oc\}<\infty$. Since $P^n$ is a martingale (under $\P$), this implies
\b*
  \delta_n&\ge& -\Esp{ |\nabla^\p    \widehat \Psi^{-1}(x,p)p_n- H^n(T) \beta^n(T) \nabla^\p    \widehat \Psi^{-1}(X^n(T),p)p| }\\
  &&- \Esp{|H^n(T) \beta^n(T)\widehat \Psi^{-1}(X^n(T),p)-\widehat \Psi^{-1}(x,p)|}
  \\
  && -|\Oc|_{\infty}  \Esp{ |H^n(T)  \beta^n(T)\nabla^\p    \widehat \Psi^{-1}(X^n(T),p)-\nabla^\p    \widehat \Psi^{-1}(x,p)|}\;,
\e*
and the required result follows from the assumptions of the proposition. 
\ep

%%%
\subsection{Discussion of the boundary condition on $\partial \Oc$}\label{sec: discussion cond bord sur P}

Since $\Oc$ is convex, it takes the form $[m,M]$ with $M,-m\in (-\infty,\infty]$. Obviously the boundary condition is meaningful only when $M<\infty$ or $m<\infty$. 

In order to recover a minimum of structure, we    impose the following conditions:
\be\label{eq: cond G pour M m}
\Psi(y,x)\ge M \Longrightarrow y\ge g(x) &\mbox{and} &\Psi(0,x)\ge m \;\;,\;\forall \;(x,y)\in \R^d\x \R\;
\ee
for some continuous function $g$. 
\vs2

For $\Psi(x,y)=\1_{y\ge g(x)}$, which corresponds to the quantile hedging problem, this holds for $M=1$ and $m=0$.  For $\Psi(x,y)=-\ell((g(x)-y)^\p)$, which corresponds to the expected loss pricing rule for $\ell$ convex non-decreasing, this holds with $M=-\ell(0)$ and $m=-\ell(\infty)=-\infty$.

\vs2

When $M$ is finite, \reff{eq: cond G pour M m} implies that $v(t,x,M)$ coincides with the super-hedging price of $g(X_{t,x}(T))$. When $m$ is finite, \reff{eq: cond G pour M m} implies that $v(t,x,m)=0$. 
Moreover, it is clear that $v$ is non-decreasing in the $p$-variable. It follows that 
$$
v_*(t,x,M)\le v^*(t,x,M)\le v(t,x,M) \mbox{ and } v^*(t,x,m)\ge v_*(t,x,m)\ge v(t,x,m)\;.
$$
However, equality may fail in the above inequalities. 

\vs2

This point being highly technical, we shall not discuss it further here. We refer to \cite{BET08} for {natural} conditions, which are typically satisfied in financial applications, under which equality holds. In particular, it is the case in the two examples of application below. 

%%%%%%%%%%%
\section{Example 1: Quantile hedging and Follmer-Leukert's formula}\label{sec: example cible en moment 1}

\subsection{Supersolution characterization of the quantile hedging price} 

In this section, we specialize the discussion to the quantile hedging problem of F\"ollmer and Leukert \cite{FoLe99}, which we already 
discussed in Chapter \ref{CHAP: quantile et loss par dualite}. We consider the non-constrained case $K=\R$ and we restrict to the one dimensional Black and Scholes model for ease of notations, see \cite{BET08} for a more general setting. 

It means that
 \be\label{musigma}
 \mu(x,a)=x\mu 
 ~\mbox{and}~
 \sigma(x,a)=x\sigma
 \ee
 where $\mu$ and $\sigma>0$ are now constants.  We fix $\rho=0$ for simplicity. 
 \vs2
 
  Then, the coefficients of the wealth process $Y$ are given by
 \be\label{eq coeff Y FK}
 \mu_Y(x,y,a)=ax \mu\;,\;
 \sigma_Y(x,a)=a x\sigma\;.
 \ee
Finally, we take 
 \be\label{Gg}
 \Psi(x,y)=\1_{\{y-g(x)\ge 0\}}
 &\mbox{for some Lipschitz function}&
 g:\R\longrightarrow\R_+.
 \ee
  
The stochastic target problem $v(t,x,p)$ corresponds to the problem of super-hedging the contingent claim $g(X_{t,x}(T))$ with probability $p$.\\

 Note that the above assumptions ensure that $v(\cdot,1)$ is continuous and is given by $v(t,x,1)=\E^{\Q_{t,x}}\left[g(X_{t,x}(T))\right]$ where $\Q_{t,x}$ is the $\P$-equivalent martingale measure defined by 
$$
d\Q_{t,x}/d\P=\exp\left(-\frac{T}{2}   |\lambda  |^2   - \lambda    W_T \right)\;,\;\lambda:=\mu/\sigma\;.
$$
For later use, let us denote by $W^{\Q_{t,x}}:=W-W_t+  \lambda (\cdot-t)$ the $\Q_{t,x}$-Brownian motion defined on $[t,T]$. 

\vs2

In Chapter \ref{CHAP: quantile et loss par dualite}, we have solved the quantile hedging problem   by means of the Neyman and Pearson's lemma from mathematical statistics. We shall see here how we can recover this result in the Markovian setting.  

\vs2

Note that, in this particular model, we have 
$$
\psi(x,p,q,b)=p\p bq/(x\sigma)\;\mbox{ and }\;\delta_K=\infty \mbox{ on } \R\setminus\{0\}\;,
 $$
so that Theorem \ref{thm: cara visco cible en moment} implies that $v_*$ should be a viscosity supersolution on $[0,T)\x (0,\infty)\x (0,1)$ of 
\be\label{eq pde bar v lower star FK}
-\partial_t \vp-\frac{\sigma^2x^2}{2}D_{xx}\vp - \inf_{b \in \R} \left(-b D_p\vp \lambda +  x\sigma b D_{xp}\vp +\frac{b^2}{2}D_{pp} \vp\right)
\ge 0\;.
 \ee 
Here, the conditions of Theorem \ref{thm: cara visco cible en moment} are not satisfied because $K$ is not compact, and we have omitted the condition $D_p\vp\ne 0$. However, the above holds for test functions such that $D_{pp}\vp(t_0,x_0,p_0)>0$ at the point where the minimum is achieved. The reason for this is that it allows to recover compactness on the set of $b$'s on which the above infimum is taken. This provides the required continuity on the operator associated to the above PDE in a neighborhood of $ D_{pp}\vp(t_0,x_0,p_0)$. Using this continuity, the proof  of Theorem \ref{thm: charact visco cible pas} in Chapter \ref{CHAPTER: cible pas} can be reproduced without difficulty. 
\vs2

Moreover, the conditions of Proposition \ref{prop: cond en T par env convexe de G-1} trivially hold with $\Psi^{-1}(x,p)=g(x)\1_{p>0}$, whose convex envelope in the $p$-variable is given by  $\widehat \Psi^{-1}(x,p)=pg(x)$. 
It follows that 
\be\label{eq: cond bord FollLeuk}
 v_*(T,x,p)\ge pg(x)\;.
\ee

\subsection{Formal explicit resolution}

The key idea for  solving \reff{eq pde bar v lower star FK}-\reff{eq: cond bord FollLeuk}  is to introduce the Legendre-Fenchel dual function of $v_*$   with respect to the $p-$variable in order to remove the non-linearity in \reff{eq pde bar v lower star FK}:
 \be\label{defU}
 w(t,x,q)
 :=
 \sup_{p\in\R}\left\{pq-  v_*(t,x,p)\right\} \,,
 & &
 (t,x,q)\in[0,T]\times(0,\infty)\times\R\;.
 \ee 
Note that  
 \be\label{defUbis}
w(\cdot,q)
=
\infty\,
\mbox{ for } q<0 \;\And\;
 w(\cdot,q)
 =
 \sup_{p\in[0,1]}\left\{pq-v_*(\cdot,p)\right\} \,
 \mbox{  for } q>0\;,
 \ee 
since 
\be\label{eq: cond ext dom v_*}
v_*\ge 0, \;v_*(\cdot,p)=0\; \mbox{ for $p< 0$ and $v_*(\cdot,p)=\infty$ for $p>1$,  }
\ee
by construction. One can actually  show, see \cite{BET08}, that 
\be\label{eq: cond bord p sur cible proba pour FollLeuk}
v_*(t,x,1)=v(t,x,1) \mbox{ and } v_*(t,x,0)=0\;,
\ee
recall the discussion of Section \ref{sec: discussion cond bord sur P}. 
\vs2

Using the PDE characterization of $v_*$  above, we shall prove below that  $w$ is an upper-semicontinuous viscosity subsolution on $[0,T)\times(0,\infty)\times (0,\infty)$ of 
 \be\label{eq Ustar}
 - \partial_t w-\frac{x^2\sigma^2}{2} D_{xx}w
     -\frac{\lambda^2 q^2}{2} D_{qq}w
     -x\sigma\lambda D_{xq}w
 &\le&
 0\;
 \ee
 with the boundary condition 
 \be\label{boundary-u}
 w(T,x,q)\le \left(q-g(x)\right)^+\;.
 \ee
Recalling the Feynman-Kac representation and comparison results of Theorems \ref{thm visco FK} and \ref{thm comp FK} of Chapter \ref{CHAP: pricing equation complet}, this implies that 
 \be
 w(t,x,q)
 \le  \bar  w(t,x,q):=\E^{\Q_{t,x}}\left[ \left(Q_{t,x,q}(T)-g\left(X_{t,x}(T)\right)\right)^+ \right]
\;,
 \ee
on $[0,T]\times(0,\infty)\times (0,\infty)$, where the process $Q_{t,x,q}$ is defined by the dynamics
 \be\label{Q}
 \frac{dQ(s)}{Q(s)} &=& \lambda  dW^{\Q_{t,x}}_s \;\;,\; Q_{t,x,q}(t)=q\in (0,\infty)\;.
 \ee 

Given the explicit representation of $ \bar  w$, we can now provide a lower bound to $v_*$ by using \reff{defUbis}. 

Clearly the function $\bar  w$ is   convex in $q$ and  there is a unique solution  $\bar q$ to the equation
 \be
\frac{\partial  \bar  w}{\partial q}\left(t,x, \bar q\right)
&=&
 \E^{\Q_{t,x}}\left[Q_{t,x,1} (T)\1_{\left\{  Q_{t,x, \bar q }(T)\ge g\left(X_{t,x}(T)\right)\right\}}\right]\nonumber\\
&=&\Pro{ Q_{t,x, \bar q }(T)\ge g\left(X_{t,x}(T)\right)}\nonumber\\
&=&
 p\;,\label{q(p)}
\ee
where we have used the fact that $d\P/d\Q_{t,x}=Q_{t,x,1}(T)$. 
It follows that the value function of the quantile hedging problem $v$ admits the lower bound
 \b*
 v(t,x,p) 
 &\ge &
 p \bar q - \bar  w\left(t,x, \bar q\right) 
  \\
 &=&
 \bar q\left[ p-   \E^{\Q_{t,x}}\left[ Q_{t,x,1}(T) 	
                \1_{\left\{ \bar qQ_{t,x,1}(T)\ge g\left(X_{t,x}(T)\right)\right\}}
          \right]\right]
 \\
 &+&  \E^{\Q_{t,x}}\left[g\left(X_{t,x}(T)\right)
                \1_{\left\{ \bar qQ_{t,x,1}(T)\ge g\left(X_{t,x}(T)\right)\right\}}
          \right]
 \\
 &=&
 \E^{\Q_{t,x}}\left[g\left(X_{t,x}(T)\right)
                \1_{\left\{ \bar qQ_{t,x,1}(T)\ge g\left(X_{t,x}(T)\right)\right\}}
          \right]=:\bar y\;.
 \e*
 On the other hand, it follows from  the martingale representation theorem, see Corollary \ref{cor: form dual prix couverture complet} of Chapter \ref{CHAP: prix dual}, that we can find $\phi \in \Ac_b$ such that 
 	$$
	Y^{\phi}_{t,x,\bar y}(T)\ge g\left(X_{t,x}(T)\right)
                \1_{\left\{ \bar qQ_{t,x,1}(T)\ge g\left(X_{t,x}(T)\right)\right\}} \;.
 	$$
Since,   $\Pro{  \bar qQ_{t,x,1 }(T)\ge g\left(X_{t,x}(T)\right)}$ $=  p$ by \reff{q(p)}, this implies that  $v(t,x,p)=\bar y$,
which corresponds exactly to the solution found in Chapter \ref{CHAP: quantile et loss par dualite}.

\subsection{Rigorous PDE characterization of the Fenchel-Legendre transform}

To conclude our argument, it remains to prove that $w$ is a   viscosity subsolution of 
\reff{eq Ustar}-\reff{boundary-u}. 

\vs2

 First note that the fact that $w$ is upper-semicontinuous on $[0,T]\x (0,\infty)\x (0,\infty)$ follows from the lower-semicontinuity of $v_*$ and the representation  in the right-hand side of \reff{defUbis}, which allows to reduce the computation of the sup  to the compact set $[0,1]$. Moreover, the boundary condition \reff{boundary-u} is an immediate consequence of   \reff{eq: cond bord FollLeuk} and \reff{eq: cond ext dom v_*}.
\\
 We now turn to the PDE characterization inside the domain. 
 Let $\vp$ be a smooth function with bounded derivatives and $(t_0,x_0,q_0)\in [0,T)\x(0,\infty)\x (0,\infty)$ be a local maximizer of $w-\vp$ such that $(w-\vp)(t_0,x_0,q_0)=0$.
\\ 
a. We first show that we can reduce to the case where the map $q\mapsto \vp(\cdot,q)$ is strictly convex. Indeed, since $w$ is convex, we necessarily have $D_{qq}\vp(t_0,x_0,q_0)\ge 0$. 
Given $\eps,\eta>0$, we now define $\vp_{\eps,\eta}$ by $\vp_{\eps,\eta}(t,x,q):=\vp(t,x,q)+\eps|q-q_0|^2+\eta|q-q_0|^2(|q-q_0|^2+|t-t_0|^2+|x-x_0|^2)$. Note that  $(t_0,x_0,q_0)$ is still a local maximizer of $w-\vp_{\eps,\eta}$. Since $D_{qq}\vp(t_0,x_0,q_0)\ge 0$, we have $D_{qq}\vp_{\eps,\eta}(t_0,x_0,q_0)\ge 2\eps>0$. 
Since $\vp$ has bounded derivatives, we can then choose $\eta$ large enough so that $D_{qq}\vp_{\eps,\eta}>0$. We next observe that, if $\vp_{\eps,\eta}$ satisfies \reff{eq Ustar} at $(t_0,x_0,q_0)$ for all $\eps>0$, then \reff{eq Ustar} holds for $\vp$ at this point  too. This is due to the fact that the derivatives up to order two of $\vp_{\eps,\eta}$ at  $(t_0,x_0,q_0)$ converge to the corresponding derivatives of $\vp$ as $\eps \to 0$.
\\
b. From now on, we thus assume that the map $q\mapsto \vp(\cdot,q)$ is strictly convex.   Let $\tilde \vp$ be the Fenchel transform of $\vp$ with respect to $q$, i.e. 
$$
\tilde \vp(t,x,p):=\sup_{q\in \R} \{pq-\vp(t,x,q)\}\,.
$$ 
Since $\vp$ is strictly convex in $q$ and smooth on its domain,   $\tilde \vp$ is strictly convex in $p$ and smooth on its domain, see e.g. \cite{R70}.  Moreover, we have 
	\be
	\vp(t,x,q)&=&\sup_{p\in \R} \{pq-\tilde \vp(t,x,p)\} \nonumber
	\\
	&=& J(t,x,q)q-\tilde \vp(t,x,J(t,x,q))\;\mbox{ on }   (0,T)\x (0,\infty)\x (0,\infty) \;\;\;\;\;\;\label{eq lien vp tilde vp J}
	\ee
where  $q\mapsto J(\cdot,q)$ denotes the inverse of $p\mapsto D_p \tilde \vp(\cdot,p)$, recall that $\tilde \vp$ is strictly convex in $p$.

We now deduce from the assumption  $q_0>0$ and \reff{defUbis} that we can find $p_0\in [0,1]$ such that 
$v(t_0,x_0,q_0)=p_0q_0-v_*(t_0,x_0,p_0)$ which,  
by using the  very definition of $(t_0,x_0,p_0,q_0)$ and $w$, implies that 
\be\label{eq txp local min}
\mbox{$(t_0,x_0,p_0)$ is a local minimizer of $v_*-\tilde \vp$ such that $(v_*-\tilde \vp)(t_0,x_0,p_0)=0$}
\ee 
and  	\be\label{eq vp tilde vp}
	&\vp(t_0,x_0,q_0)=\sup_{p\in \R} \{pq_0-\tilde \vp(t_0,x_0,p)\}=p_0q_0-\tilde \vp(t_0,x_0,p_0)\;&\\
	& \mbox{ with } p_0=J(t_0,x_0,q_0)&\nonumber
	\ee
where the last equality follows from \reff{eq lien vp tilde vp J} and the strict convexity of the map $p\mapsto pq_0-\tilde \vp(t_0,x_0,p)$ in the domain of $\tilde \vp$.

\vs2 

We conclude the proof by discussing three alternative cases depending on the value of $p_0$. 
\\
1. If $p_0\in (0,1)$, then \reff{eq txp local min} implies that $\tilde \vp$ satisfies \reff{eq pde bar v lower star FK} at $(t_0,x_0,p_0)$ and the required result follows by exploiting the  link between the derivatives of $\tilde \vp$ and  the derivatives of  its $p$-Fenchel transform $\vp$, which can be deduced from \reff{eq lien vp tilde vp J}.
\\
2.  If $p_0=1$, then the first boundary condition  in \reff{eq: cond bord p sur cible proba pour FollLeuk} and \reff{eq txp local min} imply that $(t_0,x_0)$ is a local minimizer of $v_*(\cdot,1)-\tilde \vp(\cdot,1)=v(\cdot,1)-\tilde \vp(\cdot,1)$ such that $(v(\cdot,1)-\tilde \vp(\cdot,1))(t_0,x_0)=0$. 
This implies that $\tilde \vp(\cdot,1)$ satisfies \reff{eq: PD FK} of Chapter \ref{CHAP: pricing equation complet} at  $(t_0,x_0)$, so that   $\tilde \vp$ satisfies \reff{eq pde bar v lower star FK} for $b=0$  at 
 $(t_0,x_0,p_0)$. We can then conclude as in 1. above. \\
3. If $p_0=0$, then the second boundary condition  in \reff{eq: cond bord p sur cible proba pour FollLeuk} and \reff{eq txp local min} imply  that $(t_0,x_0)$ is a local minimizer of $v_*(\cdot,0)-\tilde \vp(\cdot,0)=0-\tilde \vp(\cdot,0)$ such that $0-\tilde \vp(\cdot,0)(t_0,x_0)=0$. 
In particular, $(t_0,x_0)$ is a local maximum point for $\tilde \vp(\cdot,0)$ so that $(\partial_t \tilde \vp,D_x \tilde \vp)(t_0,x_0,0)=0$ and $D_{xx} \tilde \vp (t_0,x_0,0)\le 0$. This implies that $\tilde \vp(\cdot,0)$ satisfies  \reff{eq pde bar v lower star FK} at $(t_0,x_0,p_0)$, for $b=0$. We can then argue as in the first case.
  \ep
 
 %%%%%%%%%%%
\section{Example 2: Expected shortfall}\label{sec: example cible en moment}
 
Let us now consider the same model as above but with a risk constraint expressed through a quadratic loss function  as in Section \ref{Section: expected shortfall dualite} of Chapter \ref{CHAP: quantile et loss par dualite} (more general loss functions could obviously be considered, up to more tricky computations).

This corresponds to
$$
\Psi(x,y)=-((g(x)-y)^\p)^2\;,
$$
so that 
$$
\Psi^{-1}(x,p)=(g(x)-\sqrt{-p})^\p \;\;\mbox{ for } p\le 0\;.
$$ 

As in the previous section, we obtain that, for  any test function $\vp$ and $(t,x,p)\in [0,T)\x(0,\infty) \x (-\infty, 0)$ that achieves a minimum of $v_*-\vp$ and such that $D_{pp}\vp(t,x,p)>0$, one has 
\be\label{eq pde bar v lower star FK quantile}
-\partial_t \vp-\frac{\sigma^2x^2}{2}D_{xx}\vp - \inf_{b \in \R} \left(-b D_p\vp \lambda +  x\sigma b D_{xp}\vp +\frac{b^2}{2}D_{pp} \vp\right)
\ge 0\;.
 \ee 
By the same arguments as above, one can also show  that the Fenchel-Legendre transform
 \be\label{defU expected loss}
 w(t,x,q)
 :=
 \sup_{p\in\R}\left\{pq-  v_*(t,x,p)\right\}=\sup_{p\in(-\infty,0]}\left\{pq-  v_*(t,x,p)\right\} \,,
 \ee 
satisfies \reff{eq Ustar} on $[0,T)\x(0,\infty)\x (0,\infty)$, in the viscosity sense. As for the terminal condition, we obtain 
$v_*(T,x,p)\ge (g(x)-\sqrt{-p})^\p$, so that 
$$
w(T,x,q)= \left((4q)^{-1}-g(x)\right)\1_{\{(2q)^{-1}\le g(x)\}} \p (-qg(x)^2)\1_{\{(2q)^{-1}>g(x)\}}=:W(x,q)\;.
$$ It follows that 
$$
w(t,x,q)\ge \bar  w(t,x,q):=\E^{\Q_{t,x}}\left[ W(X_{t,x}(T),Q_{t,x,q}(T))  \right]\;,
$$
and therefore
\b*
v(t,x,p)&\ge& \sup_{q>0} \left(qp-  \E^{\Q_{t,x}}\left[ W(X_{t,x}(T),Q_{t,x,q}(T))\right]\right).
\e*
Direct computations combined with the identity  $d\P/d\Q_{t,x}=Q_{t,x,1}(T)$ then show that the optimum in the right-hand side term is achieved by 
$\bar q>0$ such that
\b*
	-p&=&-\partial_q \bar  w(t,x,\bar q)\\
	&=&\E\left[  (2Q_{t,x,\bar q}(T))^{-2}\wedge g(X_{t,x}(T))^2\right]\;.
\e*
Combining the above assertions   implies   that
\b*
v(t,x,p)&\ge& \E^{\Q_{t,x}}\left[  \left( g\left(X_{t,x}(T)\right)- (2Q_{t,x,\bar q}(T))^{-1}\right)^\p\right]=:\bar y\;.
\e*
 On the other hand, it follows from  the martingale representation theorem, see Corollary \ref{cor: form dual prix couverture complet} of Chapter \ref{CHAP: prix dual}, that we can find $\phi \in \Ac_b$ such that 
 $$
 Y^\phi_{t,x,\bar y}(T)=\left(g\left(X_{t,x}(T)\right)- (2Q_{t,x,\bar q}(T))^{-1}\right)^\p
 $$ 
 which, by the above identity, satisfies 
 $$
 \E[(g\left(X_{t,x}(T)\right)- Y^\phi_{t,x,y(t,x,p)})^\p)^2]=\E\left[  (2Q_{t,x,\bar q}(T))^{-2}\wedge g(X_{t,x}(T))^2\right]=p\;.
 $$ 
This shows that 
$$
v(t,x,p)=\E^{\Q_{t,x}}\left[  \left(g\left(X_{t,x}(T)\right)- (2Q_{t,x,\bar q}(T))^{-1}\right)^\p\right]\;,
$$
which is the result obtained in Section \ref{Section: expected shortfall dualite} of Chapter \ref{CHAP: quantile et loss par dualite}.
%%%%%%%%%%%%%%%%%%%%%%%%%%%%%%%%%%%%%%%%%%%%%%%%%%%%%%%%%%%%%%%%%%%%%%%%%%%%%%%%%%%%%%%%%%%%%%%%%%%%%%%%%%
\section{Example 3: Optimal book liquidation}

\subsection{Problem formulation and reduction}

The optimal book liquidation problem is the following. A financial agent asks a broker to sell on the market a total of $1$ stock on a time interval $[0,T]$, $1$ is taken as a normalization in order to save notations. The broker  takes the engagement that he will obtain a mean selling price which corresponds to (at least) a value $K>0$. The financial agent pays to the broker a premium $y$ at time $0$.
\vs2

The   cumulated number of stocks sold by the broker on the market since time $0$ is described by a continuous real-valued non-decreasing process $L$, we denote by $\Lc$ the set of such processes. Given $L\in \Lc$, the dynamic of the broker's portfolio $Y^L$   is given by 
$$
dY^L(t)=X^{L,1}(t)dL_t\;,\; Y^L(0)=y
$$
where $X^{L,1}$ represents the stock's price dynamics and is assumed to solve
\begin{eqnarray*}
  dX^{L,1}(t) &=&X^{L,1}(t) \mu(t,X^{L,1}(t))dt + X^{L,1}(t)\sigma(t,X^{L,1}(t))dW_t\\
  &&- X^{L,1}(t) \beta(t,X^{L,1}(t))dL_t
\end{eqnarray*}
where  $\mu,\sigma,\beta :[0,T]\x \R\mapsto \R$   are continuous
functions satisfying 
\be
&\beta \ge 0 \mbox{ , } (t,x^1) \in [0,T]\x (0,\infty) \mapsto x^1\left(\mu(t,x^1),\sigma(t,x^1),\beta(t,x^1)\right) \mbox{ is Lipschitz}\;&\nonumber\\
&\mbox{ and } x^1 \in  (0,\infty) \mapsto x\beta(t,x) \mbox{ is $C^2$ locally uniformly in $t \in [0,T]$.}&\label{eq: order book corff lipsch}
\ee
Note that we allow  the trading strategy of the broker to have an impact on the price dynamics if $\beta \ne 0$. 

In order to keep track of the cumulated number of units of asset already sold on the market, we introduce the process $X^{2,L}$ defined by the dynamics
$$
dX^{2,L}(t)=dL_t\;.
$$

\vs2

The aim of the broker is then to  find the initial premium $y=Y^L(0)$ and $L\in \Lc$  with $L_0=0$ such that $Y^L(T)  \ge K$ and $X^{2,L}(T)=1$, given that  $X^{2,L}(0)=0$.

In practice, it is clear that the above problem does not make sense and need to be relaxed. We shall therefore consider problems of the form 
\b*
&\mbox{Find $L\in \Lc$ with $L_0=0$ and $Y^L(0)$ s.t. } X^{2,L}(T)\le 1&\\
&\mbox{ and }   \Esp{\Psi(X^L(T),Y^L(T))} \ge p,&
\e* 
where
$$
\Psi(x,y):=\ell\left(y\p [x^1-x^1\beta(T,x^1)(1-x^2)](1-x^2)-K\right)\;,
$$
for $p\in  \R$ and   $\ell:\R\mapsto  \R$  is     (strictly) increasing with polynomial growth and is such that $\ell(\R)=\R$.\\
 The term  $X^{1,L}(T)[1-$ $\beta(T,X^{1,L}(T))(1-X^{2,L}(T))]$ $(1-X^{2,L}(T))$ stands for the gain of the final transaction required in order to liquidate the last units of assets at time $T$ if $X^{2,L}(T)<1$. Note that this final transaction is evaluated at the price $X^{1,L}(T)-X^{1,L}(T)\beta(T,X^{1,L}(T))(1-X^{2,L}(T))$ which already includes a possible depreciation of the stock's value due to this final trade. Obviously more sophisticated models could be considered within a similar framework.
\vs2

In order to define the associated value function, we now extend the above dynamics to arbitrary initial conditions. Given $L\in \Lc$, we write $Z_{t,x,y}^L=(X_{t,x}^L,Y_{t,x,y}^L)$, with $X_{t,x}^L=(X_{t,x}^{1,L},X_{t,x}^{2,L})$, the corresponding processes satisfying the initial condition $Z_{t,x,y}^L(t)=(x,y)=(x^1,x^2,y)$. 

The value function associated to the above stochastic target problem is then given by 
$$
v(t,x,p):=\inf\{y\in \R~:~\exists L\in \Lc\;\mbox{ s.t. }  X_{t,x}^{2,L}(T)\le 1\mbox{ and }  \Esp{\Psi(Z_{t,x,y}^L(T))} \ge p\}\;.
$$

As in the previous sections, we first convert   the above problem into a stochastic target problem. 

\begin{Proposition} For all $(t,x,p)\in [0,T]\x (0,\infty)\x [0,1]\x \R$, 
\b*
&v(t,x,p)=&
\\
&\inf\{y\in \R~:~ \exists (L,\alpha)\in \Lc\x L^2_\Pc\;\mbox{ s.t. }  X_{t,x}^{2,L}(T)\le 1\mbox{ and }  \Psi(Z_{t,x,y}^L(T)) \ge P_{t,p}^\alpha(T)\}\;,&
\e*
where 
$$
P^\alpha_{t,p}:=p\p \int_t^\cdot \alpha_s dW_s\;.
$$
\end{Proposition} 

\proof Since $\ell$ has polynomial growth, it is easily checked that $\Psi(Z_{t,x,y}^L(T))\in L^2$ for all initial condition and controls. It thus suffices to argue as in the proof of Proposition \ref{prop: reformulation cible en moment}. \ep

\subsection{PDE characterization in the domain}

The PDE characterization can be obtained by following the same arguments as in the proof of Theorem \ref{thm: cara visco cible en moment}. The main difference comes from the fact that the control $L$ is of bounded variation type which allows to play with local times in order to compensate for a lack of matching of the volatility terms. We shall come back on this important point in the proof, where it should be more clear. As a consequence, the PDE formulation is slightly different from the one obtained in the previous sections. 

In the following we denote by $v_*$ and $v^*$ the lower- and upper-semicontinuous envelopes of $v$ obtained by approximating by points in $[0,T)\x (0,\infty)\x [0,1) \x \R$.  

\begin{Theorem}\label{thm: pde dans domaine optimal liquidation} The function $v_*$ is a   viscosity supersolution on $[0,T)\x (0,\infty)\x [0,1) \x \R$ of 
\b*
 \max\{F_0\vp\;,\;x^1\p x^1\beta D_{x^1} \vp-D_{x^2}\vp  \;,\; -|D_p\vp|\}=0 
\e*
and 
the function $v^*$ is a   viscosity subsolution on $[0,T)\x (0,\infty)\x [0,1) \x \R$ of 
\b*
 \max\left\{\min\{F_0\vp\;,\;|D_p\vp|\}\;,\;x^1\p x^1\beta D_{x^1} \vp  -D_{x^2} \vp\right\}=0\;, 
\e*
where 
$$
F_0\vp:=-\Lc_X \vp -\frac{(x^1\sigma)^2}{2}\left( |D_{x^1} \vp/D_p\vp |^2 D^2_p \vp - 2(D_{x^1} \vp/D_p\vp)  D^2_{(x^1,p)} \vp  \right)\;,
$$ 
with 
$$
\Lc_X\vp:=\partial_t \vp \p x^1\mu D_{x^1}\vp \p \frac12 (x^1)^2\sigma^2 D^2_{x^1} \vp \;.
$$
\end{Theorem}

\proof {\bf Supersolution property:} We only sketch the proof as it follows from the same line of arguments as in the proof of Theorem \ref{thm: cara visco cible en moment}.   Let   $(t_0,x_0,p_0)$ be a point in $[0,T)\x (0,\infty)\x [0,1) \x \R$ which achieves a strict minimum of $v_*-\vp$ (equal to $0$ as usual). Then, if 
\b*
\max\{F_0\vp\;,\;x^1\p x^1\beta D_{x^1} \vp-D_{x^2} \vp \}(t_0,x_0,p_0)<0\;\mbox{ and } |D_p\vp(t_0,x_0,p_0)|>0\;,
\e*
we can find $\eta,r>0$ such that 
\b*
&\max\{-\Lc_{X,P}^b\vp\;,\;x^1\p x^1\beta D_{x^1} \vp-D_{x^2} \vp\;,\; - |D_p\vp|\}\le -\eta\; &\\
&\mbox{ for $b\in \R$ s.t. } |bD_p\vp\p x^1\sigma D_{x^1} \vp|\le r&
\e*
on a neighborhood of $(t_0,x_0,p_0)$,
where 
$$
\Lc_{X,P}^b\vp:=\Lc_X\vp \p \frac{1}{2}\left( b^2 D^2_p \vp \p 2 x^1\sigma b  D^2_{(x^1,p)} \vp  \right)\;.
$$
It then suffices to reproduce the arguments of the proof of Theorem \ref{thm: charact visco cible pas} in Chapter \ref{CHAPTER: cible pas}. The fact $x^1\p x^1\beta D_{x^1} \vp-D_{x^2} \vp\le 0$ allows to forget about the non-decreasing control $L$ when applying It\^o's Lemma on the difference $Y^L-\vp(\cdot,X^L,P^\alpha)$.
\vs2

 {\bf Subsolution property:}  Let   $(t_0,x_0,p_0)$ be a point in $[0,T)\x (0,\infty)\x [0,1) \x \R$ which achieves a strict maximum of $v^*-\vp$ (equal to $0$ as usual).  We have to show that 
\b*
\max\{\min\{F_0\vp\;,\; |D_p\vp|\}\;,\;x^1\p x^1\beta D_{x^1} \vp-D_{x^2}\vp \}(t_0,x_0,p_0)\le 0\;.
\e*
If 
\b*
\min\{F_0\vp\;,\; |D_p\vp|\} (t_0,x_0,p_0)> 0\;
\e*
then a contradiction to the Geometric Dynamic Programming principle is obtained by considering a control of the form 
$$
(L,\alpha)=(0,-x^1\sigma D_{x^1} \vp/D_p\vp)
$$
and by arguing as in the proof of subsolution property of Theorem \ref{thm: charact visco cible pas} in Chapter \ref{CHAPTER: cible pas}.
We next discuss the case   where 
\b*
x^1_0\p x^1_0 \beta(t_0,x^1_0) D_{x^1} \vp (t_0,x_0,p_0)- D_{x^2} \vp (t_0,x_0,p_0)> 0\;.
\e*
We shall now see how we can play with the non-decreasing control $L$ is order to compensate the fact that the Brownian diffusion parts are possibly not matched in our dynamics. 

Indeed, if the above hold, one can find  $\eps,\eta>0$ small enough so that
 	\be
	 &  x^1\p x^1\beta(t,x^1) D_{x^1} \vp(t,x,p) -D_{x^2} \vp(t,x,p) >  \eta &\label{eq contra preuve soussolution 2}\\
	 & \forall\;   (t,x,p)\in B_\eps(t_0,x_0,p_0) ,\;  |y- \vp(t,x,p)|\le \eps\;.&\nonumber
	 \ee
Set $\Oc_n:=\{(t,x,y)~:~(t,x,p)\in B_{2\eps}(t_0,x_0,p_0)\;,\;|y- \vp(t,x,p)|<2\eps\;,\;y-\vp(t,x,p)>-|\gamma_n|\}$, where 
	\be\label{eq iotan to 0 soussol}
	\gamma_n:=y_n-\vp(t_n,x_n,p_n)\to 0\;,
	\ee
with    $(t_n,x_n,p_n)_{n\ge 1}$   a sequence which converges to $(t_0,x_0,p_0)$   such that $v(t_n,x_n,$ $p_n)$ $\to$ $ v^*(t_0,x_0,p_0)$, and  $y_n:=v(t_n,x_n,p_n)-n^{-1}$.
By a simple Taylor expansion of order $1$, we then deduce from \reff{eq contra preuve soussolution 2} that, for $(t,x,p)\in B_{\eps}(t_0,x_0,p_0)$ and $y\in \R$ such that $|y- \vp(t,x,p)|<\eps$, we have, for $0< \lambda \le r$ with $r>0$ small enough,  
$$
(y-\lambda x^1)-\vp(t,x^1\p \lambda x^1\beta(t,x^1),x^2-\lambda ,p)
\le -|\gamma_n| -\lambda \eta \p O(r^2)< -|\gamma_n|
$$ 
whenever $|(t,x,p,y)-(t',x',y',p')|\le \lambda r$ for some $(t',x',y',p')$ such that  $y'-\vp(t',x',p')=-|\gamma_n|$. Otherwise stated, the direction $-(x^1,-x^1\beta(t,x^1), 1)$ is driving $(y,x^1,x^2)$ strictly out of the smooth domain $\{(y,x^1,x^2)~:$ $y-\vp(t,x^1,$ $x^2,$ $p)$ $>$ $-|\gamma_n|\}$, at least locally around $(t_0,x_0,p_0)$ and $\vp(t_0,x_0,p_0)$.  This implies that it is possible to reflect the process $Y-\vp(\cdot, X,P)$ along an inward direction by suitably pushing $(Y,X^1,X^2)$ in the direction $(X^1,-X^1\beta(t,X^1), 1)$.

More precisely,   \reff{eq: order book corff lipsch}    and the above discussion allow to apply Theorem 4.8 of \cite{DuIs93}:  there exists a continuous real-valued adapted  non-decreasing process $ L^n$ satisfying
\be
       Y^n(s\wedge \theta^n)&\ge& \vp(s\wedge \theta^n,X^n(s\wedge \theta^n),p_n)-|\gamma_n|\;\;  \mbox{ for all } s\ge t_n\;,\label{eq: soussol int Y ge vp}
\ee
where 
	\b*
	\theta^o_n
	&:=&
	\inf\left\{s\ge t_n~:~(s,X^n(s),p_n)\notin B_\eps(t_0,x_0,p_0)\right\},
	\\
	\theta_n
	&:=&
	\inf\left\{s\ge t_n~:~|Y^n(s)-\vp(s,X^n(s),p_n)|\ge \eps\right\}\wedge \theta^o_n\;,
	\e*
and  $(X^n,Y^n):=(X^{L_n}_{t_n,x_n},Y^{L_n}_{t_n,x_n,y_n})$.

 In view of \reff{eq: soussol int Y ge vp} and  \reff{eq iotan to 0 soussol}, we have $ Y^n(\theta_n)-\vp(\theta_n,X^n(\theta_n),p_n)\ge-|\gamma_n|>-\eps$ for $n$ large enough. Following the arguments  of the  proof of the subsolution property of Theorem \ref{thm: charact visco cible pas} in Chapter \ref{CHAPTER: cible pas} with the control $(0,L^n)$ then leads to   the required contradiction to the Geometric Dynamic Programming Principle. 
\ep 

\subsection{Boundary conditions}

We first discuss the boundary condition   at $t=T$. By definition 
$$
v(T-,\cdot)=\Psi^{-1}
$$
where $\Psi^{-1}$ denotes the inverse of the function $\Psi$ with respect to the $y$-variable:
$$
\Psi^{-1}(x,p):=\ell^{-1}(p) - [x^1-x^1\beta(T,x^1)(1-x^2)](1-x^2)\p K\;.
$$
However, as in Theorem \ref {thm: condition en T cible pas} of Chapter \ref{CHAPTER: cible pas},  the gradient constraint $x^1\p x^1\beta D_{x^1} \vp-D_{x^2}\vp\ge 0$ that holds inside the domain should propagate up to the boundary. This leads to the following boundary condition. 

\begin{Theorem} The function  $v_*$ is a   viscosity supersolution on $\{T\}\x (0,\infty)\x [0,1) \x \R$ of 
\b*
 \max\{\vp-\Psi^{-1}\;,\;x^1\p x^1\beta D_{x^1} \vp-D_{x^2}\vp \;,\; -|D_p\vp|\}=0 
\e*
and 
the function $v^*$ is a   viscosity subsolution on $\{T\}\x (0,\infty)\x [0,1) \x \R$ of 
\b*
\min \left\{\vp-\Psi^{-1}\;,\; \max\{|D_p\vp|\;,\;x^1\p x^1\beta D_{x^1} \vp  -D_{x^2} \vp\}\right\}=0\;.
\e*
\end{Theorem}

\proof Combine the arguments of  the proof of Theorem \ref {thm: condition en T cible pas} of Chapter \ref{CHAPTER: cible pas} with the ones used in the proof of Theorem \ref{thm: pde dans domaine optimal liquidation}. \ep

\begin{Remark}{\rm  When $\beta$ is constant, i.e. does not depend on $x^1$, and $\ell$  is $C^1$, then one easily checks that $\Psi^{-1}$ is a strong super- and subsolution of the above equations. In this case, one can actually show  that $v_*(T,\cdot)=v^*(T,\cdot)=\Psi^{-1}$.
}
\end{Remark}

It remains to study the boundary condition for $x^2=1$. Note that, when $x^2=1$, the constraint $X^{2,L}_{t,x}(T)\le 1$ and the fact that $L-L_t$ is non-decreasing imply that $L=L_t$, $X^{2,L}_{t,x}=1$ and $Y^{L}_{t,x,y}=y$. Hence, 
$$
v(t,x^1,1,p)=\inf\{y\in \R~:~\Esp{\Psi(X^{1,0}_{t,x}(T),1,y) }\ge p \}=:\bar v(t,x^1,p)\;.
$$ 
Our last result shows that the function $\bar v$ is actually the correct boundary condition at $x^2=1$. 

\begin{Proposition} We have $v_*(t,x^1,1,p)=v^*(t,x^1,1,p)=\bar v(t,x^1,p)$ for all $(t,$$x^1, $ $p) \in [0,T]\x (0,\infty) \x \R$. 
\end{Proposition}

\proof Let   $(t_n,x_n,p_n)_{n\ge 1}$ be   a sequence which converges to $(t_0,x^1_0,1,p_0)$. Fix $y_n\in \R$ and $L^n\in \Lc$ such that   $Z^n=(X^n,Y^n):=(X^{L^n}_{t_n,x_n}, Y^{L^n}_{t_n,x_n,y_n})$ satisfies $\Esp{\Psi(Z^n(T))}\ge p_n$ and $X^{2,n}(T)\le 1$. Then, the last constraint combined with the Lipschitz continuity assumption on our coefficients implies that $L^n(T)-L^n(t_n)\to 0$ so that $Z^n(T)$ $\to$ $Z^0(T):=(X^{1,0}_{t_0,x^1_0,1}(T),1,y_0)$ uniformly in  $L^q$, for any $q\ge 2$. Since $\Psi$ has polynomial growth, it follows from the dominated convergence theorem that 
$$
\Esp{\Psi(Z^n(T)) }\to \Esp{\Psi(X^{1,0}_{t_0,x^1_0,1}(T),1,y_0)}\ge p_0\;,
$$
whenever $y_n\to y_0\in \R$. 
By choosing $(t_n,x_n,p_n)_{n\ge 1}$ such that $v(t_n,x_n,p_n)\to v_*(t_0,x^1_0,1,p_0)$ and $y_n=v(t_n,x_n,p_n)\p n^{-1}$, we thus deduce that $v_*(t_0,x^1_0,1,p_0)\ge \bar v(t_0,x^1_0,p_0)$. On the other hand, one can also choose $(t_n,x_n,p_n)_{n\ge 1}$ such that $v(t_n,x_n,p_n)\to v^*(t_0,x^1_0,1,p_0)$ and $y_n=y_0=\bar v(t_0,x^1_0,p_0)$ so that, for any $\eps>0$,
$$
 \Esp{\Psi(X^{0}_{t_n,x_n}(T),y_n\p \eps)}\ge p_0\;
$$
for $n$ large enough. This follows from the convergence   $\Esp{\Psi(X^{0}_{t_n,x_n}(T),y_n\p \eps)}\to \Esp{\Psi(X^{1,0}_{t_0,x^1_0,1}(T),1,y_0\p \eps)}> p_0$, where the last inequality is a consequence of  the definition of $y_0=\bar v(t_0,x^1_0,p_0)$ and the fact that $\ell$ is strictly increasing.
This shows that $\bar v(t_0,x^1_0,p_0)\p \eps \ge v(t_n,x_n,p_n)$ for $n$ large enough, and therefore that $\bar v(t_0,x^1_0,p_0)\p \eps \ge v^*(t_0,x^1_0,1,p_0)$. We conclude by arbitrariness of $\eps>0$.
\ep 

 %%%%%%%%%%%%%%%%%%%%%%%%%%%%%%%%%%%%%%%%%%%%%%%%%%%%%%%%%%%%%%%%%%%%%%%%%%%%%%%%%%%%%%%%%%%%%%%%%%%%%%%%%%
\part{Exercices}

\small
\setcounter{section}{0}

\section{Discrete time model}
\setcounter{equation}{0}

We consider a finite probability space $(\Omega, \Fc, \P)$ equiped with a complete filtration $\F = (\Fc_n)_{n\le N}$, where $N\ge 1$, satisfying $\Fc_0 = \{ \emptyset, \Omega\}$ and $\Fc_N = \Fc$. We consider a discrete time model, with a non-risky asset $B:=(B_n)_{1\le n\le N}$ and $d$ risky securities $X:=(X_n^1, \cdots, X_n^d)_{1\le n\le N}$ where $B$ and $X$ are $\F$-adapted.\\
The dynamics of $B$ is given by $B_n = (1+r_n)B_{n-1}$ for $n\ge 1$, with $r=(r_n)_{n\ge 1}$   a positive $\F$-predictable process and $B_0=1$. A portfolio strategy is defined by a $\F$-predictable process $(\alpha, \phi)$ taking its value in $\R \x \R^d$: $\alpha_n$ (respectively $\phi_n^i$) is the quantity of non-risky asset (respectively risky asset $X^i$) held in the porfolio on the time period $[n-1, n]$.
\begin{enumerate}
\item Dynamics of the non-risky asset.
\begin{enumerate}
	\item[1.1.] If one invest 1\$ at time $n-1$ in $B$, how many shares of asset $B$ is held?
	\item[1.2.] What is the value of this portfolio at time $n$?
\end{enumerate}
\item Dynamics of the portfolio.
\begin{enumerate}
	\item[2.1.] Write the value of the portfolio $Y_n^{x,(\alpha, \phi)}$ at time $n$, where $x$ is the initial value of the portfolio.
	\item[2.2.] Write the self-financing condition.
	\item[2.3.] Write the dynamics of the portfolio value of an investor in function of $\phi$. From now on, we will use the standard notation $Y^{x,\phi}$ in place of $Y^{x, (\alpha, \phi)}$.
	\item[2.4.] We now use the following notation: $\tilde X := X/B$ and $\tilde Y=Y/B$. Give the dynamics of $\tilde X$ and $\tilde   Y^{x,\phi}$.
	\item[2.5] Denote by $\Mc (\P)$ the set of measures $\Q \sim \P$ such that $\tilde X$ is a $(\Q, \F)$-martingale. Assume that $\Mc(\P)\neq \emptyset$. Show that $\tilde Y^{x,\phi}$ is a $(\Q, \F)$-martingale for every $\Q\in \Mc (\P)$.
	\item[2.6.] Let $G$ be a $\Fc$-measurable random variable. Give the super-replication price $p(G)$ of $G$.
\end{enumerate}
\end{enumerate}

\section{Portfolio optimization}
\setcounter{equation}{0}

We use the framework of the previous exercice. We suppose now that $r\equiv0$. We denote $\Ac$ the set of previsible processes taking value in $\R^d$. Here, we search for the solution of the utility maximization problem
$$\sup_{\phi \in \Ac} \Esp{U(V_N^{x,\phi})} \, ,$$
where $U$ is a function in $\Cc^1(\R)$, strictly increasing, strictly concave, defined on the whole real line, with the Inada conditions:
$$\lim_{x \rightarrow -\infty}U'(x) = +\infty\, , \quad \lim_{x \rightarrow +\infty}U'(x) = 0 \, .$$

\begin{enumerate}
	\item[1.] Let $\tilde U$ be the function defined by
	\beq\label{eq:fenchel}
	\tilde U(y) = \sup_{x\in \R}(U(x)-xy), \quad y>0 \, .
	\eeq
	\begin{enumerate}
		\item[1.1.] When is the supremum attained?
		\item[1.2.] Deduce that $\tilde U(y) = U(\hat x (y)) - \hat x(y) y$, with $y>0$ and $\hat x(y) = (U')^{-1}(y)$. We will further admit that $\tilde U$ is $\Cc^1$.
	\end{enumerate}
	\item[2.]We suppose now that $\Mc(\P) = \{ \Q \}$. 
	\begin{enumerate}
		\item[2.1.] Is the market complete?
		\item[2.2.] Give the hedging price at time $0$ of a contingent claim $G \in L^{\infty}(\R, \Fc_N)$.
	\end{enumerate}
	\item[3.] We denote $H:= d\Q / d\P$.
	\begin{enumerate}
		\item[3.1.] Compute $\E[H V_N^{x,\phi}]$ when $\phi \in \Ac$.
		\item[3.2.] Deduce from \reff{eq:fenchel} that for all $\phi \in \Ac$ and $\lambda>0$,
		$$\Esp{U(V_N^{x,\phi})}\le \Esp{\tilde U (\lambda H)} + \lambda x \, .$$
	\end{enumerate}
	\item[4.] We admit now that there exists some $\hat \lambda>0$ such that
	$$\inf_{\lambda>0} \left(\Esp{\tilde U (\lambda H)} + \lambda x \right) = \Esp{\tilde U (\hat \lambda H)} + \hat \lambda x \, .$$
	\begin{enumerate}
		\item[4.1.] Show that $\tilde U$ is convex and deduce that
		$$\Esp{H \tilde U' (\hat \lambda H)}+ x=0 \, .$$
		\item[4.2.] From the last result, show that there exists some $\hat \phi \in \Ac$ such that $V_N^{x,\hat \phi} = - \tilde U'(\hat \lambda H)$.
	\end{enumerate}
	\item[5.] We admit that $-\tilde U' = (U')^{-1}$. Deduce from all the results you achieved that 
	$$\sup_{\phi \in \Ac} \Esp{U(V_N^{x,\phi})} = \Esp{U(V_N^{x,\hat \phi})}$$
	
\end{enumerate}

\section{Continuous time model}
\setcounter{equation}{0}

We intend here to price and hedge a European call option of maturity $T$ and strike $K$, meaning a claim such that one receives  the payoff $(X_T - K)^+$ at time $T$. We put ourselves in the framework of the Black-Scholes model, and the price of this claim will be given by the initial value of a portfolio with a strategy $\phi$ such that it returns the wealth $Y^{y,\phi}_T = (X_T - K)^+$.\\
In this model, we assume that the dynamics of the price $X_t$ of the one dimensional risky asset is given by
$$dX_t = X_t (\mu dt + \sigma dW_t)$$
where $(W_t)_{t\in [0,T]}$ is a standard Brownian motion on the complete probability space $(\Omega, \Fc, \P)$ and with $\sigma>0$ $dt\x d\P$-a.e. Let $(\Fc_t)_{t\le T}$ denotes the filtration generated by $(W_t)_{t\ge 0}$.  Let also $B$ be the risk free asset defined by
$$B_t = 1+ \int_0^t B_s r_s ds$$
where $r$ is a predictable bounded real valued process. We denote by $Y^{y,\phi}$ the portfolio value associated to the initial value $y$ and $\phi \in \Ac_b$:
$$Y_t^{y,\phi} = y + \int_0^t  \phi_s dX_s\;.$$
Here, $\Ac_b$ denotes the set of   strategies $\phi$ such that the associated wealth process is bounded from below.
\begin{enumerate}
	\item[1.] Using It\^o's lemma,
	\begin{enumerate}
		\item[1.1.] Show that $B_t = e^{\int_0^t r_s ds}$. From now on, we will use the notation $\beta_t = B_t^{-1} = e^{-\int_0^t r_s ds}$.
		\item[1.2.] Write the dynamics of $\tilde X:=\beta X$ and $\tilde Y^{y,\phi}:=\beta Y^{y,\phi}$. 
	\end{enumerate}
	\item[2.] Using Girsanov Theorem write the equivalent martingale (or risk-neutral) measure $\Q$.
	\item[3.] We suppose from now on that $r$ is constant.
	\begin{enumerate}
		\item[3.1.] What does $\E^{\Q}[e^{-rT} \1_{X_T \ge K}]$ represent?
		\item[3.2.] Compute it.
		\item[3.3.] When $X_t = K$, what happens for $t\rightarrow T$? Give an interpretation. 
	\end{enumerate}
	\item[4.] Compute $\E^{\Q}[e^{-rT} (X_T - K)^+] = v(0,X_0)=y$.
	\item[5.] 	\begin{enumerate}
		\item[5.1.] What is the PDE satisfied by $v$?
		\item[5.2.] What is the hedging strategy of the claim $(X_T - K)^+$?  
	\end{enumerate}
\end{enumerate}
 
\section{Exchange option}
\setcounter{equation}{0}

Let $B, S^1$ and $S^2$ be three assets with the following dynamics

\beqq \begin{aligned}
    dB_t & = r B_t dB_t\\
    dS^1_t & = S^1_t \left( b^1_t dt + \sigma^1_1 dW^1_t + \sigma^1_2 d W^2_t \right)\\
    dS^2_t & = S^2_t \left( b^2_t dt + \sigma^2_1 dW^1_t + \sigma^2_2 d W^2_t \right)\\
\end{aligned} \eeqq

where $W^1$ and $W^2$ are two independents Brownian motion defined on $(\Omega,\Fc,\P)$. We assume furthermore that $\sigma = \left( \sigma^i_j \right)_{i,j = 1,2}$ is a deterministic non-singular matrix.
\begin{enumerate}

    \item Provide a probability $\P^*$ such that the both processes $S^1/B$ and $S^2/B$ are $\P^*$-martingale. Provide then the dynamics of these both processes under this probability.
       
    \item Define a change of probability from $\P^*$ to $\Q^*$ such that $S^1/S^2$ is a $\Q^*$-martingale. Provide the dynamics of $S^1/S^2$ under $\Q^*$.
   
    \item Deduce then the price of the option exchange which payoff is $\left( S^1_T - S^2_T \right)^+$.

\end{enumerate}

\section{Forward Option}
\setcounter{equation}{0}

Let us consider the   Black Scholes model where the  risky asset $S_t$ and the risk free asset $S^0_t$  have the dynamics
    \beqq
        \left\{
            \begin{aligned}
                dS_t & = S_t \left( \mu dt + \sigma dW_t \right)\\
                dS^0_t & = S^0_t r dt
            \end{aligned}
        \right.
    \eeqq
in which  $W$ is a Brownian motion under the historic probability $\P$, and $\sigma$ is invertible.
\begin{enumerate}

    \item Is this market complete ? If this is the case, give the risk neutral probability measure $\Q$.
   
    \item A forward option is an option, paid at time $t_0$, which gives at time $t_1$ an option of maturity $t_2$ and strike $S_{t_1}$. Write the price of this option as an expectation under $\Q$.
       
    \item Give the value at each date of a at the money call with forward strike. One will write this price as a classical call.

\end{enumerate}

\section{Gamma hedging}

We consider a probability space $(\Omega,\Fc,\P)$ equipped with a filtration $\Fc = \left( \Fc_t \right)_{t \in [0,\infty)}$ satisfying the usual conditions, and such that $\Fc_0$ is trivial.
Let $T > 0$, and $W$ being a Brownian motion on this space.
We consider a financial market with a risk free asset of return $r=0$, and with a risky asset which the price is $S = \left( S_t \right)_{t\ge0}$ is the unique strong solution of
$$
    S_t = S_0 + \int_0^t S_s \sigma (S_s) dW_s, \ \ \ \ t \ge 0.
$$

We define $a : x \in [0,\infty) \mapsto a (x) := x \sigma (x) \in [0,\infty)$, and we assume that $a$ is uniformly Lipschitz. Denote by $\Lc$ the Dynkin operator associated to this SDE, i.e.

$$
    \Lc \vp (t,x) := \frac{\partial}{\partial t} \vp (t,x) + \demi a (x)^2 \frac{\partial^2}{\partial x^2} \vp (t,x)
$$

for $\vp \in C^{1,2}$.

Denote by $\Ac$ the set of $\F$-predictable processes $\phi$ such that $\Esp{\int_0^T \left| \phi_s a (S_s) \right|^2 ds} < \infty$ for all $T > 0$, and, we assume that, for all $T>0$ and every random variable being $\Fc_T$-measurable $X$ such that $\Esp{\left| X \right|^2} < \infty$, there exists some $\phi \in \Ac$ such that $V_T^{\Esp{X},\phi} = X \ \P$-a.s., where $V^{x,\phi}_t := x + \int_0^t \phi_s dS_s, t \ge 0, (x,\phi) \in \R \x \Ac$.

\begin{enumerate}

    \item \emph{a priori} estimations

        \begin{enumerate}
            \item Show that, for every $T > 0$ and $p \ge 1$, there exists a constant $C_{T,p} > 0$ such that $\Esp{\sup_{t\le T} \left| S_t \right|^p} \le C_{T,p}$.

            \item Show that $S$ is a $\P$-martingale.
        \end{enumerate}

    \item Let $G$ be a Borel function with polynomial growth, and $T_2 > 0$.

        \begin{enumerate}

            \item Prove the existence of a function $g : [0,T_2] \x [0,\infty) \rightarrow \R$ such that $g (t,S_t) = \Esp{\left. G \left( S_{T_2} \right) \right| \Fc_t} \ \P$-a.s. when $t \le T_2$.

            \item Assuming that $g$ is smooth enough, what is the PDE satisfied by $g$ ?

            \item Assume now that $\Esp{\int_0^{T_2} \left| \frac{\partial}{\partial x} g (t,S_t) a (S_t) \right|^2 dt} < \infty$. What is the price of an option which payoff is $G \left(S_{T_2}\right)$ paid at time $T_2$ compatible with the no-arbitrage condition. Express the hedging strategy of this option in terms of partial derivatives of both $g$ and $a$.

        \end{enumerate}

    \item Consider now an other Borel function $F$ with polynomial growth and $0 < T_1 < T_2$. We assume that $g \in C^{1,2}_b \left( [0,T_1] \x [0,\infty) \right).$\footnote{The $_b$ means that the partial derivatives are bounded on the considered set.} For $x \in \R, \phi, \alpha \in \Ac$, define

        $$
            V^{x,\phi,\alpha}_t := x + \int_0^t \phi_s dS_s + \int_0^t \alpha_s dg \left(s,S_s\right) \ \ \ \ t \in [0,T_1].
        $$

        We assume that there exists $\bar{\phi}, \bar{\alpha} \in \Ac$ such that

        \be
            0 &=& \bar{\phi}_t + \bar{\alpha}_t \frac{\partial}{\partial x} g (t,S_t) - \frac{\partial}{\partial x} f (t,S_t)\label{eq1}\\
            &=& \bar{\alpha}_t \frac{\partial^2}{\partial x^2} g (t,S_t) - \frac{\partial^2}{\partial x^2} f (t,S_t) \ \ \ \ \P \text{-a.s. } \forall \ t < T_1\label{eq2},
        \ee

        where $f \in C^{1,2} \left( [0,T_1] \x [0,\infty) \right)$ satisfies $f (t,S_t) = \Esp{\left. F \left( S_{T_1} \right) \right| \Fc_t} \ \P$-a.s. for every $t \le T_1$.

            \begin{enumerate}

                \item Give a financial interpretation of $V^{x,\bar{\phi},\bar{\alpha}}$.

                \item Find $\bar{x} \in \R$ such that $\bar{V} := V^{\bar{x},\bar{\phi},\bar{\alpha}}$ satisfies $\bar{V}_{T_1} = F \left( S_{T_1} \right)$.

            \end{enumerate}

        \item Let $n \in \N \setminus \left\{ 0 \right\}$ and $t_i := iT_1/n, i \le n.$ Denote by $\eta_t := \max\left\{ t_i, i \le n \text{ s.t. } t_i \le t \right\}$, i.e. $\eta_t = t_i$ if $t \in [t_i,t_{i+1}), t \ge 0$.
            From now on, we consider the piecewise constant strategy $\left( \we{\phi} , \we{\alpha} \right)$ defined by $\left( \we{\phi}_t , \we{\alpha}_t \right) := \left( \bar{\phi}_{\eta_t} , \bar{\alpha}_{\eta_t} \right), t \le T_1$.
            Denote by $\we{V} := V^{\bar{x}, \we{\phi}, \we{\alpha}}$. For sake of simplicity, we assume furthermore that both $g$ and $f$ are $C^\infty$ with bounded derivatives, and that the process $\bar{\alpha}$ is essentially bounded\footnote{even if it is unrealistic}.

            \begin{enumerate}

                \item By using \ref{eq1}, show that

                \b*
                    \we{V}_{T_1} - F \left( S_{T_1} \right) &=& \int_0^{T_1} \bar{\alpha}_{\eta_t} \left( \frac{\partial}{\partial x} g \left( t,S_t \right) - \frac{\partial}{\partial x} g \left( \eta_t,S_{\eta_t} \right) \right) a \left( S_t \right) dW_t\\
                    && \ \ \ \ - \int_0^{T_1} \left( \frac{\partial}{\partial x} f \left( t,S_t \right) - \frac{\partial}{\partial x} f \left( \eta_t,S_{\eta_t} \right) \right) a \left( S_t \right) dW_t\\
                    &=& \int_0^{T_1} A_t a \left( S_t \right) dW_t,
                \e*

                where $A_t := \int_{\eta_t}^t B_s a \left( S_s \right) dW_s + \int_{\eta_t}^t C_s ds$ with

                \b*
                    B_s &:=& \bar{\alpha}_{\eta_s} \frac{\partial^2}{\partial x^2} g \left( s,S_s \right) - \frac{\partial^2}{\partial x^2} f \left( s,S_s \right)\\
                    C_s &:=& \bar{\alpha}_{\eta_s} \Lc \left[ \frac{\partial}{\partial x} g \left( s,S_s \right) \right] - \Lc \left[ \frac{\partial}{\partial x} f \left( s,S_s \right) \right]
                \e*

                \item Using \ref{eq1} again, show that

                \b*
                    B_s &:=& \int_{\eta_s}^s \left( \bar{\alpha}_{\eta_s} \frac{\partial^3}{\partial x^3} g (u,S_u) - \frac{\partial^3}{\partial x^3} f (u,S_u) \right) a (S_u) dW_u\\
                    && + \int_{\eta_s}^s \left( \bar{\alpha}_{\eta_s} \Lc \left[ \frac{\partial^2}{\partial x^2} g (u,S_u) \right] - \Lc \left[ \frac{\partial^2}{\partial x^2} f (u,S_u) \right] \right) a (S_u) du.
                \e*

                \item Show that there exists $C > 0$ such that $\Esp{\left| B_s^2 \right|} \le C/n^2$ for every $t \le T_1$.

                \item Deduce from the previous question that there exists $C>0$ such that $\Esp{\left| \we{V}_{T_1} - F \left( S_{T_1} \right) \right|^2}^{\demi} \le C/n$.

            \end{enumerate}

        \item Shall we find a similar result in a stochastic volatility model ? If the answer is yes, briefly show how to proceed, and explicit the number of liquid options which must be available.

\end{enumerate}

\section{Super-hedging with constraints on proportions of wealth}
\setcounter{equation}{0}
   We let $(\Omega, \Fc,\P)$ be a complete probability space and $\F:=(\Fc_t)_{t\le T}$ be the filtration, satisfying the usual conditions, induced by a $\P$-Brownian motion $W$. We assume that $\Fc_T=\Fc$. 

Let us consider the Black-and-Scholes one dimensional model with interest rate equal to $0$, i.e. $r\equiv 0$, in which the dynamics of the risky asset is given by 
$$
X_t=X_0e^{(\mu-\sigma^2/2) t \p \sigma W_t}\;,\;\;t\le T\;,
$$
where $W$ is a Brownian motion under $\P$, $\mu \in \R$ and  $X_0,\sigma>0$. 

\vs2

The aim of this exercise is to study the super-hedging problem under constraints on the proportion of the wealth invested in $X$. Namely, we fix $m<M$, and say that a predictable process is admissible if it takes values in $[m,M]$ $dt\x d\P$-a.e. on $[0,T]$. We denote by $\Ac$ the collection of such processes. The wealth process $Y^{y,\phi}$ associated to the initial wealth $y>0$ and the strategy $\phi \in \Ac$ has the dynamics
$$
Y^{y,\phi}_t=y \p \int_0^t \frac{\phi_s Y^{y,\phi}_s}{X_s} dX_s\;,\;t\le T\;.
$$ 

\begin{enumerate}
\item Justify (in words) the above dynamics. 
\item Show that 
$$
Y^{y,\phi}_t=y \p \int_0^t \phi_s Y^{y,\phi}_s \mu ds \p  \int_0^t \phi_s Y^{y,\phi}_s \sigma dW_s \;,\;t\le T\;.
$$
\end{enumerate}

From now on, we fix a bounded random variable $G \in L^0(\Fc_T)$ satisfying $G>0$ $\Pas$ The super-hedging price is defined as 
$$
p(G):=\inf\{y> 0~:~\exists\; \phi \in \Ac\mbox{ s.t. }  Y^{y,\phi}_T\ge G\}\;.
$$
We set 
$$
\delta(\zeta)=\zeta^\p M -\zeta^- m \mbox{ with }  \zeta^\p=\zeta\1_{\zeta>0} \mbox{ and } \zeta^-=-\zeta\1_{\zeta<0}  \mbox{ for } \zeta\in \R \;.
$$

We denote by $\Uc$ the set of predictable processes $\nu$ such that $|\nu|\le c$ $dt\x d\P$-a.e. on $[0,T]$ for some $c>0$ which depends on $\nu$. We finally define 
$$
\Ec^\nu:=e^{-\int_0^\cdot \delta(\nu_s) ds} e^{-\frac12 \int_0^\cdot |\lambda^\nu_s|^2 ds -\int_0^\cdot \lambda^\nu_s dW_s}\;,
$$
where 
$$
\lambda^\nu:=(\mu-\nu)/\sigma\;.
$$

\begin{enumerate}
\setcounter{enumi}{2}
\item Show that for any $(\phi,\nu)\in \Ac\x \Uc$, 
$$
Y^{y,\phi}\Ec^\nu=y \p \int_0^\cdot  Y^{y,\phi}_s\Ec^\nu_s\left(\phi_s \sigma- \lambda^\nu_s   \right)dW_s
\p \int_0^\cdot Y^{y,\phi}_s\Ec^\nu_s \left(\phi_s\nu_s-\delta(\nu_s)   \right) ds\;. $$ 
\item Using the definitions of $\delta$ and $\Ac$, show that $Y^{y,\phi}\Ec^\nu$ is a $\P$ super-martingale for any $(\phi,\nu)\in \Ac\x \Uc$.
\item Let $y>0$ and $\phi\in \Ac$. Show that, if  $Y^{y,\phi}_T\ge G$, then
$$
y\ge \bar p(G):=\sup_{\nu \in \Uc} \Esp{\Ec^\nu_T G}\;.
$$
\item Show that this implies that $p(G)\ge \bar p(G)$.
\end{enumerate}

We now aim at proving the converse inequality. We first assume that there exists a cadlag adapted process $P$ such that 
$$
P_t=\esssup_{\nu\in \Uc}J^\nu_t  \;\mbox{ for all } t\le T\;,
$$
where
$$
J^\nu_t:=\Esp{\Ec^\nu_T G~|~\Fc_t}/\Ec^\nu_t \mbox{ for } \nu \in \Uc\mbox{ and } t\le T\;.
$$

\begin{enumerate}
\setcounter{enumi}{6}
\item  Show that the family $\{J^\nu_t,\;\nu\in \Uc\}$ is directed upward for all $t\le T$. 
\item Show that for any $\nu^1, \nu^2\in \Uc$ and $s\le t \le T$, there exists $\nu^3 \in \Uc$ such that 
$$
\frac{\Ec^{\nu^1}_t}{\Ec^{\nu^1}_s}\frac{\Ec^{\nu^2}_T}{\Ec^{\nu^2}_t}=\frac{\Ec^{\nu^3}_T}{\Ec^{\nu^3}_s}\;.
$$
\item Deduce that $\Ec^\nu P$ is a $\P$-supermatingale for any $\nu \in \Uc$. 
\end{enumerate}

In view of the last question, and the multiplicative Doob-Meyer decomposition, it follows that,  we can find a family of martingales $\{M^\nu,\;\nu \in \Uc\}$ and a non-increasing process $\{A^\nu,\;\nu \in \Uc\}$ such that 
\be\label{eq: def exo ec nu a nu}
\Ec^\nu P =M^\nu A^\nu\;,\; A^\nu>0\;,\;M^\nu>0 \mbox{ and } A^\nu_0=1 \mbox{ for all } \nu \in \Uc\;.
\ee
In the following, we denote by $0$ a process $\nu$ such that $\nu=0$ $dt\x d\P$-a.e. on $[0,T]$.

\begin{enumerate}
\setcounter{enumi}{9}
\item  Show that $M^0_T\ge \Ec^0_TP_T\ge \Ec^0_TG>0$.  
\item Deduce that there exists a predictable process, $\Pas$  square integrable,  $\psi^0$ such that 
$$
M^0_T=M_0^0\p \int_0^T M^0_s \psi^0_s dW_s \ge  \Ec^0_T G\;.
$$
\item Show that $ M^0/\Ec^0$ can be rewritten as 
$$
M^0/\Ec^0=M^0_0\p \int_0^\cdot \frac{M^0_s}{\Ec^0_s}\left(\lambda^0_s \p \psi^0_s  \right) dW \p \int_0^\cdot \frac{M^0_s}{\Ec^0_s} \left(\lambda^0_s \psi^0_s \p |\lambda^0_s |^2\right) ds 
$$
\item Deduce that 
$$
Y^0:=M^0/\Ec^0=Y^{M_0,\phi^0} \mbox{ and } Y^{M_0,\phi^0}_T\ge G\;,
$$
for some $\Pas$  square integrable  predictable process $\phi^0$.
\item Deduce from the equality $M^\nu  =\Ec^\nu Y^0 A^0/A^\nu$ that 
$$
\int_0^\cdot F^\nu_s \left(\phi^0 \nu_s - \delta(\nu_s)\right) \frac{A^0_s}{A^\nu_s}ds\p \int_0^\cdot \frac{ F^\nu_s}{A^\nu_s} dA^0_s-\int_0^\cdot \frac{ F^\nu_sA^0_s}{|A^\nu_s|^2} dA^\nu_s=0\;\;\;
$$
where 
$$
 F^\nu:=\Ec^\nu Y^0\;
$$
for $\nu \in \Uc$.
\item By using \reff{eq: def exo ec nu a nu} and the fact that $A^0$ and $A^\nu$ are non-increasing, deduce from the previous result   that 
\b*
1\ge A^\nu&=& \int_0^\cdot A^\nu_s \left(\phi^0_s \nu_s - \delta(\nu_s)\right)ds\p \int_0^\cdot \frac{ A^\nu_s}{A^0_s} dA^0_s
\\
&\ge & \int_0^\cdot A^\nu_s \left(\phi^0_s \nu_s - \delta(\nu_s)\right)ds\p \int_0^\cdot \frac{1}{A^0_s} dA^0_s
\e*
for all $\nu \in \Uc$. 
\item Deduce from the above inequality and a formal argument that 
$$\sup_{\zeta\in \R} (\phi^0 \zeta - \delta(\zeta))<\infty\;\; \mbox{ $dt\x d\P$-a.e. }$$
\item Deduce that $\phi^0\in [m,M]$  $dt\x d\P$-a.e. 
\item Show that $\bar p(G)\ge p(G)$ and conclude. 
\end{enumerate}

 %%%%%%%%%%%%%%%%%%%%%%%%%%%%%%%%%%%%%%%%%%%%%%%%%%%%
\section{Super-hedging with impact on the volatility} 
\setcounter{equation}{0}

We let $(\Omega, \Fc,\P)$ be a complete probability space and $\F:=(\Fc_t)_{t\le T}$ be the filtration, satisfying the usual conditions, induced by a one dimensional $\P$-Brownian motion $W$. We assume that $\Fc_T=\Fc$. 

We consider a simple Black-Scholes type model in which the volatility of the risky asset is influenced by the strategy of the trader. More precisely, given a financial strategy $\phi\in \Ac$, the set of square integrable predictable processes (i.e. $\Esp{\int_0^T|\phi_s|^2ds}<\infty$), the evolution of the stock process is given by 
$$
X^\phi_{t,x}(s)=x\p \int_t^s   \sigma(\phi_r) dW_r\;,\;\;t\le s \le T\;,
$$
where $x\in \R$ is the value of the stock at time $t$, and $\sigma$ is assumed to be continuous such that 
\be\label{mideterm 2 hyp vol} 
a\in \R\mapsto (a \sigma(a),\sigma(a)) \;\mbox{ is bounded.}
\ee
Here, $\phi$ represents the number of stocks held in the portfolio and we assume that the interest rate is $0$, so that the associated wealth process starting at $y\in \R$ at time $t$ is given by 
\be\label{mid2: dyna Y}
Y_{t,x,y}^\phi(s)=y\p  \int_t^s \phi_r dX^\phi_{t,x}(r) \;,\;t\le s\le T\;.
\ee

The aim of this exercise is to study the super-hedging problem of a European option of payoff $g(X^\phi_{t,x}(T))$ paid at time $T$:
$$
v(t,x):=\inf\{y\in \R~:~ \exists\;\phi \in \Ac\mbox{ s.t. } Y_{t,x,y}^\phi(T)\ge g(X^\phi_{t,x}(T))\}.
$$
We assume that $g$ is bounded.

\begin{enumerate}
\setcounter{enumi}{0}
\item  By using \reff{mideterm 2 hyp vol}, show  that $Y_{t,x,y}^\phi$ is a martingale on $[t,T]$ for all $\phi \in \Ac$. 
\item  Assume that the infimum in the definition of $v$ is achieved and show that this implies that $v(t,x)\ge \Esp{g(X^{\hat \phi}_{t,x}(T))}=:\bar p(t,x)$ for at least one $\hat \phi\in \Ac$. 
\item Show that  there exists $\psi \in \Ac$ such that $\bar p(t,x)\p \int_t^T\psi_s dW_s= g(X^{\hat \phi}_{t,x}(T))$. 
\item Is there a chance that $\hat \phi$ and $\psi$ are such that $\psi=\hat \phi\sigma(\hat \phi)$, i.e. $Y_{t,x,\bar p(t,x)}^{\hat \phi}(T)\ge g(X^{\hat \phi}_{t,x}(T))$ ?
\end{enumerate}

From now on, we assume that $v$ is a bounded function in $ C^{1,2}([0,T)\x \R)$. We will show that $v$ should then solve
\be\label{mid2: pde}
-F\vp(t,x):=-\Lc^{\psi(x,D\vp(t,x))} \vp(t,x)=0\mbox{ on } [0,T)\x \R\;,
\ee
where, for $a\in \R$,
$
\Lc^{a} \vp(t,x)=\partial_t \vp(t,x)\p \frac12 \sigma(a)^2 D^2\vp(t,x)\;,
$
and $\psi(x,p)$ is the unique solution of 
$$
a\sigma(a)=\sigma(a) p\; \mbox{ for some } a \in \R\;,
$$
i.e. $\psi(x,p)\sigma(\psi(x,p))=\sigma(\psi(x,p)) p$. 
In the following, we shall assume that $\psi$ is Lipschitz continuous. 

\vs2

{\bf Part 1:} In this part we prove  the subsolution property. Fix $(t_0,x_0)\in [0,T)\x \R$ and assume that 
\be\label{mid2: contra sous sol}
-Fv(t_0,x_0)>0\;.
\ee
\begin{enumerate}
\setcounter{enumi}{4}
\item Show that \reff{mid2: contra sous sol} implies that $-F\vp>0$ on $B_\eps(t_0,x_0)$ for some $\eps>0$, where $\vp(t,x)=v(t,x)\p|t-t_0|^2\p|x-x_0|^4$. 
\end{enumerate}

Set $y_0=v(t_0,x_0)-(\eps\wedge \zeta)/2$ where 
$
-\zeta:=\max_{\partial B_\eps(t_0,x_0)} v-\vp<0\;.
$
Let $(X^0,Y^0)$ be the solution of 
\b*
X^0_t&=&x_0 \p \int_{t_0}^t   \sigma\left(\psi(X^0_s,D  \vp(s,X^0_s) )\right) dW_s\\
Y^0_t&=&y_0 \p \int_{t_0}^t \psi(X^0_s,D  \vp(s,X^0_s) )  \sigma\left(\psi(X^0_s,D  \vp(s,X^0_s) )\right) dW_s
 \;,\;t_0\le t\le \theta\;,
\e*
where 
$
\theta:=\inf\{s\ge t_0~:~(s,X^0_s)\notin   B_\eps(t_0,x_0) \mbox{ or } |Y^0_s-  \vp(s,X^0_s)|\ge \eps\}\;.
$
\begin{enumerate}
\setcounter{enumi}{5}
\item Show that $Y^0_\theta-v(\theta,X^0_\theta)\ge Y^0_\theta-\vp(\theta,X^0_\theta)\ge -(\eps\wedge \zeta)/2 > -\eps$. 
\item Deduce that $Y^0_\theta-  \vp(\theta,X^0_\theta)\ge  \eps>0$ if $|Y^0_\theta-  \vp(\theta,X^0_ \theta)|\ge \eps$.
\item Also deduce that $Y^0_\theta-v(\theta,X^0_\theta)\ge Y^0_\theta-  \vp(\theta,X^0_\theta)\p \zeta \ge \zeta/2>0$ if  $(\theta,X^0_\theta)\in \partial B_\eps(t_0,x_0)$.
\item Conclude that $Y^0_\theta-v(\theta,X^0_\theta)>0$. 
\item Conclude from the last assertion that \reff{mid2: contra sous sol} can not hold.
\end{enumerate}

\vs2

{\bf Part 2:} We now prove the supersolution property. Fix $(t_0,x_0)\in [0,T)\x \R$ and assume that 
\be\label{mid2: contra sur sol}
-Fv(t_0,x_0)<0\;.
\ee
We now set $\vp(t,x)=v(t,x)-|t-t_0|^2-|x-x_0|^4$ and admit that the above implies that 
\be\label{mid2: contra sur sol decal}
-\Lc^a\vp  <-\eta \mbox{ for $(t,x,a)\in B_\eps(t_0,x_0)\x \R $ s.t. } |a\sigma(a)-D\vp(t,x)\sigma(a)|\le \eps\;,
\ee
for some $\eps,\eta>0$. Let $\phi \in \Ac$ and set $(X^0,Y^0):=(X^\phi_{t_0,x_0},Y^\phi_{t_0,x_0,y_0})$ for $y_0:=v(t_0,x_0)\p (\zeta\wedge\eps)/2$ where 
$$
\zeta:=\min_{\partial B_\eps(t_0,x_0)} v-\vp>0\;.
$$
Also set 
$
\theta:=\inf\{s\ge t_0~:~(s,X^0_s)\notin   B_\eps(t_0,x_0) \mbox{ or } |Y^0_s-  \vp(s,X^0_s)|\ge \eps\}\;.
$
Given a bounded predicable  process $\lambda$, let us finally define the local martingale $L$ by 
$$
L_t:=1-\int_{t_0}^{t\wedge \theta} L_s \lambda_s\delta_s dW_s\;\mbox{ with } \delta:= \phi\sigma(\phi)-D\vp(\cdot,X^0)\sigma(\phi)\;.
$$
\begin{enumerate}
\setcounter{enumi}{10}
\item  Show that, for $t\in [t_0,\theta]$,
\b*
d\left(L_t[Y^0_t-\vp(t,X^0_t)]\right)&=& L_t \left( - \Lc^{\phi_t}\vp(t,X^0_t) - \lambda_t |\delta_t|^2  \right) dt \p \gamma_t dW_t
\e*
where $\gamma:=L\delta(1-\lambda[Y^0-\vp(\cdot,X^0)])$.
\item Deduce from \reff{mid2: contra sur sol decal} and \reff{mideterm 2 hyp vol}  that we can choose $\lambda$ such that, on $[t_0,\theta]$,  
$$
d\left(L_t[Y^0_t-\vp(t,X^0_t)]\right)\le  \gamma_tdW_t
$$
\item Deduce from the later and \reff{mideterm 2 hyp vol}  that  $L(Y^0-\vp(\cdot,X^0))$ is a  supermartingale on $[t_0,\theta]$. 
\item By the geometric dynamic programming principle, we should be able to find $\phi\in \Ac$ such that 
$$
Y^0_{ \theta}-v({ \theta},X^0_{ \theta})\ge 0\;.
$$
Show that this implies that 
$
Y^0_{ \theta}-\vp({ \theta},X^0_{  \theta})\ge (\eps \wedge \zeta) \;. 
$
\item Deduce that $L_{ \theta}\left(Y^0_{  \theta}-\vp({  \theta},X^0_{  \theta}) \right)\ge  L_{  \theta}( \eps\wedge \zeta ) $.
\item Deduce that $ (\zeta\wedge\eps)/2 \ge \Esp{L_{  \theta}\left(Y^0_{  \theta}-\vp({ \theta},X^0_{  \theta})\right)}\ge \zeta\wedge \eps$.
\item  Conclude from the last assertion that \reff{mid2: contra sur sol} can not hold.
\end{enumerate}

%%%%%%%%%%%%%%%%%%%%%%%%%%%%%%%%%%%%%%%%%%%%%%%%%%%%%%%%%%%%%%%%%%%%%%%%%%%%%%%%%%%%%%%%%%%%%%%%%%%%%%%%%%%%%

\end{document}